\documentclass[12pt,twoside,reqno]{amsproc}
\usepackage[centertags]{amsmath}
\usepackage{amssymb,amsfonts,amsthm}
\usepackage{amscd}  
\usepackage{graphicx}
\newcommand{\TryPackage}[3]{\IfFileExists{#1.sty}{\usepackage{#1}
#2}{#3}}
\TryPackage{mathrsfs}{\renewcommand{\mathcal}{\mathscr
}}{%
    \TryPackage{eucal}{}{}}
\usepackage{a4mlnewusa,math40,local,paulsmacros}

\newcommand\mcut{M^{\rm cut}}

\newcommand\tga{\tilde{\gamma}}
\newcommand\tA{\widetilde{A}}

\newcommand\etab{\widetilde{\eta}}
\newcommand\mcU{ \mathcal{U} }
\newcommand\del{\partial}

\newcommand\tsig{\tilde{\sigma}}

\newcommand\Calderon{Calder{\'o}n}
\newcommand{\specess}{\operatorname{spec}_{\operatorname{ess}}}

\newcommand{\Gr}{\operatorname{Gr}}
\newcommand{\Fred}{\operatorname{Fred}}
\newcommand{\detz}{\operatorname{det}_\zeta}
\newcommand{\detf}{\operatorname{det}_{\operatorname{F}}}

\begin{document}

\title[$\eta$--invariant, Maslov index, and spectral flow]{ The
$\eta$--invariant, Maslov index, and spectral flow for Dirac--type operators on
manifolds with boundary}

\author{Paul Kirk}
\address{Department of Mathematics\\ Indiana University\\
Bloomington, IN, 47405\\
USA}
\email{pkirk@indiana.edu}
\urladdr{http://php.indiana.edu/$\sim$pkirk}

\author{Matthias Lesch}
\address{The University of Arizona\\
Department of Mathematics\\
617 N. Santa Rita\\
Tucson, AZ, 85721--0089\\
USA}
\email{lesch@math.arizona.edu}
\urladdr{http://www.math.arizona.edu/$\sim$lesch}
\thanks{The first named author gratefully acknowledges the support of
the National Science Foundation under   grant no. DMS-9971020.  The second
named author was supported by a Heisenberg fellowship of Deutsche
Forschungsgemeinschaft and by the National Science
Foundation under grant no. DMS-0072551
}


\begin{abstract} Several proofs have been published of the $\mod \Z$ gluing
formula for the  $\eta$--invariant of a Dirac operator. However, so far the
integer contribution to the gluing formula for the $\eta$--invariant is left
obscure in the literature.  In this article we present
a gluing formula for the   $\eta$--invariant which expresses the integer
contribution  as a triple index involving the
boundary conditions and the \Calderon \ projectors of the two parts of the
decomposition. The main ingredients
of our presentation are the Scott--Wojciechowski theorem for the
determinant of a
Dirac operator on a manifold with boundary and the approach of
Br\"uning--Lesch to the $\mod \Z$ gluing formula.

Our presentation includes   careful constructions  of the Maslov index and triple
index in a symplectic Hilbert space.  As a  byproduct we give intuitively
appealing  proofs of two theorems of Nicolaescu on the spectral flow of Dirac
operators.

As an application of our methods, we carry out a detailed analysis of the
$\eta$--invariant of the odd signature operator coupled to a flat
connection using adiabatic methods. This is used to extend the definition
of  the  Atiyah--Patodi--Singer
$\rho$--invariant to manifolds with boundary. We derive a ``non--additivity''
formula for the Atiyah--Patodi--Singer
$\rho$--invariant and relate it to Wall's non-additivity formula for the
signature of even--dimensional manifolds.
\end{abstract}

\maketitle
\tableofcontents

\section{Introduction}\label{sec1}
An intriguing feature of certain spectral invariants is that they behave
nicely with respect to cutting and pasting. Such a feature has several
advantages, in particular with respect to computations. For example,  the index
of a Dirac operator behaves additively with respect to gluing of
manifolds. This is not surprising due to the locality of the index.
For higher spectral invariants (e.g. analytic torsion and the $\eta$--invariant)
cutting and pasting properties came as a surprise and proofs
are nontrivial. The gluing formula for the $\eta$--invariant has a long
history (cf. \cite{BruLes:EIN} for a historical account). Basically, there
are four different types of proof due to Bunke \cite{Bun:GPE},
Wojciechowski
\cite{Woj:AEII,Woj:AEIII}, M\"uller \cite{Mul:IDO} and Br\"uning and Lesch
\cite{BruLes:EIN}. Bunke's argument was simplified and generalized by Dai and
Freed
\cite{DaiFre:EID}.

While the articles \cite{Woj:AEII,Woj:AEIII,Mul:IDO,DaiFre:EID}
contain proofs of
the gluing formula only in $\R/\Z$, the original formula of Bunke
\cite{Bun:GPE}
offers a formula for the integer contribution in terms of indices of certain
projections. Unfortunately, these projections are not intrinsically defined
and therefore Bunke's formula is difficult to work with.
In \cite{BruLes:EIN} it is shown (though not explicitly stated)
that the integer contribution can be expressed
as the spectral flow of a naturally defined family of self--adjoint operators.

In the current paper we present another formula for the integer contribution
in terms of \Calderon\ projectors. This is very satisfactory from
a theoretical point of view since all ingredients of the formula are defined
intrinsically. Moreover, using adiabatic techniques our formula can be made
rather explicit; we carry out a detailed analysis for the odd
signature operator.

Given an appropriate orthogonal projection $P$ in the Hilbert space of sections
over the boundary, the domain of a Dirac operator $D$ can be restricted to
those sections whose restriction to the boundary lie in the kernel of $P$.
Denote the resulting operator
$D_P$.  The self--adjoint Fredholm Grassmannian $\Gr(A)$ (see
Definition \plref{SAFG})
consists of those projections $P$ so that $D_P$ is a self--adjoint discrete
Fredholm operator. It contains a distinguished element, namely the Calder\'on
projector for the Dirac operator $D$. Denote by $\etab$ the reduced
$\eta$--invariant, $\etab(D)=(\eta(D)+\dim\ker D)/2$.    Our main result
is the following. (See Theorem
\plref{symsplit}, Theorem
\plref{symsplit1}, and Lemma
\plref{splitbundles}.)

\theoremstyle{plain}
\newtheorem*{plaintheorem}{Theorem}
\begin{plaintheorem} Let $D$ be a Dirac operator on the closed manifold
$M$ and let $N\subset M$ split $M$ into $M^+$ and $M^-$. Assume
that $D$ is in   product form $D=\gamma(\tfrac{d}{dx}+A)$ in a collar of
$N$, with $A$ self--adjoint. Let
$P\in\Gr(A)$ and let $P_t$ be a smooth path in $\Gr(A)$
from $P$ to the \Calderon\ projector $P_{M^+}$ for $D$ acting on  $M^+$.   Then
\[\begin{split}
  \etab(D,M)
    &= \etab(D_{P},M^+)+\etab(D_{I-P},M^-)+ \SF(D_{P_t},M^+)_{t\in [0,1]}
+ \SF(D_{I-P_t},M^-)_{t\in [0,1]}\\
&= \etab(D_{P},M^+)+\etab(D_{I-P},M^-) -\tau_\mu(I-P_{M^-}, P, P_{M^+}).
\end{split}
\]
In particular, taking $P=P_{M^+}$,
$$\etab(D,M)=\etab(D_{P_{M^+}},M^+) +
\etab(D_{I-P_{M^+}},M^-).$$
\end{plaintheorem}

In these formulas     SF denotes the spectral flow, and $\tau_\mu$ refers to
a Maslov triple index we define for appropriate triples of
projections.  We also prove a more general formula,  Theorem
\plref{invertiblecase},  which holds for any boundary conditions $(P,Q)$, rather
than the special case $(P,I-P)$.

  It is well--known that spectral flow and
$\eta$--invariants are intimately related. It is therefore an interesting
feature of our approach that it can be used to give new and conceptually simple
proofs of  Nicolaescu's formulas for the spectral flow of a family of Dirac
operators
\cite{Nic:MIS}.
  (See Theorems \plref{ML-S6.5} and \plref{ML-S6.6}.)

For purposes of computation it is usually convenient to  use the
positive spectral projection   of the tangential operator, $P^+$,
rather than the
Calder\'on projector as boundary conditions.  According to our    theorem
   this requires computing
$\SF(D_{P_t},M^+)_{t\in [0,1]}+ \SF(D_{I-P_t},M^-)_{t\in [0,1]}$ where $P_t$ is
   a path starting at $P^+$ and ending at the Calder\'on
projector.  In favorable circumstances, such a path (actually its reverse) is
obtained by stretching the collar neighborhood of the separating hypersurface.
More precisely, replacing
$M^+$ by $M^+\cup (N\times[-r,0])$ gives a continuous path (as $r\to\infty$) of
projections starting at the Calder\'on projector and limiting essentially to
$P^+$.  This gives a method to obtain computationally useful splitting formulas,
and sheds light on the mechanism of adiabatic stretching.

We carry out this analysis in detail in Section \plref{sec8}
for  the odd signature operator. Given a flat connection   with holonomy
$\alpha$ over an odd--dimensional manifold, we take $D$ to be   the odd
signature operator   in the corresponding flat bundle. The adiabatic limit of the
Calder\'on projectors for
$D$ as  the collar is
stretched is identified in Theorem \plref{thmonadlim}. We use this
identification along with the topological invariance   of the kernel of $D$
to establish the formula (cf. \eqref{dobeedo2}):
\begin{equation*} 
\eta(D,M)=\eta(D_{P^+(V_{+,\al}) }, M^+)+
\eta(D_{P^-(V_{-,\al}) }, M^-)
+m(V_{+,\al},V_{-,\al},\alpha,g).
\end{equation*}
In this expression   $V_{\pm,\al}=\im H^*(M^\pm;\C^n_\al)\to H^*(N;
\C^n_\al)$, and
$m(V_{+,\al},V_{-,\al},\alpha,g)$ is a real--valued symplectic invariant which
depends only on the subspaces $V_{\pm,\al}\subset H^*(N;
\C^n_\al)$  and a choice of
Riemannian metric on the separating hypersurface $N$.  The projections
$P^\pm(V_{\pm,\al})$ are the sum of the  positive/negative  spectral
projections
of the tangential operator and the finite--dimensional projection to
$V_{\pm,\al}$. In particular if $H^*(N;\C^n_\al)=0$
the formula simplifies to
\[ \etab(D,M) =\etab(D_{ P^+}, M^+
)+
\etab(D_{P^-}, M^-).\]

  These formulae motivate  a definition for the
$\rho$--invariant  of a manifold  with boundary, $\rho(X,\alpha,g)$ (Definition
\plref{rhoXnotclosed}), which is shown to depend only on the smooth structure of
$X$, the conjugacy class of the representation $\alpha$, and the choice of
Riemannian metric $g$ on $\del X$. We then  prove the following theorem, and
discuss its relation  to Wall's non--additivity theorem
\cite{Wal:NAS} for the signature of even--dimensional manifolds.

{\renewcommand{\thetheorem}{\ref{there}}
\begin{theorem} Suppose the closed, odd--dimensional manifold $M$
contains a hypersurface $N$
separating
$M$ into $M^+$ and $M^-$. Fix a Riemannian metric $g$ on $N$.   Suppose
that
$\al:\pi_1(M)\to U(n)$ is a representation, and let $\tau:\pi_1(M)\to
U(n)$ denote
the trivial  representation. Then
\[
\rho (M,\al)
    =\rho(M^+,\al,g)+\rho(M^-,\al,g)
   +
 m(V_{+,\al},V_{-,\al},\alpha,g) -m(V_{+,\tau},V_{-,\tau},\alpha,g).
\]
\end{theorem}}

The paper is organized as follows:

In Section \plref{sec2} we review the basic facts about Dirac operators on
manifolds with boundary and the Grassmannian of their boundary value problems.

In Section \plref{sec3} we introduce the $\eta$--invariant and review its basic
features. Using the Scott--Wojciechowski Theorem \cite{ScoWoj:DQD} we prove
in Section \plref{sec4} a formula describing   the dependence on the choice
of boundary condition of the
$\eta$--invariant of a Dirac operator on a manifold with boundary   (Theorem
\plref{invertible}).

Section \plref{sec5} deals with splittings of manifolds. We prove a result on
the behavior of the spectral flow under splittings (Corollary \plref{ML-S5.6})
and the gluing formula for the $\eta$--invariant (Theorem
\plref{invertiblecase}).

Section \plref{sec6} contains careful constructions of  various forms of the
Maslov index for families of self--adjoint projections in a Hermitian symplectic
Hilbert space. Conventions must be set to deal with degenerate situations when
defining symplectic invariants, and we  carefully construct the various
invariants consistently and in such a way that they match our choice of
convention for the spectral flow.

A byproduct of our considerations are  new proofs of (generalizations of)
two theorems by Nicolaescu \cite{Nic:MIS} identifying the spectral flow of a
family of Dirac operators   with a Maslov index involving the \Calderon\
projectors and boundary conditions. These results (Theorem
\plref{ML-S6.5} for manifolds with boundary and Theorem \plref{ML-S6.6} for split
manifolds), together with an improvement (Theorem \plref{symsplit2})  of our
gluing formula for the
$\eta$--invariant  which allows more general boundary conditions,  are presented
in  Section
\plref{sec7}.

Finally, in Section \plref{sec8} we apply our splitting results for the
$\eta$--invariant to the special case of the odd signature operator
coupled to a flat connection. By making use of the method of adiabatic
stretching of the collar of a separating hypersurface and the fact that
the dimension of the kernels of these operators are topological, i.e.
independent of the Riemannian metric, we obtain a splitting formula
for the Atiyah--Patodi--Singer $\rho_\alpha$ invariant. The main tool introduced
in this section is Theorem \plref{thmonadlim}, which gives  a precise identification
of the adiabatic limit of the
\Calderon
\ projectors in this setting. We end the paper with an examination of the role
adiabatic stretching plays in addition formulas for the $\eta$--invariants of
general Dirac operators.

\section{Dirac operators on manifolds with boundary and the self--adjoint
Fredholm Grassmannian}
\label{sec2}

We begin by describing the set--up of Dirac operators on a manifold
with boundary.

Let $X$ denote a compact smooth Riemannian manifold with boundary  $\del X$.
We fix an identification of a neighborhood of $\del X$ in $X$ with
$\del X\times[0,\eps)$.   Let $E\to X$ be a complex Hermitian vector bundle and
suppose that
$D:C^\infty(E)\to C^\infty (E)$ is a symmetric Dirac operator,
 i.e. a symmetric first--order operator whose square is a generalized
Laplacian (the square of the leading symbol of $D$ is scalar and given
  by the metric tensor).
 The symmetry is measured with
respect to the $L^2$ inner product; thus we assume  that if
$\phi_1,\phi_2\in C^\infty(E)$ are supported in the interior of $X$ then
$$\int_X\scalar{D\phi_1}{\phi_2}_{E_x}dx=\int_X\scalar{\phi_1}{D\phi_2}_{E_x}dx.$$
A Dirac operator satisfies the unique continuation property \cite{BooWoj:EBP}.

In this paper we will deal only with the \textit{product case},
i.e. we assume that the restriction of $D$
to the collar takes the form
$D=\gamma(\frac{d}{dx}+A)$, where
   $\gamma:E_{|\del X}\to E_{|\del X}$ is a bundle endomorphism  and
$A: C^\infty(E_{|\del X})\to C^\infty(E_{|\del X})$ is a first--order self--adjoint
elliptic differential operator on the closed manifold
$\del X$ (called the {\it tangential operator})  satisfying
\begin{equation}
\gamma^2=-I, \quad\gamma^*=-\gamma,\quad \text{and}\quad \gamma A=-A\gamma.
\label{ML-G2.1}
\end{equation}
Note that $A$ is assumed to be independent of $x$ for $x\in [0,\eps)$.

The operator    $D:C^\infty(E)\to C^\infty(E)$   can be extended to an
unbounded self--adjoint operator on $L^2(E)$    by imposing appropriate
boundary conditions.  Since $D$ is a first order operator, it
can be extended to a bounded operator $H_1(E)\to L^2(E)$,
where
$H_s(E)$ denotes the Sobolev space of sections of
$E$ with $s$ derivatives in $L^2$. Given an orthogonal projection
$P:L^2(E_{|\del X})\to L^2(E_{|\del X})$  define $D_P$ to be $D$ acting on the domain
$$ \mathcal{D}(D_P):=\bigsetdef{\phi\in L^2(E)}{\phi\in
H_1(E)\;\text{and}\; P(\phi_{|\del X})=0}\subset L^2(E).$$
We will consider the operators
$D_P$   for a certain class of projections
$P$ which we now introduce.
Let $$P_{>0}:L^2(E_{|\del X})\longrightarrow L^2(E_{|\del X})$$ denote  the
positive spectral
projection for the self--adjoint tangential operator
$A:C^\infty(E_{|\del X})\to  C^\infty(E_{|\del X})$; thus if
$\{  \psi_\gl \}$ is a  basis of $L^2(E_{|\del X})$ with
$A\psi_\gl=\gl\psi_\gl$, then $P_{>0}(\sum  a_\gl
\psi_\gl)=\sum_{\gl>0}a_\gl
\psi_\gl$.

\begin{definition} \label{SAFG} Define the {\it self--adjoint Fredholm
Grassmannian}
$\Gr(A)$ to be the set of maps
$P:L^2(E_{|\del X})\to L^2(E_{|\del X})$  so that

\begin{thmenum}
\item\label{grassm1} $P$ is pseudo--differential of order $0$,
\item\label{grassm2} $P=P^*, P^2=P$, i.e. $P$ is an orthogonal projection,
\item\label{grassm3}  $\gamma P\gamma^*=I-P$,
\item\label{grassm4} $(P_{>0}, P)$ form a Fredholm pair, that is,
$$P_{>0|\im P}:\im P\to \im P_{>0}$$ is
Fredholm.
\end{thmenum}

The Grassmannian $\Gr(A)$ is topologized using the norm topology on bounded
operators.
\end{definition}

  \begin{remark}\hfill
\begin{enumerate}
\item
We note that a $P\in\Gr(A)$ also acts as a (non--orthogonal) projection
in the Sobolev space $H_s(E)$ for all $s\in\R$. This follows from
\eqref{grassm1}.

\item We obtain the same Grassmannian if we replace
$P_{>0}$ in \eqref{grassm4} by any pseudo--differential orthogonal
projection $Q$
such that $P_{>0}-Q$ is smoothing. This follows immediately from the following
general fact:
\begin{quote} Let $P,Q,R$ be orthogonal projections in the Hilbert space $H$
such that $Q-R$ is compact. Then $(P,Q)$ is a  Fredholm pair if
and only if $(P,R)$  is a Fredholm pair.
\end{quote}
This fact can be seen as follows: by \cite[Prop. 3.1]{AvrSeiSim:IPP} $(P,Q)$
is Fredholm if and only if $\pm 1\not\in\specess(P-Q)$. Since $Q-R$ is compact
this is equivalent to $\pm 1\not\in\specess(P-R)$. Applying again
\cite[Prop. 3.1]{AvrSeiSim:IPP}
the latter is the case if and only if  $(P,R)$ is Fredholm.
\end{enumerate}
\end{remark}

If $P\in \Gr(A)$, then $D_P$ is self--adjoint, Fredholm, and has compact
resolvent; in particular its spectrum is discrete and each eigenvalue has
finite multiplicity. These facts follow since $(D,P)$ is a
well--posed boundary value problem in the sense of R. T. Seeley
\cite{See:TPO}. A general reference for boundary value problems
for Dirac type operators is the monograph \cite{BooWoj:EBP}.
A different approach is presented in \cite{BruLes:STB,BruLes:BVPI}.

It will be necessary to consider a more restricted class of
projections, those that differ from $P_{>0}$ by a smoothing operator.
Define $\Gr_{\infty}(A)\subset \Gr(A)$ by
   \begin{equation}\label{grinfty} \Gr_\infty(A)=\bigsetdef{ P\in
\Gr(A)}{P-P_{>0}
\text{ is a smoothing operator}}.
\end{equation}
Again, in \eqref{grinfty} we can replace $P_{>0}$ by any
pseudo--differential orthogonal projection $Q$ such that $P_{>0}-Q$ is
smoothing.

The projection $P_{>0}$ does not lie in $\Gr(A)$ unless $\ker A=0$, since
the third condition does not hold for $P=P_{>0}$ if $\ker A\not=0$.
It is convenient to
specify  a finite rank perturbation of
$P_{>0}$ which does lie in $\Gr(A)$.

Notice that
$\gamma$ leaves $\ker A$ invariant. It is well--known that since
$(\del X,A)$ ``bounds'' $(X,D)$, the $i$ and $-i$ eigenspaces  of
$\gamma$ acting on $\ker A$ have the same dimension \cite[Chap. XVII]{Pal:SAS}.   This
implies that there are subspaces $L\subset \ker A$ satisfying
$\gamma(L)=L^\perp\cap\ker A$ (such subspaces are called {\it Lagrangian
subspaces}; see Definition \plref{defofsymp} below).
   Given a Lagrangian subspace
$L\subset \ker A$   define
\begin{equation}\label{pplus} P^+(L)=\hbox{proj}_{L}+P_{>0}.
\end{equation}
    Then
$P^+(L)$ differs from
$P_{>0}$ by the projection onto
$L$, a subspace of  $\ker A$,  which consists only of smooth
functions. Since $P_{>0}$ is a 0th order pseudo--differential projection,
so is $P^+(L)$.  It is straightforward to check that $P^+(L)\in
\Gr_\infty(A)$.

   We   call $P^+(L)$ the
   {\it Atiyah--Patodi--Singer}  projection corresponding to the
Lagrangian subspace $L$.   Notice  that $P^+(L)$ depends only
on the tangential operator $A$ and the choice of
$L$; in particular it is unchanged if $D$ is altered in the interior of
$X$.

\vskip.3in

There is a canonical  projection in
$\Gr(A)$ determined by the operator
$D$ which will play a special role in what follows, namely the  {\it
\Calderon\  projector} $P_X$. It  is defined as the orthogonal projection onto
the {\it Cauchy data space}
\begin{equation}\label{calderon}
L_X:=r\left(\ker D : H_{1/2}(E)\longrightarrow
H_{-1/2}(E)\right)\subset L^2(E_{|\del X}).
\end{equation} Here $r$ denotes the
restriction to the boundary. The trace operator $r$ is a priori only defined
on $H_s(E)$ for $s>1/2$ but one can show that $r$ defines a bounded
map from the
$H_{1/2}$--kernel of $D$ into $L^2(E_{|\del X})$ (see \cite{BooWoj:EBP}
for a proof).

The \Calderon\ projector
$P_{X}= \proj_{L_X}$ lies in
$\Gr_\infty(A)$ \cite[Prop. 2.2]{Sco:DDB}, \cite[Prop. 4.1]{Gru:TEP}.
The unique continuation property for $D$ implies that
$$r:\left(\ker D :H_{1/2}(E)\longrightarrow
H_{-1/2}(E)\right)\longrightarrow L^2(E_{|\del X})$$  is injective,
so that to any vector $x$ in the image of $P_X$ we can assign a unique
solution to $D\phi=0$ on $X$ with $\phi\in
H_{1/2}$ and $r(\phi)=x$.  This makes it possible to identify the kernel of
$D$ with  boundary condition given by a projection $P$ and the
intersection of the Cauchy data space with the kernel of $P$,
as in the following lemma.

\begin{lemma}\label{kerdp} Let $P\in\Gr(A)$. Then
$$  \ker P_{| \im P_X}=\im P_X\cap \ker P=\gamma(\ker P_X)\cap  \ker  P$$
and this space is isomorphic to  the kernel of
$D_P$.  Thus $D_P$ is invertible if and only if $\im P_X\cap \ker P=0$. In
particular
$D_{P_X}$ is invertible.
\end{lemma}
\begin{proof} If $\phi\in \ker D$, then by definition the restriction of
$\phi$ to the boundary of $X$ lies in the image of the \Calderon\  projector
$P_X$. In particular, if $\phi\in\ker D_P$, then the restriction of
$\phi$ to the boundary  lies in the intersection of $\ker P$ and the
image of $P_X$.  The unique continuation property for $D$ implies that
this intersection is exactly the kernel of
$D_P$, i.e. the kernel of $D_P$ is isomorphic to $\ker P_{|\im P_X}$.

As a discrete self--adjoint operator, $D_P$ is invertible if and only
if $\ker D_P=\{0\}$.
Moreover, $P_X$ is a self--adjoint projection satisfying the equation
$\gamma P_X
\gamma=-(I-P_X)$. Thus $\im P_X=\ker (I-P_X)=\gamma(\ker P_X)$.
   \end{proof}

In a rough sense the Atiyah--Patodi--Singer projection $P^+(L)$ and the
\Calderon\ projector $P_X$ are opposites: $P^+(L)$ is determined entirely by the
boundary data, i.e. the tangential operator $A$ acting on $\del X$ (and the
choice of $L$), whereas
$P_X$ depends on all of
$X$ and
$D$.

\vskip.3in

 For future reference we note the following special
case of a result due to K. Wojciechowski.

\begin{prop}\label{connected} The Grassmannians $\Gr(A), \Gr_\infty(A)$ are path
connected. For a fixed $P\in\Gr_\infty(A)$ (resp. $\Gr(A)$) the space
$\setdef{Q\in\Gr_\infty(A)}{\ker Q\cap\im P=0}$
(resp. $\setdef{Q\in\Gr(A)}{\ker Q\cap\im P=0}$) is path connected.
\end{prop}
\begin{remark} This result could also be proved using Proposition
\plref{homequ} below (resp. its analog for pseudo--differential Grassmannians)
and properties of the unitary group.
\end{remark}
\begin{proof} The first statement is a special case of
\cite[Appendix B]{DouWoj:ALE}, where the homotopy groups of $\Gr_\infty(A)$
and $\Gr(A)$ are computed. The path connectedness of
$\bigsetdef{Q\in\Gr_\infty(A)}{\ker Q\cap\im P=0}$ was proved
in Proposition 5.1 of \cite{ScoWoj:DQD}.
The path connectedness of
$\bigsetdef{Q\in\Gr(A)}{\ker Q\cap\im P=0}$ can be proved along the
same lines: if $\ker Q\cap \im P=0$ then $\|Q-P\|<1$ and hence
$Q_s:=Z_sPZ_s^{-1}, 0\le s\le 1$, where $Z_s:=I+s(Q-P)(2P-I)$, is a path
in $\Gr(A)$ connecting $P$ and $Q$ (cf. e.g. \cite[Sec. 3]{BruLes:BVPI}).
\end{proof}

Notice that $\Gr(A)$ and $\Gr_\infty(A)$ can also be defined by replacing
$P_{>0}$ by
$P^+(L)$ or $P_X$ in the fourth condition defining $\Gr(A)$, and in
\eqref{grinfty}.

\vskip.3in

We next discuss two alternative perspectives on the Grassmannian
$\Gr(A)$, identifying this space with the space of  certain  unitary
operators on a Hilbert space, and also with certain Lagrangian subspaces
of a symplectic Hilbert space.

The bundle endomorphism
$\gamma:E_{|\del X}\to E_{|\del X}$ induces a decomposition of
$E_{|\del X}=E_i\oplus E_{-i}$ into the $\pm i$ eigenbundles and consequently
we get a decomposition of
$L^2(E_{|\del X})$ into the $\pm i$ eigenspaces,
\begin{equation}\label{decompL2}
   L^2(E_{|\del X})= L^2(E_i)\oplus L^2(E_{-i})=:\mathcal{E}_{i}\oplus
\mathcal{E}_{-i}.
\end{equation}

Given $P\in \Gr(A)$,  write
$$P=\frac{1}{2}\begin{pmatrix} A&B\\ C&D \end{pmatrix} $$ with respect to
the decomposition \eqref{decompL2}. Then $P=P^*$ implies $C=B^*$. The
conditions $\gamma P\gamma^*=I-P$ and
$\gamma^*=-\gamma$ imply that $A=D=I$, and the condition $P^2=P$ implies that
$BB^*=I=B^*B$. This proves the first part of the following lemma.

\begin{lemma}\label{pairs} If $P\in \Gr(A)$, then with respect to the
decomposition \eqref{decompL2}, $P$ can be written in the form
   \begin{equation} P=\frac{1}{2}\begin{pmatrix} I&T^*\\ T&I
\end{pmatrix},
\end{equation} where $T$ is a 0th order pseudo--differential isometry from
$\mathcal{E}_i$ onto $\mathcal{E}_{-i}$. Conversely, given such an isometry
$T$,  then
$$\frac{1}{2}\begin{pmatrix} I&T^*\\ T&I  \end{pmatrix} $$
is a pseudo--differential projection satisfying \eqref{grassm1}, \eqref{grassm2},
\eqref{grassm3} of Definition \plref{SAFG}.

Given $$P=\frac 12 \begin{pmatrix}
I&T^*\\T&I\end{pmatrix},\;Q=\frac{1}{2}\begin{pmatrix} I&S^*\\ S&I
\end{pmatrix},$$
satisfying \eqref{grassm1}, \eqref{grassm2},
\eqref{grassm3} of Definition \plref{SAFG}, then:
\begin{thmenum}
\item $(P,Q)$ form a Fredholm pair if and only if $-1\not\in
\specess T^*S$,
\item $(P,Q)$ is invertible if and only if $-1\not\in
\spec T^*S$,
\item $\ker P\cap \im Q$ is canonically isomorphic to $\ker (I+T^*S)$,
\item $P-Q$ is smoothing if and only if $T^*S-I$ is smoothing.
\end{thmenum}
In particular, if $Q=P^+(L)$ for some Lagrangian $L\subset \ker A$, then
$P\in\Gr(A)$ if and only if
$-1\not\in\specess T^*S$.
\end{lemma}
\begin{proof} The first part was proved above. Since $S^*S=I=SS^*$, any
element in $L^2(E_{|\del X})=\mathcal{E}_{i}\oplus \mathcal{E}_{-i}$ can be
written  in the form $\begin{pmatrix}x\\ Sy\end{pmatrix} $ for $x,y\in
\mathcal{E}_i$. Since
$$Q\begin{pmatrix}x\\Sy\end{pmatrix}=\frac{1}{2}
\begin{pmatrix}x+y\\S(x+y)\end{pmatrix},$$  it follows that $ \im
Q=\bigsetdef{\begin{pmatrix}x\\Sx\end{pmatrix}}{x\in
\mathcal{E}_i}.$  Thus the restriction of $P$ to the image of
$Q$ is
$$P\begin{pmatrix}x\\Sx\end{pmatrix}=\frac{1}{2}
\begin{pmatrix}(I+T^*S)x\\T(I+T^*S)x\end{pmatrix}.$$

It follows that $(P,Q)$ is Fredholm (i.e. $P$ restricts to a Fredholm
operator on the image of $Q$) if and only if $I+T^*S$ is Fredholm, which
occurs precisely when $-1$ is not in the essential spectrum of $T^*S$.
Similarly $(P,Q)$ is invertible (i.e. the restriction of $P$ to the image
of $Q$ defines an isomorphism onto the image of $P$) if and only if $-1$ is
not in the spectrum of $T^*S$. The same argument also shows (3).

Finally, since
$$P-Q=\begin{pmatrix}0&T^*-S^*\\T-S&0\end{pmatrix},$$
$P-Q$ is smoothing if and only if $T-S$ is smoothing. Here we use that the
projections $\frac {1}{2i}(i\pm \gamma)$ onto $\mathcal{E}_{\pm i}$ are
differential operators of order $0$. Since $T,S$ are pseudo--differential and
unitary the operator $T-S$ is smoothing
if and only if $T^*S-I$ is smoothing.
\end{proof}

Let $\mathcal{U}(\mathcal{E}_i,\mathcal{E}_{-i})  $ denote the set of
unitary isomorphisms from  $\mathcal{E}_i$ to
$\mathcal{E}_{-i}$. Then $P\mapsto T$ defines a map
   \begin{equation}\label{phimap}
\Phi: \Gr(A)\longrightarrow \mathcal{U}(\mathcal{E}_i,\mathcal{E}_{-i}),
\end{equation}
i.e.,
$$P=\frac{1}{2}\begin{pmatrix} I&\Phi(P)^*\\\Phi(P)&I\end{pmatrix}.$$

More abstractly, consider   the group $\mathcal{U}$
    of unitary pseudo--differential isomorphisms
$\mathcal{E}_i\to \mathcal{E}_i$.  Let
\begin{equation}
\mcU_{\Fred}=\bigsetdef{U\in \mcU}{-1\not\in\specess U},
\label{ML-G2.7}
\end{equation}
and
\begin{equation}
\mcU_{\infty}=\bigsetdef{U\in \mcU_{\Fred}}{U-I \text{ is a smoothing
operator}}.
\label{ML-G2.8}
\end{equation}
Then given any $P\in \Gr_\infty(A)$,
the map
\begin{equation}U\mapsto \frac{1}{2}\begin{pmatrix}
I&(\Phi(P)U) ^*\\ \Phi(P)U &I
\end{pmatrix}
\label{UtoP}
\end{equation}
defines homeomorphisms
$$\mcU_{\Fred}\longrightarrow \Gr(A)$$ and
$$\mcU_{\infty}\longrightarrow \Gr_\infty(A).$$

\vskip.4in

Another useful description  of $\Gr(A)$ and $\Gr_\infty(A)$ is in terms of
Lagrangian subspaces.

\begin{lemma}\label{hermitiansymplectic}
Let $(H, \langle\ , \
\rangle)$ be a separable complex Hilbert space  and  $\gamma:H\to H$ an
isomorphism satisfying $\gamma^2=-I$,
$\gamma^*=-\gamma$.
Then there exists a subspace $L\subset H$ such that $\gamma(L)=L^\perp$
if and only if $\dim\ker(\gamma-i)=\dim\ker(\gamma+i)$.
\end{lemma}
\begin{proof} Suppose $L\subset H$ is a subspace with $\gamma(L)=L^\perp$.
Then it is easy to check that the orthogonal projections $\pi_\pm:L\to \ker(\gamma\pm i)$
are isomorphisms and hence $\dim\ker(\gamma+i)=\dim\ker(\gamma-i)$.

Conversely, if  $\dim\ker(\gamma+i)=\dim\ker(\gamma-i)$ then let
$T:\ker(\gamma-i)\to\ker(\gamma+i)$
be a unitary isomorphism. Then $L=\bigsetdef{\begin{pmatrix} x\\
Tx\end{pmatrix}}{x\in \ker(\gamma-i)}$ is a subspace satisfying
$\gamma(L)=L^\perp$.\end{proof}

\begin{dfn}\label{defofsymp} A {\it Hermitian symplectic Hilbert space} is a
separable complex Hilbert space together with an isomorphism $\gamma:H\to H$
satisfying $\gamma^2=-I, \gamma^*=-\gamma$ and such that the
$i$ and $-i$ eigenspaces of $\gamma$ have the
same dimension (i.e. if $H$ is infinite--dimensional we require that both
eigenspaces are infinite--dimensional). The symplectic form is the skew-Hermitian
form
\begin{equation}\label{symplectic}
\omega(x,y):=\scalar{x}{\gamma y}.
\end{equation}

A {\em Lagrangian subspace} $L\subset H$  of a Hermitian
symplectic Hilbert space is a subspace so that $\gamma(L)=L^\perp$.
A Lagrangian subspace is automatically closed.
\end{dfn}

The space $L^2(E_{|\del X})$ together with the map $\gamma$ is a Hermitian
symplectic Hilbert space.   The space $\ker A$ is a finite--dimensional
Hermitian symplectic Hilbert space since $(\del X,A)$ bounds $(X,D)$.

  Given
$P\in \Gr(A) $, the kernel of $P$ is a Lagrangian subspace, since $\ker P$ is
orthogonal to $\gamma(\ker P)$.
Notice that the kernel of $P$ can be
expressed as the graph of
$-\Phi(P)$,
$$\ker P=\bigsetdef{
\begin{pmatrix}x\\-\Phi(P)x\end{pmatrix}}{x\in\mathcal{E}_i}\subset
L^2(E_{|\del X}).$$

This gives a third characterization of
$\Gr(A)$ as follows.
   We  define $\mathcal{L}$ to be the set of Lagrangian subspaces of
$L^2(E_{|\del X})$ whose associated projections are pseudo--differential of
order $0$. The Cauchy data space, $L_X$, (the image  of the \Calderon\
projector)   is a Lagrangian subspace of $L^2(E_{|\del X})$.

Define
\begin{equation}
\mathcal{L}_{\Fred}=\bigsetdef{L\in \mathcal{L}}{(L,\gamma(L_X))
\text{ is a Fredholm pair of subspaces}},
\end{equation}
and
\begin{equation}
\mathcal{L}_{\infty}=\bigsetdef{L\in \mathcal{L}_{\Fred}}{\proj_L-\proj_{L_X}
\text{ is a smoothing operator }}.
\end{equation}
Then we have homeomorphisms
$$\mathcal{L}_{\Fred}\longrightarrow \Gr(A)$$ and
$$\mathcal{L}_{\infty}\longrightarrow \Gr_\infty(A).$$

The identifications of $\mathcal{L}$, $\Gr(A)$, and
$\mathcal{U}$ are determined by the conditions that
$$L\in \mathcal{L}_{\Fred} , \ P\in \Gr(A), \text{ and }  T\in
\mcU_{\Fred}(\mathcal{E}_i,
\mathcal{E}_{-i})$$ correspond if
$$L=\im P=\hbox{graph of } T  \text{ and }
    T= \Phi(P).$$

\section{The $\eta$--invariant and spectral flow}
\label{sec3}

It was mentioned  in the last section that $D_P$ is the self--adjoint
realization of a well--posed boundary value problem and hence it
is a discrete operator in the Hilbert space $L^2(E)$. For the discussion
of $\zeta$-- and  $\eta$--functions we need the more refined analysis of
the heat trace of $D_P$. The $\zeta$-- and $\eta$--functions of $D_P$ are
defined, for $\Re(s)>>0$, by
\begin{equation}
   \begin{split}
     \eta(D_P;s)&:=\tr(D_P |D_P|^{-s-1})=\sum_{\gl\in\spec D_P\setminus\{0\}}
     \sign(\gl) |\gl|^{-s},\\
     \zeta(D_P;s)&:=\tr(D_P^{-s})=\sum_{\gl\in\spec D_P\setminus
\{0\}} \gl^{-s}\\
     &= \frac 12 \bigl(\zeta(D_P^2;s/2)+\eta(D_P;s)\bigr)+e^{-i\pi
s}\frac 12\bigl(\zeta(D_P^2;s/2)-\eta(D_P;s)\bigr).
   \end{split}\label{ML-G3.1}
\end{equation}

\begin{theorem}\label{ML-S3.1} For $P\in\Gr(A)$ the functions
$\zeta(D_P;s),\eta(D_P;s)$
extend meromorphically to the whole complex plane with poles of order
at most $2$.
If $P\in\Gr_\infty(A)$ then $\eta(D_P;s)$ and $\zeta(D_P;s)$ are regular
at $s=0$. Moreover $\zeta(D_P;0)$ is independent of $P\in\Gr_\infty(A)$.
\end{theorem}
That the $\zeta$-- and $\eta$--functions extend meromorphically has been proved
in increasing generality in \cite{DouWoj:ALE},  \cite{GruSee:WPP},
\cite{GruSee:ZEF},
\cite{BruLes:EIN}, \cite{Woj:DAI}, and \cite{Gru:TEP}. The definitive
treatment of all well--posed boundary value problems is given
in \cite{Gru:TEP}. The proof of the statement about regularity at $s=0$
can be found in \cite{Woj:DAI}. The methods of \cite{Gru:TEP} show that
the assumption $P\in\Gr_\infty(A)$ can be somewhat relaxed \cite{Gru:PZE}.
Finally, that $\zeta(D_P;0)$ is independent of $P\in\Gr_\infty(A)$ is
\cite[Prop. 0.5]{Woj:DAI}.

\begin{definition}\label{ML-S3.2} The $\eta$-{\it invariant of} $D_P$,
$\eta(D_P)$, is defined to be the constant term in the Laurent
expansion of
$\eta(D_P;s)$ at $s=0$, i.e.
$$\eta(D_P;s)= a s^{-2} + b s^{-1} + \eta(D_P) + O(s).$$
\end{definition}

We also give a symbol to a convenient normalization of the
$\eta$--invariant.

\begin{definition}\label{ML-S3.3} The {\it reduced $\eta$--invariant}
is defined to be
   \begin{equation}\label{xidef}
\etab(D_P)=\left(\eta(D_P) + \dim \ker D_P\right)/2.
\end{equation}

\end{definition}

We continue  with a discussion of the spectral flow and its
relation to the $\eta$--invariant.
Suppose one is given  a smooth path
of  Dirac type operators  $D_t:C^\infty(E)\to C^\infty(E) , t\in [0,1]$,
over $X$  so that $D_t=\gamma(\frac{d}{dx}+A_t)$ on the collar.  Choose
a smooth path of projections
$P_t$ so that $P_t\in \Gr(A_t)$ for $t\in [0,1]$. Then the family
$D_P(t):=(D_t)_{P_t}$ is in particular a graph continuous family  of
self--adjoint discrete operators. As a consequence, the eigenvalues
of $D_{P_t}$ vary continuously (as a general reference see
\cite{Kat:PTL}). The {\it spectral flow} of the family
$D_P(t)$, which we denote by  $\SF(D_P(t))_{t\in[0,1]} $ or just
$\SF(D_P(t))$, is  the integer defined  (roughly) to be the difference in
the number of eigenvalues that start negative and end non-negative and
the number of eigenvalues that  start non-negative and end negative
(see \cite{BooWoj:EBP, CapLeeMil:SA2} for a
precise definition).
Notice that we have chosen a particular convention for dealing with zero
eigenvalues.  This convention is often called the $(-\eps,-\eps)$ spectral
flow in the literature, since it corresponds to intersecting the graph of the
spectrum as a function of $t$ with the line  from $(0,-\eps)$ to  $(1,-\eps)$.

The $1$--parameter family of $\eta$--invariants $\eta(D_P(t))\in\R$ will in
general not vary smoothly with respect to $t\in [0,1]$. However,
it follows from the work of G. Grubb \cite{Gru:TEP} that the reduction modulo
integers of the reduced $\eta$--invariant $\tilde\eta(D_P(t))$
varies smoothly with $t$.
In particular, the real valued function
$$u\mapsto \int_0^u\frac{d}{dt}\bigl(\eta(D_P(t))\bigr) dt
   $$ is smooth.

   In general, given a smooth function $f:[0,1]\to \R/\Z=S^1$, the
expression $u\mapsto c+\int_0^u \frac{df}{dt}dt$ is just an explicit
formula for the unique smooth lift of $f$ to the universal cover $\R$ of
$S^1$ starting at $c\in\R$. Thus if $f$ and $g$ are (possibly
discontinuous) functions from
$[0,1]$ to
$\R$ so that the reductions of $f$ and $g$ modulo $\Z$ are smooth and
agree,  then the smooth real--valued functions
$u\mapsto \int_0^u\frac{df}{dt}dt$ and
$u\mapsto\int_0^u\frac{dg}{dt}dt$  coincide.

\begin{lemma}\label{sf-eta} Suppose that $D_t, t\in [0,1]$, is a smooth
path of symmetric Dirac type operators as above, and
$P_t\in \Gr(A_t)$ is a smooth path, giving a smooth path of self--adjoint
discrete operators $D_P(t)$.

Then
$$\etab(D_P(1))- \etab(D_P(0))= \SF(D_P(t))_{t\in[0,1]}+
\frac{1}{2}\int_0^1\frac{d}{dt}(\eta (D_P(t))) dt . $$  Moreover, if the
dimension of the kernel of $D_P(t)$ is independent of $t$, then the
function $t\mapsto \etab(D_P(t))$ is smooth.
\end{lemma}
\begin{proof} We only sketch the proof, since this fact is well--known,
at least when the $\eta$--function is regular at $s=0$, and the general
case is proven by the same argument, because the pole  of the $\eta$--function
at $s=0$ is determined by the  asymptotics of the spectrum,
whereas the spectral flow depends only on the small eigenvalues.

   Given
$r\in [0,1]$,  choose an $\eps>0$ so that  $\pm\eps$ does not lie in the
spectrum of $D_P(r)$. Applying standard results from perturbation theory
\cite{Kat:PTL} we infer that $\pm\eps$ does not lie in the spectrum of
$D_P(t)$ for $t$ close enough to  $r$, say
$t\in [t_0,t_1]$. Moreover the span of those eigenvectors of
$D_P(t)$ whose eigenvalues lie in $(-\eps,\eps)$ varies continuously for
$t\in[t_0,t_1]$.

Thus we can write $\eta(D_P(t);s)$ for $t\in[t_0,t_1]$ and
$\Re(s)>>0$ as a sum
\begin{equation}\label{etasum}
\begin{split}
\eta(D_P(t);s)&=\sum_{\gl\in\spec D_P(t), 0<|\gl|<\eps}\sign(\gl)|\gl|^{-s} +
\sum_{\gl\in\spec D_P(t),|\gl|>\eps }\sign(\gl)|\gl|^{-s}\\
     &=\eta_{<\eps}(D_P(t);s)+\eta_{>\eps}(D_P(t);s).
\end{split}
\end{equation}

The sum $\eta_{<\eps}(D_P(t);s)$ is finite, and so its analytic
continuation to $s=0$ is integer valued:
$$\eta_{<\eps}(D_P(t);0)=\sum_{\gl\in\spec D_P(t),0<|\gl|<\eps}\sign(\gl). $$
Thus
\begin{equation}\label{eq1}
\begin{split}
         \eta_{<\eps}&(D_P(t_1);0)-\eta_{<\eps}(D_P(t_0);0)\\
                                  &= 2\SF(D_P(t))_{t\in [t_0, t_1]} +
\dim\ker D_P(t_0) -\dim
\ker D_P(t_1).
\end{split}
\end{equation}
Notice that this equation depends on our
choice of convention for defining the spectral flow.
The function $\eta_{>\eps}(D_P(t);s)=\eta(D_P(t);s)-\eta_{<\eps}(D_P(t);s)$
varies smoothly in $t\in [t_0,t_1]$ since we
have subtracted the eigenvalues that cross zero, and since no eigenvalues
equal $\pm \eps$ in this interval.  If we define $\eta_{>\eps}(D_P(t))$
similarly
to Definition \plref{ML-S3.2} then $\eta_{>\eps}(D_P(t))$ is smooth and
$$\eta_{>\eps}(D_P(t))\equiv \eta(D_P(t)) \mod\Z.$$

Therefore, using \eqref{eq1} we obtain
\begin{equation}\label{eq2}
    \begin{split}
      \int_{t_0}^{t_1}\frac{d}{dt}&\bigl(\eta(D_P(t))\bigr) dt\\
 &=\int_{t_0}^{t_1}\frac{d}{dt}\bigl(\eta_{>\eps}(D_P(t))\bigr) dt\\
                     &=\eta_{>\eps}(D_P(t_1))-\eta_{>\eps}(D_P(t_0))\\
                     &=\eta(D_P(t_1))- \eta_{<\eps}(D_P(t_1))-\eta(D_P(t_0))+
                        \eta_{<\eps}(D_P(t_0))\\
                     &= \eta( D_P(t_1)) -\eta(D_P(t_0))-
                         2\SF(D_P(t))_{t\in[t_0,t_1]}\\
                     &\quad+\dim\ker D_P(t_1) -\dim\ker D_P(t_0).
    \end{split}
\end{equation}
Dividing by $2$  proves the lemma over the interval $[t_0, t_1]$. The
general case is obtained by covering the interval
$[0,1]$ by small subintervals and adding the  results.

For the last assertion, notice that if the dimension of $\ker D_P(t)$ is
independent of $t$, then $\SF(D_P(t))_{t\in [0,s]}=0$ for all $s\in [0,1]$.
\end{proof}

\section{The Scott--Wojciechowski theorem}
\label{sec4}

The theorem of Scott and Wojciechowski \cite{ScoWoj:DQD} identifies the
regularized $\zeta$--de\-ter\-minant     of a
boundary--value problem for a Dirac operator with a Fredholm determinant
of the associated boundary projection. In this section we summarize and
slightly extend that part of their result which we need, in the language
of this article. Briefly, their theorem shows that the reduction mod $\Z$ of the
$\eta$--invariant of $D_P$ for $P\in \Gr_\infty(A)$  and the Fredholm
determinant of the unitary map
$\Phi(P)$ which corresponds to $P$ via \eqref{UtoP} agree up to a
constant independent of $P$.  The important consequence for this
article is that
the  $\mod \Z$ reduction of the $\eta$--invariant for a manifold with boundary
depends only on the boundary data and the \Calderon\ projector.

In this section $D$ denotes a fixed Dirac type operator on a manifold $X$
with boundary $\del X$ and $A$ denotes its tangential operator.

Before stating the Scott--Wojciechowski theorem, let us briefly recall
the $\zeta$--regularized determinant. Let $P\in\Gr_{\infty}(A)$. Then
$\zeta(D_P;s)$ is regular at $s=0$ and one puts
\begin{equation}
     \detz D_P:=\begin{cases}\exp\bigr(-\zeta'(D_P;0)\bigr), &
0\not\in\spec D_P,\\
                            0, & 0\in\spec D_P.
               \end{cases}
\label{ML-G4.1}
\end{equation}
In view of \eqref{ML-G3.1} and Theorem \plref{ML-S3.1} a straightforward
calculation shows for $D_P$ invertible
\begin{equation}
    \detz D_P= \exp\Bigl( i\frac{\pi}{2} \bigl(\zeta(D_P^2;0)-\eta(D_P)\bigr)
                  -\frac 12 \zeta'(D_P^2;0)\Bigr).
    \label{ML-G4.2}
\end{equation}
We emphasize that the regularity of $\eta(D_P;s)$ and $\zeta(D_P;s)$ at
$s=0$ is essential for \eqref{ML-G4.2} to hold.
\eqref{ML-G4.2} implies that in general $(\detz D)^2\not= \detz
(D^2).$ Note that Fredholm determinants are multiplicative, i.e. if
$S,T$ are operators of determinant class in a Hilbert space then
$\detf(ST)=\detf(S)\detf(T)$, where $\detf$ denotes the Fredholm determinant.

With these preparations, the Scott--Wojciechowski theorem
reads as follows.
\begin{theorem}\label{ML-S4.1} Let $P\in\Gr_\infty(A)$. Then
\begin{equation}    \detz(D_P)=\detz(D_{P_X})
       \detf\bigl(\frac{I+\Phi(P_X)\Phi(P)^*}{2}\bigr).
       \label{ML-G4.3}
\end{equation}
\end{theorem}
This result was proved for $ M$ odd--dimensional in \cite[Thms. 0.1,
1.4]{ScoWoj:DQD}. An alternative proof which applies to all dimensions and to
slightly more general operators will be presented in \cite{Les:IP}.

In view of \eqref{ML-G4.2} the Scott--Wojciechowski theorem can
be applied to express the dependence of $\tilde\eta(D_P)$ on $P$
in terms of Fredholm determinants.

Let $P, Q\in \Gr_\infty(A)$. Since $\Phi(Q)\Phi(P)^*-I$ is a smoothing
operator, it is of trace class and hence
\begin{equation}
     \frac{I+\Phi(Q)\Phi(P)^*}{2}=I+\frac{\Phi(Q)\Phi(P)^*-I}{2}
\end{equation}
is of determinant class. In particular, $\frac{I+\Phi(P_X)\Phi(P)^*}{2}$
is of determinant class and thus the right hand side in \eqref{ML-G4.3}
is well--defined.

Also, $\Phi(P)\Phi(Q)^*$ is of determinant class. Hence the
determinant
$\detf(\Phi(P)\Phi(Q)^*)$ is defined and lies in $ U(1)$ since
$\Phi(P)\Phi(Q)^*$ is unitary.

\begin{theorem}\label{ML-S4.2}  Let $P,Q\in \Gr_\infty(A)$. Then
\begin{equation} e^{2\pi
i(\etab(D_P)-\etab(D_{Q}))}=\detf(\Phi(P)\Phi(Q)^*).\label{ML-G4.5}
\end{equation}
\end{theorem}
\begin{proof} Assume first that $P$ is the
\Calderon\ projector $P_X$ and that the pair $(P_X,Q)$ is invertible. By Lemma
\plref{kerdp} this means that $D_{P_X}$ and $D_Q$ are invertible.
Putting Theorem \plref{ML-S4.1}
and \eqref{ML-G4.2} together and taking into account that $\zeta(D_P^2;0)$ is
independent of $P$ (Theorem \plref{ML-S3.1}), we obtain
\begin{equation}
    e^{i\frac{\pi}{2}\bigl(\eta(D_{P_X})-\eta(D_Q)\bigr)} e^{\frac
12\bigl(\zeta'(D_{P_X}^2;0)-\zeta'(D_Q^2;0)\bigr)}
=\frac{\detz(D_Q)}{\detz(D_{P_X})}=
\detf\bigl(\frac{I+\Phi(P_X)\Phi(Q)^*}{2}\bigr),
    \label{ML-G4.6}
\end{equation}
and thus
\begin{equation}
\frac{\detf\bigl(\frac{I+\Phi(P_X)\Phi(Q)^*}{2}\bigr)}{\bigl|\detf\bigl(\frac{I+\Phi(P_X)\Phi(Q)^*}{2}\bigr)\bigr|}
    = e^{i\frac{\pi}{2}\bigl(\eta(D_{P_X})-\eta(D_Q)\bigr)}.
\label{ML-G4.7}
\end{equation}

Since $\Phi(P_X)\Phi(Q)^*-I$ is of trace class we may choose a self--adjoint
trace class operator $H$ such that
$e^{iH}=\Phi(P_X)\Phi(Q)^*$. Then
\begin{equation}\begin{split}
     \detf\bigl(\frac{I+\Phi(P_X)\Phi(Q)^*}{2}\bigr)^2&=
     \detf\bigl(\frac{I+e^{iH}}{2}\bigr)^2\\
        &= \detf\bigl(e^{iH}\cosh^2(H/2)\bigr)\\
        &= \detf\bigl(\Phi(P_X)\Phi(Q)^*\bigr)\detf\bigl(\cosh^2(H/2)\bigr),
        \end{split}
\label{ML-G4.8}
\end{equation}
where we have used the multiplicativity of the Fredholm determinant in the last
line. Consequently
\begin{equation}
  \frac{\detf\bigl(\frac{I+\Phi(P_X)\Phi(Q)^*}{2}\bigr)^2}{\bigl|\detf\bigl(\frac{I+\Phi(P_X)\Phi(Q)^*}{2}\bigr)\bigr|^2}
       = \detf\bigl(\Phi(P_X)\Phi(Q)^*\bigr).
\label{ML-G4.9}
\end{equation}
Putting together
\eqref{ML-G4.7} and \eqref{ML-G4.9} we obtain \eqref{ML-G4.5} for $P=P_X$ and
$Q\in\Gr_\infty(A)$ such that $(P_X,Q)$ is an invertible pair. However, since
both sides of \eqref{ML-G4.5} depend continuously on $Q$,
\eqref{ML-G4.5} remains valid
for all $Q\in\Gr_\infty(A)$. Finally, if $P,Q\in\Gr_\infty(A)$ are arbitrary
then
\begin{equation}\begin{split}
     e^{2\pi i(\tilde \eta(D_P)-\tilde\eta(D_Q))}&=e^{2\pi
i(\tilde\eta(D_P)-\tilde\eta(D_{P_X}))}e^{2\pi
     i(\tilde\eta(D_{P_X})-\tilde\eta(D_Q))}\\
      &=
     \detf\bigl(\Phi(P)\Phi(P_X)^*\bigr)\detf\bigl(\Phi(P_X)\Phi(Q)^*\bigr)\\
      &= \detf\bigl(\Phi(P)\Phi(Q)^*\bigr).
        \end{split}
\end{equation}
\end{proof}

We will use the following convenient form of the Scott--Wojciechowski
theorem. We consider the reals $\R$ as the universal cover of $U(1)$ via
the map $r\mapsto e^{2\pi i r}$.

\begin{cor}\label{lifting}  Let $P_t, \  t\in [0,1]$, be a smooth
path in
$\Gr_\infty(A)$. Then the map
$$s\mapsto \frac{1}{2}\int_0^s \frac{d}{dt}(\eta(D_{P_t})) dt$$ is the
unique lift to $\R$ of the  map $[0,1]\to U(1)$
$$s\mapsto \detf\bigl(\Phi(P_s)\Phi(P_0)^*\bigr).$$
\end{cor}

   In preparation for the next theorem, suppose that $P\in \Gr_\infty(A)$.
  From Lemma \plref{kerdp} we know that  $D_P$ is invertible if and only if
$\ker P_X\cap \gamma (\ker P)=0$ where $P_X$ denotes the \Calderon\
projector; by Lemma
\plref{pairs} this happens if and only if
$-1\not\in\spec(\Phi(P)\Phi(P_X)^*)$. In fact,
the kernel of $D_P$ is canonically isomorphic to $\ker(I+\Phi(P)\Phi(P_X)^*)$.

Using the functional calculus we can define the operator
$\log(\Phi(P)\Phi(P_X)^*)$. The choice of the branch of $\log$
will be essential in what follows. We define $\log:\C\setminus\{0\}\to \C$
as follows
\begin{equation}
        \log(re^{it})=\ln r+i t, \quad r>0, -\pi<t\le \pi.
\label{ML-G4.11}
\end{equation}
Since
$-1\not\in\specess(\Phi(P)\Phi(P_X)^*)$, $-1$ is an isolated point in the
spectrum
of $\Phi(P)\Phi(P_X)^*$ and thus we can choose a holomorphic branch of the
logarithm which coincides on $\spec(\Phi(P)\Phi(P_X))^*$ with $\log$ defined in
\eqref{ML-G4.11}.
The so defined
$\log(\Phi(P)\Phi(P_X)^*)$ is of trace class and
$$\tr
\log(\Phi(P)\Phi(P_X)^*)\equiv\log\detf\bigl(\Phi(P)\Phi(P_X)^*\bigr)\mod
2\pi i\Z.$$

After these preparations we can improve Theorem \plref{ML-S4.2} as follows.

\begin{theorem}\label{invertible} Let $X$ be a compact 
manifold with boundary   and let $D$ be a Dirac type operator such
that in a collar
$\del X\times [0,\eps)$ of the boundary $D$ takes the form
$D=\gamma(\frac{d}{dx}+A)$ with $A, \gamma$ as in \eqref{ML-G2.1}. Let $\Phi$
be the map defined in \eqref{phimap}. Then for
$P\in \Gr_\infty(A)$ we have
$$\etab(D_P)-\etab(D_{P_X})=\frac{1}{2\pi i} \tr\log
(\Phi(P)\Phi(P_X)^*).$$
\end{theorem}
\begin{proof} We assume first that $D_P$ is invertible. $D_{P_X}$ is
invertible by Lemma \plref{kerdp}. In view of Proposition \plref{connected}
and Lemma \plref{kerdp}
the space of those
$P\in \Gr_\infty(A)$ so that $D_P$ is invertible is path connected.
Choose a smooth path $P_t$ in $\Gr_\infty(A)$ starting at $P_X$ and ending
at $P$ so that
$D_{P_t}$ is invertible for all $t$.

    The spectral flow of $D_{P_t}$ equals zero since the kernel is zero
along the path and so Lemma
\plref{sf-eta} shows that
   $t\mapsto \etab(D_{P_t})$  is smooth.  Hence
   \begin{equation}\label{smooth1} t\mapsto \etab(D_{P_t})-\etab(D_{P_X})
\end{equation}
   is smooth.  Also, the map
\begin{equation}\label{smooth2} t\mapsto
\frac{1}{2\pi i}\tr\log(\Phi(P_t)\Phi(P_X)^*)
\end{equation}
   is smooth since
$-1\not\in\spec(\Phi(P_t)\Phi(P_X)^*)$ for all $t$ and hence $\log$ is
holomorphic on $\spec(\Phi(P_t)\Phi(P_X)^*)$.

Theorem \plref{ML-S4.2} states that the  two smooth   functions of
\eqref{smooth1} and \eqref{smooth2} are the lifts to $\R$ of the same
function to
$U(1)=\R/\Z$, and they both start at $0$.  Hence they coincide for all
$t$.

Now let $P\in\Gr_\infty(A)$ be arbitrary. We may choose a path $(P_t)_{-\eps\le
t\le \eps}$ in $\Gr_\infty(A)$ such that $(P_t,P_X)$ is invertible
for $t\not=0$,
$P_0=P$,
and such that at $t=0$  exactly $k=\dim\ker D_{P_0}$ eigenvalues of
$\Phi(P_t)\Phi(P_X)^*$
cross $-1$ from the upper half plane to the lower half plane and no
eigenvalues cross from the lower half plane to the upper half plane.
To see this let $R$ be the orthogonal projection onto
$\ker(I+\Phi(P)\Phi(P_X)^*)$. The projection  $R$ is a pseudo--differential
operator. Now put
\begin{equation}
\Phi(P_t):=\bigl(e^{i(\pi+t)}R\oplus
(I-R)\Phi(P)\Phi(P_X)^*\bigr)\Phi(P_X).
\label{ML-G4.14}
\end{equation}
By our choice of $\log$ we then have
\begin{equation}
      \frac{1}{2\pi i}\tr\log(\Phi(P_0)\Phi(P_X)^*)=\lim_{t\to
      0-}\frac{1}{2\pi i}\tr\log(\Phi(P_t)\Phi(P_X)^*).
\end{equation}
Moreover, for $t\not=0$ we have from the first part of this proof
\begin{equation}
     \tilde\eta(D_{P_t})-\tilde\eta(D_{P_X})=\frac{1}{2\pi
     i}\tr\log(\Phi(P_t)\Phi(P_X)^*).
\label{ML-G4.16}
\end{equation}
 From Lemma \plref{sf-eta}, \eqref{ML-G4.14} and \eqref{ML-G4.16} one
infers $\SF(D_{P_t})_{-\eps\le t\le\eps}=-k$ and since
$\dim\ker D_{P_0}=k$ at $t=0$ exactly $k$ eigenvalues of $D_{P_t}$
cross $0$ from $+$ to $-$ and no eigenvalues cross from $-$ to $+$.
Hence
\begin{equation}\begin{split}
   \tilde\eta(D_{P_0})-\tilde\eta(D_{P_X})
&=\lim_{t\to 0-}\tilde\eta(D_{P_t})-\tilde\eta(D_{P_X})\\
&=\lim_{t\to 0-} \frac{1}{2\pi i}\tr\log(\Phi(P_t)\Phi(P_X)^*)\\
        &= \frac{1}{2\pi i}\tr\log(\Phi(P_0)\Phi(P_X)^*),
        \end{split}
\end{equation}
completing the proof.
\end{proof}

\section{Splittings of manifolds and the $\eta$--invariant I}
\label{sec5}

We consider now the gluing problem for the $\eta$--invariant. Suppose we
are given a closed manifold $M$ containing a separating hypersurface
$N\subset M$. We consider only Dirac operators $D$ on $M$ so that in a
collar neighborhood $[-\eps,\eps]\times N$ of $N$,
$D$ has the form $D=\gamma(\frac{d}{dx}+A)$ as in \eqref{ML-G2.1}.

Let $\mcut$ denote the compact manifold with boundary obtained by cutting
$M$ along $N$. Thus $\mcut$ is the disjoint union of two submanifolds
$M^+$ and $M^-$, with $\partial M^+$ and $\partial M^-$ canonically
identified with $N$.
To apply the results of the previous section, we reparameterize the
collar of $M^-$ as $\partial M^-\times [0,\eps]$ with $x=0$ corresponding to
the boundary. See the following figure.

\vskip2ex

 \begin{center}\includegraphics[width=4in]{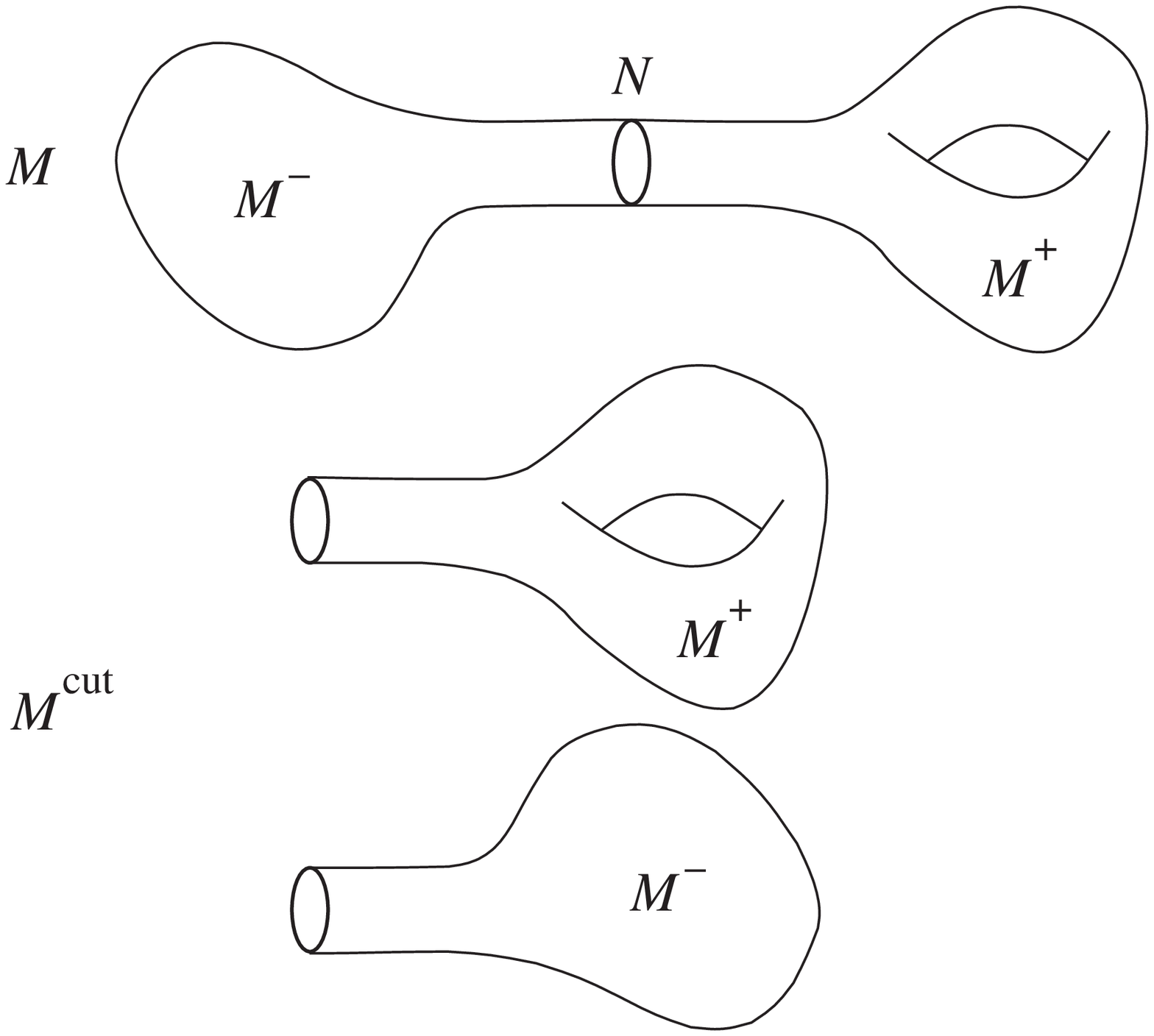}
 \end{center}
\vskip2ex
\centerline{{\bf The manifolds} {$M$} {\bf and }{$\mcut$}}
\vskip2ex

The $H_1$--Sections of a bundle $E$ over $M$ correspond to sections
$f\in H_1(E_{|{\mcut}})$ over
$\mcut$ so that $f_{|\del M^+}=f_{|\del M^-}$ with respect to the
canonical identifications $\del M^\pm =N$.
More precisely, the restriction of the section $f$ to the boundary of
$\mcut$ lies in $H_{1/2}(E_{|\del \mcut})\subset L^2(E_{|\del \mcut})$.
The identification of $\del
\mcut$ with two copies of $N$ gives a canonical decomposition
   \begin{equation}\label{decompl2} L^2(E_{|\del
\mcut})=L^2(E_{|N})\oplus L^2(E_{|N})
\end{equation} where the first factor corresponds to $\del M^+$ and the
second to $\del M^-$.  The restriction of a section $f$ over $ \mcut$  to
the boundary can thus be written as
$(f_+,f_-)$, and the sections over $M$ correspond exactly to those $f$ so
that $f_+=f_-$.

On the collar of $\mcut$,  The operator $D$  takes the form
\begin{equation}D=\begin{pmatrix}\gamma&0\\0&-\gamma\end{pmatrix}
\left(\frac{d}{dx}+\begin{pmatrix}A&0\\0&-A\end{pmatrix}\right) =:\tga
(\frac{d}{dx}+\tA) \label{ML-G5.2}
\end{equation}
with respect to the decomposition
\eqref{decompl2} (this is because of the change of parameterization of the
collar of
$M^-$).

Note that the (closed) diagonal subspace
$$\Delta=\{(f,f)\ | \ f\in L^2(E_{|N})\}\subset L^2(E_{|N})\oplus
L^2(E_{|N})$$ is Lagrangian. In fact:
\begin{enumerate}
\item $\Delta$ is orthogonal to $\tga(\Delta)$ since
$$\langle (f,f),\tga(g,g)\rangle=\langle f,\gamma g \rangle +\langle f,-\gamma
g \rangle=0.$$
\item $\Delta+\tga(\Delta)=L^2(E_{|N})\oplus L^2(E_{|N})$ since
$$(f,g)=\frac{1}{2}\left( (f+g,f+g)+\tga(-\gamma f+\gamma g ,-\gamma f+\gamma
g)\right).$$
\end{enumerate}

The orthogonal projection to $\Delta^\perp$ will be denoted by
$P_\Delta$. It is called the {\it continuous transmission}  projection.
By construction an $H_1$--section $f$ over $\mcut$ defines an
$H_1$--section over $M$ if
and only if the restriction $(f_+,f_-)$ of $f$ to the boundary satisfies
$P_\Delta(f_+,f_-)=0$.   Note that
$P_\Delta\in \Gr(\tA)$ since   $D_{P_\Delta}$ is canonically identified
with the (Fredholm) operator $D$ acting on the closed manifold $M$.

With respect to the decomposition \eqref{decompl2} the operator
$P_\Delta$ takes the form
\begin{equation}
P_\Delta=\frac{1}{2}\begin{pmatrix}{\hskip.13in 1}&-1\\-1&{\hskip.13in
1}\end{pmatrix}.\label{ML-G5.2.5}
\end{equation}

In a strict sense, the projection $P_\Delta$ is not in $\Gr(\tilde A)$
since it does not act as a pseudo--differential operator on $E_{|\del\mcut}$.
Namely, since \eqref{ML-G5.2.5} contains off--diagonal terms the two copies of
of $N\subset \del\mcut$ interact and hence $P_\Delta$ is a
Fourier integral operator. However, $P_\Delta$ is pseudo--differential
on the bundle $E_{|N}\oplus E_{|N}$ over $N$.
This is only a mild generalization
of the situation of Section \plref{sec2} and we refrain from formalizing
it. From now on $\Gr(\tilde A)$ is to be understood as the set of those
pseudo--differential operators on the bundle $E_{|N}\oplus E_{|N}$ over $N$
which satisfy \eqref{grassm2}, \eqref{grassm3}, \eqref{grassm4} of
Definition \plref{SAFG}.
It is fairly clear that the results of the previous sections also
apply to this situation.

There is a natural map 
\begin{equation}\label{eq5}
\Gr(A)\times \Gr(-A)\longrightarrow \Gr(\tA),\quad
(P,Q)\mapsto\begin{pmatrix}P&0\\0&Q\end{pmatrix}
\end{equation}
with respect to the decomposition
\eqref{decompl2}. In particular the \Calderon\  projector for $\mcut$ takes
the form
\begin{equation}
P_{\mcut}=\begin{pmatrix}P_{M^+}&0\\ 0&P_{M^-}\end{pmatrix}.
\label{ML-G5.3}
\end{equation}

\theoremstyle{definition}
\newtheorem*{warning}{Warning}
\begin{warning}
1. There are two different decompositions of
$L^2(E_{|\del \mcut})$, one coming from the $\pm i$ eigenspace
decomposition of $\tga$ \eqref{decompL2}, and the second from
the decomposition $\del\mcut= N\coprod N$ \eqref{decompl2}.
This leads to  two different  matrix representations of
$P\in \Gr(\tA)$.  These two decompositions are compatible since
$\tga=\gamma\oplus (-\gamma)$, and so in fact one can write
$$ L^2(E_{|\del \mcut})=(\mathcal{E}_i\oplus \mathcal{E}_{-i})
\oplus( \mathcal{E}_{-i}
\oplus \mathcal{E}_i).$$

2. Although $P_\Delta\in\Gr(\tilde A)$, it is \emph{not} in $\Gr_\infty(\tilde
    A)$. This fact causes technical difficulties.

3. It follows from \eqref{ML-G5.2} that if one parameterizes the collar
of $M^-$ as $\partial M^-\times[0,\eps)$ then $\gamma$ is replaced by
$-\gamma$ and $A$ is replaced by $-A$. This in particular means that
the natural symplectic structure on $L^2(E_{|\partial M^-})$ is
induced by $-\gamma$.
Sometimes it will be crucial to distinguish between the map $\Phi_\gamma$
and the map $\Phi_{-\gamma}$ (cf. \eqref{phimap}). The relation between
the two is
\begin{equation}
     \Phi_{-\gamma}(P)=-\Phi_\gamma(I-P)^*,\quad P\in\Gr(A).\label{ML-G5.5}
\end{equation}
\end{warning}

\begin{lemma}\label{splitbundles} Let $D$ be a Dirac operator on
$\mcut$ and  suppose that $ P_t, 0\le t\le 1,$
is a continuous path in $\Gr(A)$ and $Q_t, 0\le t\le 1,$ is
a continuous path
in $\Gr(-A)$. Let
$$B_t=\begin{pmatrix} P_t&0\\0&Q_t\end{pmatrix}$$ be the corresponding
path in $\Gr(\tA)$.
Then
$$\SF(D_{B_t},\mcut)_{t\in[0,1]}=\SF(D_{P_t},M^+)_{t\in[0,1]}+
\SF(D_{Q_t},M^-)_{t\in[0,1]}.$$
\end{lemma}
\begin{proof} This follows from the fact that
\begin{equation}
L^2(E,\mcut)=L^2(E,M^+)\oplus L^2(E, M^-),
\label{ML-G5.4}
\end{equation}
and $D, B_t$ preserve this
splitting. Hence, $D_{B_t}=D_{P_t}^{M^+}\oplus D_{Q_t}^{M^-}$.
Note that the splitting \eqref{ML-G5.4} induces the splitting
\eqref{decompl2} by restricting to the boundary.\end{proof}

Notice that $P\in \Gr(A)$ if and only if $I-P\in \Gr(-A)$. Therefore,  a
particularly symmetric  family of  boundary conditions for $D$
acting on
$\mcut$ is given by the image of $(P,I-P)$ under the map \eqref{eq5}.

Corollary \plref{lifting} to the Scott--Wojciechowski theorem implies
the following
lemma.

\begin{lemma}\label{swconc} Let $D$ be a Dirac operator over $M$ and  let
$P_0, P_1\in \Gr_\infty(A)$. Choose a smooth path $P_t\in\Gr_\infty(A), 0\le
t\le 1,$ from $P_0$ to $P_1$ and put
$$Q_t:=\begin{pmatrix}P_t&0\\ 0& I-P_t\end{pmatrix}\in \Gr_\infty(\tA).$$

Then
$$\etab(D_{Q_1},\mcut)-\etab(D_{Q_0},\mcut)=
\SF(D_{P_t},M^+)_{t\in[0,1]}+\SF(D_{I-P_t},M^-)_{t\in[0,1]}.$$
In particular, the quantity $\SF(D_{P_t},M^+)_{t\in[0,1]}+
\SF(D_{I-P_t},M^-)_{t\in[0,1]}$ is
independent of the choice of the path $P_t$.
\end{lemma}
\begin{proof} We know from Proposition \plref{connected} that
$\Gr_\infty(A)$ is path connected.
This assures the existence of
a path $P_t$. Furthermore, notice that
\begin{equation}\label{version3}\etab(D_{Q_t},\mcut)=\etab(D_{P_t},M^+)+\etab(D_
{I-P_t},M^-)\end{equation} since $D$ and $Q_t$ preserve  the splitting of
$L^2(E_{|\mcut})=L^2(E_{|M^+})\oplus L^2(E_{|M^-})$.

Lemma \plref{sf-eta} and Corollary \plref{lifting} imply that
\begin{equation}\label{version1}
\begin{split}
\etab(D_{P_1}&,M^+)-\etab(D_{P_0},M^+)- \SF(D_{P_t},M^+)_{t\in[0,1]}\\
&=\frac{1}{2}\int_0^1\frac{d}{dt}(\eta(D_{P_t}))dt
   = \frac{1}{2\pi
i}\int_0^1\frac{d}{dt}\log\detf\bigl(\Phi(P_t)\Phi(P_0)^*\bigr)dt,
\end{split}
\end{equation}
and
\begin{equation}\label{version2}
\begin{split}
\etab(D_{I-P_1}&,M^-)-\etab(D_{I-P_0},M^-)- \SF(D_{I-P_t},M^-)_{t\in[0,1]}\\
&=\frac{1}{2}\int_0^1\frac{d}{dt}(\eta(D_{I-P_t}))dt
   =\frac{1}{2\pi
i}\int_0^1\frac{d}{dt}\log\detf\bigl(\Phi_{-\gamma}(I-P_t)\Phi_{-\gamma}(I-P_0)^*\bigr)dt.
\end{split}
\end{equation}
Note that in \eqref{version1} $\Phi$ is taken with respect to $\gamma$ and in
\eqref{version2} $\Phi$ is taken with respect to $-\gamma$ (cf. 3. of the
warning above). In view of \eqref{ML-G5.5} we find
\begin{equation}\label{version4}
 \detf(\Phi_{-\gamma}(I-P_t)\Phi_{-\gamma}(I-P_0)^*)=
\detf(\Phi(P_t)^*\Phi(P_0))
=\overline{\detf(\Phi(P_t)\Phi(P_0)^*)},
\end{equation}
and consequently,
\begin{equation}
\frac{d}{dt}\log\detf\bigl(\Phi_{-\gamma}(I-P_t)\Phi_{-\gamma}(I-P_0)^*\bigr)
=- \frac{d}{dt}\log\detf\bigl(\Phi(P_t)\Phi(P_0)^*\bigr).
\end{equation}
Adding \eqref{version1} and \eqref{version2}  and using
\eqref{version3} gives the desired formula.
\end{proof}

For any $P\in\Gr(A)$ a natural path connecting $P_\Delta$ and
$\begin{pmatrix} P&0\\ 0&I-P\end{pmatrix}$ is given by (cf.
\cite[Sec. 3]{BruLes:EIN})
\begin{equation}\label{defofptheta}P(\theta,P):=
\begin{pmatrix}
\cos^2(\theta)P+\sin^2(\theta)(I-P)&-\cos(\theta)\sin(\theta)\\
-\cos(\theta)\sin(\theta)& \cos^2(\theta)(I-P)+\sin^2(\theta)P
\end{pmatrix}.\end{equation}
A straightforward calculation shows that $\xi=\begin{pmatrix}\xi_+\\
\xi_-\end{pmatrix}\in
\ker P(\theta,P)$ if and only if
\begin{equation}\begin{split}
          \cos(\theta)P\xi_+&=\sin(\theta)P\xi_-,\\
          \sin(\theta)(I-P)\xi_+&=\cos(\theta)(I-P)\xi_-.
        \end{split}
\label{ML-G5.10}
\end{equation}

\begin{lemma} Let $P\in\Gr(A)$. If $\cos(\theta)\not=0$ then the map
$P(\theta,P)$ lies in $\Gr(\tA)$.
Furthermore
$$P(0,P)= \begin{pmatrix} P&0\\ 0&I-P\end{pmatrix}  \ \ \hbox{ and } \ \
P(\tfrac{\pi}{4},P)= P_\Delta.$$
\end{lemma}
\begin{proof}  Fix a Lagrangian subspace $L\subset \ker A$ and let
$P^+=P^+(L)$. The only part which is not straightforward is the claim that
$(P(\theta,P),
\tilde P^+)$,
$\tilde P^+:=\begin{pmatrix} P^+&0\\0 &I-P^+\end{pmatrix}$ is a Fredholm
pair.  We use the following criterion
(cf. \cite[Remark 3.5]{BruLes:BVPI}). \begin{quote}
Two orthogonal projections $Q,R$ in
a Hilbert space form a Fredholm pair (invertible pair) if and only if
the operator
$QRQ+(I-Q)(I-R)(I-Q)$ is Fredholm (invertible).\end{quote}
One calculates
\begin{equation}\begin{split}
    &\tilde P^+ P(\theta)\tilde P^++(I-\tilde
P^+)(I-P(\theta,P))(I-\tilde P^+)\\
&= \bigl(\cos^2(\theta)(P^+PP^++(I-P^+)(I-P)(I-P^+))+\\
&\quad        \sin^2(\theta)(P^+(I-P)P^++(I-P^+)P(I-P^+))\bigr)\otimes
      \begin{pmatrix}I&0\\0&I\end{pmatrix}
\\ &\ge  \cos^2(\theta)(P^+PP^++(I-P^+)(I-P)(I-P^+))\otimes
      \begin{pmatrix}I&0\\0&I\end{pmatrix}.\end{split}\end{equation}
Hence if $\cos(\theta)\not=0$ then the pair $(P(\theta,P),\tilde P^+)$ is
Fredholm
(invertible) if the pair $(P,P^+)$ is Fredholm (invertible).
\end{proof}

We emphasize that even if $P\in\Gr_\infty(A)$ then
$P(\theta,P)\not\in\Gr_\infty(A)$
if $\sin(\theta)\not=0$. The significance of the family $P(\theta,P)$ stems
from the fact that $D_{P_\Delta}$ is naturally unitarily equivalent to $D$
acting on the closed manifold $M$.

We first note some consequences of the existence of the family $P(\theta,P)$
which do not make use of $\eta$--functions.

\begin{prop}\label{ML-S5.4} Let $D$ be a Dirac operator over $M$ and
let $P_0,P_1\in\Gr(A)$.
Choose a smooth path $P_t\in\Gr(A), 0\le t\le 1$, from $P_0$ to $P_1$ and put,
as in Lemma \plref{swconc},
\begin{equation}
Q_t:=\begin{pmatrix}P_t&0\\ 0& I-P_t\end{pmatrix}\in \Gr_\infty(\tA).
\label{ML-G5.12}
\end{equation}
Then
\[\SF(D_{P(\theta,P_1)},\mcut)_{\theta\in[0,\frac{\pi}{4}]}
   -\SF(D_{P(\theta,P_0)},\mcut)_{\theta\in[0,\frac{\pi}{4}]}
   +\SF(D_{Q_t},\mcut)_{t\in[0,1]}=0.\]
\end{prop}
\begin{proof} Note again that in view of Proposition \plref{connected}
the space $\Gr(A)$ is path connected.
Using $P_t$ one obtains a map $H: [0,\frac{\pi}{4}]\times [0,1]\to
\Gr(\tA)$:
\begin{equation}
H(\theta,t)= P(\theta,P_t).
\end{equation}
Since $H(\frac{\pi}{4},t)$ is the constant map at $P_\Delta$, one sees
that the path
\[\theta\mapsto H(\theta,0)=P(\theta,P_0),\quad 0\le \theta\le \frac{\pi}{4},
\]
is homotopic to the composite of the paths
\[
t\mapsto H(0,t)=Q_t,\quad 0\le t\le 1,
\]
and
\[
\theta\mapsto H(\theta,1)= P(\theta, P_1).
\]
The claim now follows from the homotopy invariance and additivity of
the spectral flow.
\end{proof}

\begin{prop}\label{caldconst} For the \Calderon\  projectors
$P_{M^+}$ of $M^+$ and $P_{\mcut}$ of $\mcut$, the space  $\ker
P(\theta, P_{M^+})\cap \ker
\tga(P_{\mcut})$ is canonically isomorphic to $\im P_{M^+}\cap\im P_{M^-}$.
In particular, its dimension is independent of $\theta\in
[0,\frac{\pi}{4}]$. Moreover,
\begin{equation}
\SF(D_{P(\theta,P_{M^+})},\mcut)_{\theta\in [0,\frac{\pi}{4}]}=0.
\label{ML-G5.18}
\end{equation}
Furthermore, if $P\in\Gr(A)$ and $P_t, 0\le t\le 1,$ is a smooth path in
$\Gr(A)$ from $P$ to the \Calderon\ projector $P_{M^+}$ then
\begin{equation}
\SF(D_{P(\theta,P)},\mcut)_{\theta\in[0,\frac{\pi}{4}]}=\SF(D_{Q_t},\mcut)_{t\in[0,1]},
\label{ML-G5.19}
\end{equation}
where $Q_t=P_t\oplus (I-P_t)$ as in \eqref{ML-G5.12}.
\end{prop}
\begin{proof} By Lemma \plref{kerdp} the space $\ker P(\theta, P_{M^+})\cap \ker
\tga(P_{\mcut})$ is isomorphic to $\ker D_{P(\theta,P_{M^+})}$.
Hence, if we can
show that  $\ker P(\theta, P_{M^+})\cap \ker \tga(P_{\mcut})$ is independent
of $\theta$, then $\SF(D_{P(\theta,P_{M^+})},\mcut)_{\theta\in
[0,\frac{\pi}{4}]}=0.$

Consider $\xi=\begin{pmatrix} \xi_+\\ \xi_-\end{pmatrix}\in \ker
P(\theta, P_{M^+})\cap \ker
\tga(P_{\mcut})$. In view of \eqref{ML-G5.3} and \eqref{ML-G5.10} this means
\begin{equation}\begin{split}
    &\cos(\theta)P_{M^+}\xi_+=\sin(\theta)P_{M^+}\xi_-,\\
    &0=\sin(\theta)(I-P_{M^+})\xi_+=\cos(\theta)(I-P_{M^+})\xi_-.
        \end{split}
\label{ML-G5.14}
\end{equation}
Since $\cos(\theta)\not=0$ we infer $\xi_-\in\im P_{M^+}\cap\im P_{M^-}$.

Conversely, given $\xi_-\in\im P_{M^+}\cap\im P_{M^-}$ put
$\xi_+:=\tan(\theta)\xi_-$. Then \eqref{ML-G5.14} implies that $\begin{pmatrix}
\xi_+ \\ \xi_- \end{pmatrix}\in\ker P(\theta, P_{M^+})\cap \ker
\tga(P_{\mcut})$.

\eqref{ML-G5.19} is an immediate consequence of \eqref{ML-G5.18} and
Proposition
\plref{ML-S5.4}.
\end{proof}

\begin{cor}\label{ML-S5.6} Let $M$ be a split manifold and let $D(t), a\le t\le
b,$ be a smooth path of Dirac type operators such that in a collar of the
separating
hypersurface we have $D(t)=\gamma(\frac{d}{dx}+A(t))$. Let $P_{M^+}(t)$ be
the corresponding family of \Calderon\ projectors. Then
\[ \SF(D(t))_{t\in[a,b]}=\SF(D_{I-P_{M^+}(t)}(t),M^-)_{t\in[a,b]}.\]
\end{cor}
\begin{proof} We note that it was proved in \cite{Nic:MIS} that
$P_{M^+}(t)$ is smooth. Consider the two parameter family of operators on
$\mcut$
\[
(D_{P(\theta,P_{M^+}(t))}(t),\mcut)_{0\le\theta\le\frac{\pi}{4},a\le
t\le b}.\]
By Proposition \plref{caldconst} for fixed $t$ the dimension of the kernel
of $D_{P(\theta,P_{M^+}(t))}(t)$ is independent of $\theta$. By the homotopy
invariance
of the spectral flow this implies
\[
\SF(D_{P(0,P_{M^+}(t))}(t),\mcut)_{t\in[a,b]}=\SF(D_{P(\frac{\pi}{4},P_{M^+}(t))}(t),\mcut)_{t\in[a,b]}.\]
Since $P(\frac{\pi}{4})=P_\Delta$ the right hand side equals
$\SF(D(t))_{t\in[a,b]}$. The left hand side equals
\[
\SF(D_{P_{M^+}(t)}(t),M^+)_{t\in [a,b]}+\SF(D_{I-P_{M^+}(t)},M^-)_{t\in[a,b]}
\]
and since $D_{P_M^+(t)}(t)$ is invertible its spectral flow vanishes and we
reach the desired conclusion.
\end{proof}
\begin{remark} We emphasize that we did not use $\eta$--invariants
to prove Proposition \plref{ML-S5.4}, Proposition \plref{caldconst}, and
Corollary \plref{ML-S5.6}. The only ingredients of the proof are the
family $P(\theta,P)$ and basic properties of the spectral flow.
\end{remark}

We now return to the discussion of $\eta$--invariants.
Since $D_{P_\Delta}$ is naturally unitarily equivalent to $D$ acting on
the closed manifold $M$,  we have for any $P\in \Gr(A)$
\begin{equation}\label{eta1}
\eta(D_{P(\frac{\pi}{4},P)},\mcut)=\eta(D_{P_\Delta},\mcut)=\eta(D,M).
\end{equation}
On the other hand, $D_{P(0,P)}$ is the
direct sum of $D_P$ acting on $M^+$ and $D_{I-P}$ acting on $M^-$. Therefore,
\begin{equation}\label{eta2}
\eta(D_{P(0,P)},\mcut)=\eta(D_P,M^+)+\eta(D_{I-P},M^-).\end{equation}

Hence, by Lemma \plref{sf-eta} we have
\begin{equation}\begin{split}
\eta(D,M)&= \eta(D_P,M^+)+\eta(D_{I-P},M^-)\\
&\quad+\frac12
\int_0^{\frac{\pi}{4}}\frac{d}{d\theta}\eta(D_{P(\theta,P)},\mcut)dt
+\SF(D_{P(\theta,P)})_{\theta\in[0,\frac{\pi}{4}]}.
        \end{split}
\label{ML-G5.13}
\end{equation}

Thus, in order to obtain a splitting theorem for the $\eta$--invariant
one needs to understand the last two terms on the right hand side
of \eqref{ML-G5.13}. If $P$ is the \Calderon\ projector of $M^+$ or
$M^-$ then by Proposition \plref{caldconst} the spectral flow term vanishes.

Consider now the Atiyah--Patodi--Singer projection
$P^+=P^+(L)$ of \eqref{pplus}.  The following theorem is the main
result of the article  \cite{BruLes:EIN} by J. Br\"uning and M. Lesch
(\cite[Theorem 3.9]{BruLes:EIN}, see also  (3.68) of loc.cit. with $T_+=-T_-^*$
determined by the choice of Lagrangian $L\subset\ker A$).

\begin{theorem}\label{brule} Let
$P^+= P^+(L)$ be the Atiyah--Patodi--Singer projection and let
$P(\theta,P^+)_{\theta\in[0,\frac{\pi}{4}]}$ the deformation
\eqref{defofptheta} to the
continuous transmission projection. Then
$$\frac{d}{d\theta}(\eta(D_{P(\theta,P^+)}, \mcut)=0.$$\qed
\end{theorem}

In view of \eqref{ML-G5.13} we conclude from Theorem
\plref{brule} that
\begin{equation}\label{p+eqn}
\etab(D,M)-\etab(D_{P^+},M^+)-\etab(D_{I-P^+},M^-)= \SF(D_{P(\theta,P^+)},
\mcut)_{\theta\in[0,\frac{\pi}{4}]}.
\end{equation}
Since the right hand side of \eqref{p+eqn} is an integer this formula
implies
the $\mod \Z$ gluing formula for the $\eta$--invariant
(see \cite{BruLes:EIN} for a discussion of the history of this result).
Note that \eqref{p+eqn} is slightly weaker than Theorem \plref{brule}.

Our strategy to obtain a useful splitting theorem for the
$\eta$--invariant can now be explained. On the one hand \eqref{p+eqn}
gives a complete splitting formula for the
$\eta$--invariant with respect to Atiyah--Patodi--Singer boundary
conditions, but   it contains  the (in general) uncomputable term $
\SF(D_{P(\theta,P^+)})_{\theta\in[0,\frac{\pi}{4}]}$. On the other hand if
we were to replace $P^+$ by the \Calderon\  projector
$P_{M^+}$, Proposition \plref{caldconst} shows that the corresponding
spectral flow term vanishes. Thus Theorem \plref{brule} (or at least
\eqref{p+eqn}) needs to be
extended to a more general class of projections, including the \Calderon\
projector.
One possible strategy would be to generalize the
arguments of \cite{BruLes:EIN} to more general projections. This
might be manageable but technically tedious.
Here, we will use a simpler approach which shows slightly less.
Lemma \plref{swconc} and Proposition \plref{caldconst}
lead to a generalization of \eqref{p+eqn} to
projections in $\Gr_\infty(A)$. This is less than a generalization
of Theorem \plref{brule} since the variation of
the
$\eta$--invariant with respect to the path $P(\theta,P)$ might be
non-zero.

\begin{theorem}\label{symsplit} Let $D$ be a Dirac operator on
$M$ and let $N\subset M$ split $M$ into $M^+$ and $M^-$. We assume
that in a collar neighborhood $[-\eps,\eps]\times N$ of $N$, $D$ has
the form $D=\gamma(\frac{d}{dx}+A)$ as in \eqref{ML-G2.1}.
Let
$P\in\Gr_\infty(A)$ and let $P_t$ be a smooth path in $\Gr_\infty(A)$
from $P$ to the \Calderon\ projector $P_{M^+}$. As in \eqref{ML-G5.12}
put $Q_t:=P_t\oplus (I-P_t)$. Then
\begin{equation}
\begin{split}\etab(D,M)&=\etab(D_{P},M^+)+\etab(D_{I-P},M^-)+
\SF(D_{P(\theta,P)},
\mcut)_{\theta\in[0,\frac{\pi}{4}]}\\
   &=\etab(D_{P},M^+)+\etab(D_{I-P},M^-)+ \SF(D_{Q_t},\mcut)_{t\in [0,1]}.
   \end{split}\label{ML-G5.22}
\end{equation}
In particular, if $P_{M^+}$  is the \Calderon\  projector  for $M^+$ then
$$\etab(D,M)=\etab(D_{P_{M^+}},M^+) +
\etab(D_{I-P_{M^+}},M^-).$$
\end{theorem}

\begin{proof}  Fix a Lagrangian subspace $L\subset \ker A$ and choose a
smooth path
$R_t\in\Gr_\infty(A)$ from
$P$ to the Atiyah--Patodi--Singer projection $P^+= P^+(L)$. Set $\tilde
R_t:=R_t\oplus (I-R_t)$. By Proposition \plref{ML-S5.4} we have
$$\SF(D_{P(\theta,P)},\mcut)_{\theta\in [0,\frac{\pi}{4}]}=
\SF(D_{P(\theta,P^+)},\mcut)_{\theta\in [0,\frac{\pi}{4}]}
+\SF(D_{\tilde R_t},\mcut)_{t\in[0,1]}.$$
Using Lemma \plref{swconc} and \eqref{p+eqn} we obtain
\begin{equation}\label{proof1}
\begin{split}
    &\etab(D,M)-\etab(D_{P},M^+)-\etab(D_{I-P},M^-)\\
    &=
    \etab(D,M)-\etab(D_{P^+},M^+)-\etab(D_{I-P^+},M^-)+\etab(D_{\tilde
    R_1},\mcut)-\etab(D_{\tilde R_0},\mcut)\\
    &=\SF(D_{P(\theta,P^+)},\mcut)_{\theta\in
    [0,\frac{\pi}{4}]}+\SF(D_{\tilde R_t},\mcut)_{t\in[0,1]}\\
    &=\SF(D_{P(\theta,P)},\mcut)_{\theta\in [0,\frac{\pi}{4}]}.
\end{split}
\end{equation}
This proves the first line of \eqref{ML-G5.22}. The second line of
\eqref{ML-G5.22} and the last assertion follow from Proposition
\plref{caldconst}.
\end{proof}

Notice that by symmetry the same argument also shows that
$$\etab(D,M)=\etab(D_{I-P_{M^-}},M^+) +
\etab(D_{P_{M^-}},M^-).$$

Applying Theorem \plref{invertible} allows us to extend Theorem
\plref{symsplit} as follows.

\begin{theorem}\label{invertiblecase}Let $D$ be a Dirac operator on
$M$ and let $N\subset M$ split $M$ into $M^+$ and $M^-$.
Then for $P\in \Gr_\infty(A)$ and $Q\in \Gr_\infty(-A)$
we have, with $\Phi=\Phi_\gamma$,
\[
\begin{split}
\etab(&D,M)-\etab(D_P,M^+) - \etab(D_{Q},M^-)\\
&=-\tfrac{1}{2\pi i}\tr\log(\Phi(P)\Phi(P_{M^+})^*)-\tfrac{1}{2\pi
i}\tr\log(\Phi(P_{M^-})\Phi(Q)^*)\\
&\quad +\tfrac{1}{2\pi i}\tr\log(\Phi(I-P_{M^-})\Phi(P_{M^+})^*).
\end{split}
\]
In particular,
 \begin{equation}\label{iceq1}\etab(D,M)=\etab(D_{P_{M^+}},M^+) +
\etab(D_{P_{M^-}},M^-)+\tfrac{1}{2\pi
i}\tr\log(\Phi(I-P_{M^-})\Phi(P_{M^+})^*).
\end{equation}
\end{theorem}
\begin{proof}  It was remarked after \eqref{calderon} that
 $P_{M^+}\in\Gr_\infty(A)$
and $P_{M^-}\in\Gr_\infty(-A)$. Consequently, $I-P_{M^-}\in\Gr_\infty(A)$ and
hence $I-P_{M^-}-P_{M^+}$ is trace class.

Theorem \plref{symsplit} implies that
$\etab(D,M)-\etab(D_P,M^+) -
\etab(D_{Q},M^-)$ is equal to
\begin{equation}
\label{diff1}
\begin{split}
  & -(\etab(D_P,M^+)-\etab(D_{P_{M^+}},M^+))
-(\etab(D_Q,M^-)-\etab(D_{P_{M^-}},M^-))\\
  & -(\etab(D_{P_{M^-}},M^-)-\etab(D_{I-P_{M^+}},M^-)).
\end{split}
\end{equation}
Applying Theorem \plref{invertible} to the three summands
in \eqref{diff1} and taking \eqref{ML-G5.5} into account gives the assertion.
\end{proof}

The formula \eqref{iceq1} expresses the $\etab$--invariant of $D$ on $M$ in terms
of two $\etab$--invariants  intrinsic to the two pieces $M^+$ and $M^-$
of the decomposition and an ``interaction'' term.

\section{Maslov index and winding number}
\label{sec6}

In this section we compile the necessary material about the Maslov index
and the winding number.  One important comment is that in constructing the
various invariants (winding number, Maslov index, triple index,   spectral
flow, and the branch of the logarithm) conventions must be chosen to set signs
and  to handle degenerate cases. In particular, care must be taken to ensure
that the different possible conventions are chosen compatibly.  Thus, although
some of the material we present here is a generalization of ideas which appear
in the literature, the subtleties arising in  organizing the conventions
compatibly and extending the constructions from the finite--dimensional to the
infinite--dimensional context require the careful exposition we present.

\newcommand\wind{\operatorname{wind}}
\subsection{Winding number}
\newcommand{\mck}{\mathcal{K}}
Let $H$ be a complex Hilbert space  and denote by
$\cu(H)$ the group of unitary operators on $H$. Similarly to \eqref{ML-G2.7},
\eqref{ML-G2.8} we introduce the following subspaces:
\begin{equation}\begin{split}
     \cu_*(H)&:=\bigsetdef{U\in\cu(H)}{-1\not\in\spec U},\\
     \cu_{\Fred}(H)&:=\bigsetdef{U\in\cu(H)}{-1\not\in\specess U},\\
     \cu_{\mck}(H)&:=\bigsetdef{U\in\cu(H)}{U-I \text{ is compact }},\\
     \cu_{\tr}(H)&:=\bigsetdef{U\in\cu(H)}{U-I \text{ is trace class }}.
        \end{split}
\end{equation}
The spaces $\cu_*(H)$ and $\cu_{\Fred}(H)$ are not groups.    It is well--known
that the inclusion
$\cu_{\tr}(H)\hookrightarrow\cu_{\mck}(H)$ is a homotopy equivalence and that
$\cu_\mck(H)$ is homotopy equivalent to the infinite unitary group
$\cu(\infty)=\lim\limits_{n\to\infty}\cu(n)$. Therefore, one has by Bott
periodicity
\begin{equation}\begin{split}
     \pi_{2k}(\cu_{\mck}(H))&=\pi_{2k}(\cu_{\tr}(H))=0,\\
     \pi_{2k+1}(\cu_{\mck}(H))&=\pi_{2k+1}(\cu_{\tr}(H))\simeq\Z,
        \end{split}
\qquad
     k=0, 1, 2, ...
\end{equation}
Furthermore, the isomorphism $\pi_1(\cu_{\tr}(H))\to \Z$ is given
by the winding number. I.e. if $f:[0,1]\to \mcU_{\tr}(H)$
is a closed $C^1$--path then
\begin{equation}
     \wind(f):=\frac{1}{2\pi i}\int_0^1 \tr(f(t)^{-1}f'(t))dt.
  \label{winddef}
\end{equation}

\begin{lemma}\label{ML2-S6.1}\hfill
\begin{enumerate}
\item The inclusion $\cu_{\mck}(H)\to \cu_{\Fred}(H)$ is a weak homotopy
equivalence,
\item For any $U\in\cu_{\Fred}(H)$ there exists a smooth path
$f:[0,1]\to\cu_{\Fred}(H)$ such that $f(0)-I$ is of finite rank,
$f(1)=U$, and such that $\dim\ker(f(t)+I)$ is independent of $t$.
\end{enumerate}
\end{lemma}
\begin{proof} 1. Let $\cq(H):=\cb(H)/\mck(H)$ be the Calkin algebra. Then
the quotient map $\sigma:\cb(H)\to \cq(H)$ sends $\cu_{\Fred}(H)$
onto $\bigsetdef{u\in \cq(H)}{-1\not \in \spec u}=:\mcU_*\cq(H)$. Moreover,
$\cu_{\mck}(H)$ acts freely (from left and right) on the fibers. Thus
we obtain a fibration $\cu_{\mck}(H)\to \cu_{\Fred}(H)\to \cu_*\cq(H)$.
The claim now follows since $\cu_*\cq(H)$ is contractible.

To see the latter we note that for any $C^*$--algebra $\ca$ the set
$\setdef{u\in \ca}{u \text{ unitary}, -1\not\in\spec u}$ is contractible.
The contraction is given by $H_t(u):=\exp(t\log u)$, $0\le t\le 1$.
This is well--defined since $-1\not\in\spec u$.

2. Let $H=\ker(U+I)\oplus H_1=:H_0\oplus H_1$. Then $U$ splits into
$U=-I_{H_0}\oplus
\tilde U$ and $-1\not\in\spec\tilde U$. Now put $f(t):=-I_{H_0}\oplus \exp(t
\log\tilde U)$.
\end{proof}

In view of this lemma the winding number \eqref{winddef} extends to
a group isomorphism
\[\wind:\pi_1(\cu_{\Fred}(H))\to \Z.\]
Next we define the winding number
for not necessarily closed paths in $\cu_{\Fred}(H)$. Namely, as it was
noted in the previous proof the space $\cu_*(H)$ is contractible. Hence
the natural map $\pi_1(\cu_{\Fred(H)})\to \pi_1(\cu_{\Fred}(H),\cu_*(H))$
is a bijection and thus we obtain a winding number defined for
curves $f: ([0,1],\{0,1\})\to (\cu_{\Fred}(H),\cu_*(H))$.
More concretely, if $f$ is such a curve then one chooses
$\tilde f:[0,1]\to \cu_*(H)$ with $\tilde f(0)=f(1)$ and
$\tilde f(1)=f(0)$. Then $f*\tilde f$ is a closed curve
in $\cu_{\Fred}(H)$ and one puts $\wind(f):=\wind(\tilde f*f)$.
Since $\cu_*(H)$ is contractible it is clear that $\wind(f)$ is
well--defined independently of the choice of $\tilde f$.

Finally, we choose a convention to define the winding number for a
curve whose endpoints do not lie in $\cu_*(H)$: let $f:[0,1]\to
\cu_{\Fred}(H)$ be a continuous curve. $-1$ is an isolated point
in the spectrum of $f(t)$ since $-1\not\in\specess(f(t))$. We may
therefore choose an $\eps>0$ such that for all
$\varphi\in[-\eps,\eps], \varphi\not=0,$ we have
$-1\not\in\spec(f(j)e^{i\varphi}), j=0,1$. Now define
\begin{equation}
     \wind(f):=\wind(f e^{-i\eps}).
\end{equation}
The   winding number has the following properties:

\begin{enumerate}
\item \textit{Path Additivity}: Let $f_1, f_2:[0,1]\to\cu_{\Fred}(H)$
be continuous paths with $f_2(0)=f_1(1)$. Then
\[\wind(f_1*f_2)=\wind(f_1)+\wind(f_2).\]
\item \textit{Homotopy invariance:} Let $f_1,f_2$ be continuous paths
in $\cu_{\Fred}$. Assume that there is a homotopy
$H:[0,1]\times[0,1]\to\cu_{\Fred}$
such that $H(0,t)=f_1(t), H(1,t)=f_2(t)$ and such that
$\dim \ker(H(s,0)+I), \dim \ker(H(s,1)+I)$ are independent of $s$.
Then
$\wind(f_1)=\wind(f_2)$.
\item If $f:[0,1]\to\cu_{\tr}(H)$ is a $C^1$--curve then
\begin{equation}\wind(f)=\frac{1}{2\pi i}\Bigl(\int_0^1
\tr(f(t)^{-1}f'(t))dt
              -\tr(\log f(1))+\tr(\log
f(0))\Bigr),\label{ML-G6.10}\end{equation}
where the logarithm is normalized as in \eqref{ML-G4.11}.
\end{enumerate}

We note in passing that the winding number may be interpreted as a
spectral flow
across $-1$ \cite{BooFur:MIF}, \cite{Phi:SAF}. Namely, the winding
number of a path
$f:[0,1]\to\cu_{\Fred}(H)$ can be calculated as follows: choose a subdivision
$0=t_0<t_1<...<t_n=1$ and $0<\eps_j<\pi$, $j=0,...,n-1$, such that
$-e^{i\varphi}\not\in\specess(f(t))$ for $t\in[t_j,t_{j+1}]$ and
$|\varphi|\le\eps_j$ and moreover $-e^{\pm i\eps_j}\not\in\spec f(t)$
for $t\in[t_j,t_{j+1}]$. Then put
\begin{equation}\begin{split}
    \wind(f(t))_{t_j\le t\le t_{j+1}}&:= \#(\spec(f(t_{j+1}))\cap
    \bigsetdef{-e^{i\varphi}}{0<\varphi<\eps_j})\\
&\quad-\#(\spec(f(t_j))\cap \bigsetdef{-e^{i\varphi}}{0<\varphi<\eps_j},
        \end{split}\label{ML-G6.11}
\end{equation}
where eigenvalues are counted with multiplicity. Finally,
\begin{equation}
    \wind(f)=\sum_{j=0}^{n-1}\wind(f(t))_{t_j\le t\le t_{j+1}}.
\label{ML-G6.12}
\end{equation}

\begin{dfn}\label{doubleindexdfn} Let $U\in\cu_{\mck}(H)$ and
$V\in\cu_{\Fred}(H)$. Then
the {\em double index} $\tau_w(U,V)\in\Z$ is defined as follows:
choose continuous paths $f:[0,1]\to\cu_{\mck}(H), g:[0,1]\to\cu_{\Fred}(H)$
such that $f(0)=g(0)=I$ and $f(1)=U$, $g(1)=V$. Then put
$\tau_w(U,V):=\wind(f)+\wind(g)-\wind(fg)$.
$\tau_w(U,V)$ is defined accordingly if $U\in\cu_{\Fred}(H), V\in\cu_{\mck}(H)$
or $U,V\in\cu_{\Fred}(H), UV\in\cu_{\mck}(H)$.
\end{dfn}

\begin{prop}\label{doubleindex} The double index $\tau_w$ is
well--defined. It has the following properties:
\begin{enumerate}
\item \emph{(Homotopy invariance)}  If $f:[0,1]\to\cu_{\mck}(H), g:[0,1]\to
\cu_{\Fred}(H)$ are continuous paths then
\[   \tau_w(f(1),g(1))-\tau_w(f(0),g(0))=\wind(f)+\wind(g)-\wind(fg).\]
In particular, if $\dim\ker(f(t)+I),\dim\ker(g(t)+I)$  and $\dim\ker
(f(t)g(t)+I)$ are independent of
$t$ then $\tau_w(f(1),g(1))=\tau_w(f(0),g(0))$.
\item
If $U,V\in\cu_{\tr}(H)$ then
\begin{equation}
    \tau_w(U,V)=\frac{1}{2\pi i}\bigl( \tr\log UV-\tr\log U-\tr\log V\bigr).
\label{ML-G6.13}
\end{equation}
\item For any $U\in\cu_{\Fred}(H)$ we have
$$\tau_w(I,U)=\tau_w(U,I)=0 \text{ and }
\tau_w(U,U^{-1})=-\dim\ker(U+I).$$
\end{enumerate}
\end{prop}

\begin{proof} First note that if $U\in\cu_{\mck}(H)$ and $V\in\cu_{\Fred}(H)$
then since $U-I$ is compact one has $\specess(UV)=\specess(V)$,
in particular $UV\in\cu_{\Fred}(H)$. If $\tilde f:[0,1]\to\cu_{\mck}(H),
\tilde g:[0,1]\to\cu_{\Fred}(H)$ are different paths with $\tilde f(0)=\tilde
g(0)=I,
\tilde f(1)=U, \tilde g(1)=V$ then consider the closed
paths $f*\tilde f_-$
and $g*\tilde g_-$, where $\tilde f_-$ denotes the path $\tilde f$ traversed
in the opposite direction. Since the pointwise product of closed paths
$(f*\tilde f_-)(g*\tilde g_-)$
is homotopic to $(f*\tilde f_-)*(g*\tilde g_-)$ we find
\begin{equation}\begin{split}
     0&=\wind(f*\tilde f_-)+\wind(g*\tilde g_-)-\wind((f*\tilde
f_-)(g*\tilde g_-))\\
      &=-\wind(\tilde f)+\wind(f)-\wind(\tilde g)+\wind(g)+\wind(\tilde f\tilde
      g)-\wind(fg).
        \end{split}\label{Ml-G6.14}
\end{equation}
This shows that $\tau_w$ is well--defined. The homotopy invariance is
straightforward from the definition and the homotopy invariance of
the winding number.

  2. This assertion is a consequence of \eqref{ML-G6.10}.

3. That $\tau_w(I,U)=\tau_w(U,I)=0$ follows immediately from the definition.

For $U\in\cu_{\tr}(H)$ the third identity follows from Assertion 2. (note the
normalization
\eqref{ML-G4.11} of $\log$). If $U$ is arbitrary we apply Lemma
\plref{ML2-S6.1} 2. and choose
a continuous path $f:[0,1]\to\cu_{\Fred}(H)$ such that $f(1)=U,
f(0)\in\cu_{\tr}(H)$ and such that $\dim\ker(f(t)+I)$ is independent of $t$.
The claim now follows from the homotopy invariance 1.
\end{proof}

A priori $\tau_w$ cannot be defined on $\cu_{\Fred}(H)\times\cu_{\Fred}(H)$
(which might be desirable) since for $U,V\in\cu_{\Fred}(H)$ in general
$UV\not\in\cu_{\Fred}(H)$. Even if one assumes $UV\in\cu_{\Fred}(H)$ it is
in general not possible to choose paths $f,g$ as above such that
$f(t)g(t)\in\cu_{\Fred}(H)$
for all $t$.

\begin{cor}\label{windfinv} Let $f:[0,1]\to\cu_{\Fred}(H)$ be a continuous
path.
Then
\[ \wind(f)+\wind(f^{-1})= \dim\ker(f(0)+I)-\dim\ker(f(1)+I).\]
\end{cor}
\begin{proof} We apply Proposition \plref{doubleindex} 1. with $g=f^{-1}$
and obtain using Proposition \plref{doubleindex} 3.
   \[\begin{split}\wind(f)+\wind(f^{-1})&=
   \tau_w(f(1),f(1)^{-1})-\tau_w(f(0),f(0)^{-1})\\
      &=-\dim\ker(f(1)+I)+\dim\ker(f(0)+I).
     \end{split}
\]
\end{proof}

\subsection{Maslov Index}\label{secsixtwo} 
Let $(H, \langle \ , \ \rangle, \gamma)$ be a
Hermitian symplectic Hilbert space (cf. Def. \plref{defofsymp}). 
Thus $\gamma:H\to H$ is a unitary
map satisfying
$\gamma^2=-1$ and   the eigenspaces $\ce_{\pm i}:=\ker(\gamma\mp i)$
have the same Hilbert space dimension.   As in  Section \plref{sec2}
we denote by
\begin{equation*}
    \cl:=\bigsetdef{L\subset H}{L \textrm{ closed subspace, }  \gamma
L=L^\perp}
\end{equation*}
the set of Lagrangian subspaces. As usual $L\in \cl$ will be identified
with the orthogonal projection $P_L$ onto $L$. The image of an orthogonal
projection $P$ in $H$ is Lagrangian if and only if $\gamma P\gamma^*=I-P$.
Similarly as in Section
\plref{sec2} we put
\begin{equation}\begin{split}
     \Gr(H)&:=\bigsetdef{P\in\cb(H)}{P=P^*, P^2=P, \gamma P\gamma^*=I-P},\\
     \Gr^{(2)}_{\Fred}(H)&:=\bigsetdef{(P,Q)}{P,Q\in\Gr(H), (P,Q) \textrm{ are a
     Fredholm pair}},\\
     \Gr^{(2)}_*(H)&:=\bigsetdef{(P,Q)}{P,Q\in\Gr(H), (P,Q) \textrm{ is an
     invertible pair}},\\
     \Gr^{(2)}_{\mck}(H)&:=\bigsetdef{(P,Q)}{P,Q\in\Gr(H), P-Q
\textrm{ is compact}}.
        \end{split}\label{ML-G6.16}
\end{equation}
Notice that, in contrast to the definition of $\Gr(A)$, there is no Fredholm
assumption about elements of $\Gr(H)$. The corresponding spaces of Lagrangians
are
\begin{equation}\begin{split}
     \cl^{(2)}_{\Fred}&:=\bigsetdef{(L_1,L_2)}{L_1,L_2\in\cl,
(L_1,L_2)\textrm{ is
     Fredholm}},\\
     \cl^{(2)}_*&:=\bigsetdef{(L_1,L_2)\in\cl^{(2)}}{(L_1,L_2)
\textrm{ is invertible}}.
        \end{split}\label{ML-G6.17}
\end{equation}
Recall that a pair of Lagrangian spaces $(L_1,L_2)$ is Fredholm if $L_1\cap
L_2$ is finite--dimensional and if $L_1+L_2$ is closed with finite codimension,
and that the pair
$(L_1,L_2)$ is invertible if $L_1\cap L_2=\{0\}$ and $L_1+L_2=H$.

We emphasize the confusing fact that $(L_1,L_2)$ is Fredholm (resp. invertible)
if and only if the pair of projections $(I-P_{L_1},P_{L_2})$ is Fredholm
(resp. invertible). Therefore, a Fredholm pair of projections $(P,Q)$ will
sometimes be identified with the pair $(\ker P,\im Q)$ of Lagrangian subspaces.

As in Lemma \plref{pairs} one sees that with respect to the decomposition
$H=\ce_i\oplus \ce_{-i}$ each $P\in\Gr(H)$ takes the form
\begin{equation}
        P=\frac 12\begin{pmatrix}I& \Phi(P)^*\\ \Phi(P) &
I\end{pmatrix},\label{ML-G6.18}
\end{equation}
where $\Phi(P)\in\cu(\ce_i,\ce_{-i})$. Moreover, the map
\begin{equation}
     \Phi:\Gr(H)\longrightarrow \cu(\ce_i,\ce_{-i})\label{ML-G6.19}
\end{equation}
is a diffeomorphism. Furthermore, the pair $(P,Q)$ is Fredholm if and only if
$\Phi(P)\Phi(Q)^*\in\cu_{\Fred}(\ce_{-i})$ and $(P,Q)$ is invertible
if and only if $\Phi(P)\Phi(Q)^*\in\cu_*(\ce_{-i})$ (cf. Lemma \plref{pairs}).
Finally,
$P-Q$ is compact (resp. trace class) if and only if
$\Phi(P)\Phi(Q)^*\in\cu_{\mck}(\ce_{-i})$ (resp. $\cu_{\tr}(\ce_{-i})$).

\begin{prop}\label{homequ} There are diffeomorphisms
\begin{eqnarray*}\Gr^{(2)}_{\Fred}(H)&\cong& \cu_{\Fred}(\ce_{-i})\times
\cu(\ce_i,\ce_{-i}),\\
  \Gr^{(2)}_{\mck}(H)&\cong& \cu_{\mck}(\ce_{-i})\times
\cu(\ce_i,\ce_{-i}),\\
\Gr^{(2)}_*(H)&\cong& \cu_*(\ce_{-i})\times \cu(\ce_i,\ce_{-i}),\\
  \Gr^{(2)}_{\mck}(H)\cap\Gr^{(2)}_*(H)&\cong&
\bigl(\cu_{\mck}(\ce_{-i})\cap\cu_*(\ce_{-i})\bigr)\times \cu(\ce_i,\ce_{-i}).
\end{eqnarray*}

In particular the identifications induce homotopy equivalences
\begin{eqnarray*}(\Gr^{(2)}_{\Fred}(H),\Gr^{(2)}_*(H))&\simeq&
(\cu_{\Fred}(\ce_{-i}),
\cu_*(\ce_{-i}))\\
  (\Gr^{(2)}_{\mck}(H),\Gr^{(2)}_*(H)\cap\Gr_{\mck}(H))&\simeq&
(\cu_{\mck}(\ce_{-i}),
\cu_*(\ce_{-i})\cap\cu_{\mck}(\ce_{-i})).\end{eqnarray*}
\end{prop}

\begin{proof} In all four cases the diffeomorphism is given by
$$(P,Q)\mapsto (\Phi(P)\Phi(Q)^*, \Phi(P)).$$

By Kuiper's Theorem \cite{Kui:HTU} the space $\cu(\ce_i,\ce_{-i})$ is
contractible and hence we obtain the claimed homotopy equivalences.
\end{proof}

The {\em Maslov Index} \cite{CapLeeMil:OMI}, \cite{Nic:MIS} is an integer
invariant of Fredholm pairs of paths of Lagrangian subspaces. We discuss it in
terms of the projection
picture of Lagrangian subspaces. Let $(f,g):[0,1]\to\Gr^{(2)}_{\Fred}(H)$ be a
continuous path
(i.e. a pair of paths in $\Gr(H)$ such that $(f(t),g(t))$ is Fredholm for all
$t$). The Maslov index $\Mas (f,g)$ is the algebraic count of how many times
$\ker f(t)$ passes through $\im g(t)$ along the path. We use the notation
$\Mas(f,g), \Mas(\ker f,\im g), \Mas(\im f,\ker g)$ interchangeably. Indeed
$\Mas(\ker f,\im g)=\Mas(\gamma \ker f,\gamma \im g)=\Mas(\im f,\ker g)$.

The Maslov index has the
following properties (cf. \cite{Nic:MIS}, \cite{CapLeeMil:OMI}):
\begin{enumerate}
\item \textit{Path Additivity:} Let $(f_j,g_j):[0,1]\to
\Gr^{(2)}_{\Fred}(H), j=1,2$,
be continuous paths with $f_2(0)=f_1(1), g_2(0)=g_1(1)$ then
\[ \Mas((f_1,g_1)*(f_2,g_2))=\Mas(f_1,g_1)+\Mas(f_2,g_2).\]
\item \textit{Homotopy Invariance:} Let $(f_j,g_j):[0,1]\to
\Gr^{(2)}_{\Fred}(H)$, $j=0,1$, such that $(f_0,g_0)$ is homotopic
$(f_1,g_1)$ rel
endpoints then
\[ \Mas(f_0,g_0)=\Mas(f_1,g_1).\]

More generally, suppose that $(F,G)$ is a homotopy from
$(f_0,g_0)=(F(-,0),G(-,0))$ to
$(f_1,g_1)=(F(-,1),G(-,1))$ and suppose that $\dim(\ker F(0,s)\cap \im G(0,s))$
and $\dim(\ker F(1,s)\cap \im G(1,s))$ are independent of $s\in [0,1]$. Then
$ \Mas(f_0,g_0)=\Mas(f_1,g_1).$

\item \textit{Normalization:} This is done in two steps. First one requires
that on paths with endpoints in $\Gr^{(2)}_*(H)$ the Maslov index
induces a group isomorphism
$\pi_1(\Gr^{(2)}_{\Fred}(H),\Gr^{(2)}_*(H))\to \Z$
(since $\Gr^{(2)}_*(H)\simeq \cu_*(\ce_{-i})\times \cu(\ce_i,\ce_{-i})$
is contractible $\pi_1(\Gr^{(2)}_{\Fred}(H),\Gr^{(2)}_*(H))$ is
indeed a group).
This determines $\Mas$ on paths with endpoints in $\Gr^{(2)}_*(H)$ up
to a sign.
The sign is chosen as follows: if $(P,Q)\in\Gr^{(2)}(H)$ then
$\Mas(e^{t\gamma }Pe^{-t\gamma },Q)_{-\eps\le t\le \eps}=\dim (\ker
P\cap\im Q)$
for $\eps$ small enough.

Secondly, if $(f,g):[0,1]\to\Gr^{(2)}_{\Fred}(H)$ is an arbitrary continuous
path then one chooses $\eps$ small enough such that the pairs
$(e^{s\gamma }f(j)e^{-s\gamma },g(j))$ are invertible for $j=0,1, 0<s\le\eps$.
Then one puts
\begin{equation} \Mas(f,g):= \Mas(e^{\eps \gamma }fe^{-\eps \gamma
},g).\label{ML-G6.20}\end{equation}
\end{enumerate}

Actually, the normalization property 3. determines the Maslov index completely
and it may be viewed as its definition. Properties 1. and 2. follow
from 3. There
exist other conventions for dealing with paths whose endpoints do not lie in
$\Gr^{(2)}_*(H)$ and not all of these conventions satisfy path additivity.

\vskip.4in

The discussion of the Maslov index works as well when the Hermitian symplectic
Hilbert space  $H$ is finite-dimensional. In this case  the Fredholm condition
is vacuous and the Maslov index is defined on
$\Gr^{(2)}(H)=\Gr(H)\times \Gr(H)$. We
will   use the Maslov index in both contexts; in the
infinite--dimensional setting with
$H=L^2(E_{|\del X})$ and in the finite--dimensional context with $H=\ker A$.

\vskip.4in

\begin{theorem}\label{MaslovWinding} For a continuous path $(f,g)$ in
$\Gr^{(2)}_{\Fred}(H)$ the Maslov index is related to the winding number by the
equation
\begin{equation}
     \Mas(f,g)=-\wind(\Phi(f)\Phi(g)^*). \label{ML-G6.21}
\end{equation}
\end{theorem}
\begin{proof} In view of Proposition \plref{homequ} the right hand side of
\eqref{ML-G6.21} induces a group isomorphism
$\pi_1(\Gr^{(2)}_{\Fred}(H),\Gr^{(2)}_*(H))\to \Z$. It remains therefore to
check the sign convention and the convention for paths with endpoints not in
$\Gr^{(2)}_*(H)$.
Let $(P,Q)\in\Gr_{\Fred}^{(2)}(H)$. Then, by definition, $\Mas(e^{t\gamma
}Pe^{-t\gamma },Q)_{-\eps\le t\le \eps}=
\dim (\ker P\cap\im Q)$ for $\eps$
small enough. A straightforward calculation shows
\begin{equation}
     \Phi(e^{s\gamma }Pe^{-s\gamma })=e^{-2s\gamma }\Phi(P)\label{ML-G6.22}
\end{equation}
and thus for $\eps>0$ small enough we have, in view of Lemma
\plref{pairs} \eqref{grassm3},
\begin{equation}\begin{split}
        \wind(\Phi(e^{s\gamma }Pe^{-s\gamma })\Phi(Q)^*)_{-\eps\le
s\le\eps}&=\dim (\ker
        P\cap\im Q) \wind(-e^{-2is})_{-\eps\le s\le \eps}\\&=-\dim
(\ker P\cap\im Q).
        \end{split}\label{ML-G6.23}
\end{equation}
To check the convention for paths with endpoints not in $\Gr_*^{(2)}(H)$ we
consider the paths $(e^{s\gamma }Pe^{-s\gamma },Q), -\eps\le s\le 0$,
resp. $0\le s\le \eps$.
By definition we have for $\delta>0$ small enough
\begin{equation}\begin{split}
      \Mas(e^{s\gamma }Pe^{-s\gamma },Q)_{-\eps\le s\le 0}&=\Mas(e^{(\delta+s)
\gamma }Pe^{-(\delta+s)
\gamma },Q)_{-\eps\le s\le 0}\\
&=     \Mas(e^{t\gamma }Pe^{-t\gamma },Q)_{-\eps+\delta\le
s\le+\delta}=\dim(\ker P\cap\im Q),\end{split}\label{ML-G6.24}
\end{equation}
and, analogously,
\begin{equation}
     \Mas(e^{s\gamma }Pe^{-s\gamma },Q)_{0\le s\le \eps}=0.\label{ML-G6.25}
\end{equation}
According to our convention for the winding number we have, on the other hand,
\begin{equation}\begin{split}
      \wind(-e^{-2is})_{-\eps\le s\le 0}&=-1,\\
      \wind(-e^{2is})_{0\le s\le \eps}&=0.
        \end{split}\label{ML-G6.26}
\end{equation}
In view of \eqref{ML-G6.23} the proof is complete.\end{proof}

\begin{cor}\label{Maslovorientation} Let $(f,g)$ be a continuous path
in $\Gr_{\Fred}^{(2)}(H)$.
\begin{thmenum}\item
The Maslov index $\Mas_{-\gamma}$ with respect to the opposite
symplectic structure is related to $\Mas_{\gamma}$ as follows:
$\Mas_{-\gamma}(f,g)=\Mas_{\gamma}(g,f)$.

\item $\Mas_\gamma(f,g)+\Mas_{\gamma}(g,f)=\dim(\ker f(1)\cap\im g(1))-\dim(\ker
f(0)\cap\im g(0))$.
\end{thmenum}
\end{cor}
\begin{proof} (1) In view of \eqref{ML-G5.5} and the previous theorem we
find
$\Mas_{-\gamma}(f,g)=-\wind(\Phi_{-\gamma}(f)\Phi_{-\gamma}(g)^*)
=-\wind(\Phi_{\gamma}(f)^*\Phi_{\gamma}(g))=
\wind(\Phi_{\gamma}(g)\Phi_{\gamma}(f)^*)=\Mas_{\gamma}(g,f).$

(2) Using the previous Theorem and Corollary \plref{windfinv} we obtain
(we write $\Mas$ instead of $\Mas_\gamma$)
  \[
\begin{split}
    \Mas(f,g)+\Mas(g,f)&=-\wind(\Phi(f)\Phi(g)^*)-\wind((\Phi(f)\Phi(g)^*)^{-1})\\
       &=\dim(\ker f(1)\cap\im g(1))-\dim(\ker f(0)\cap\im g(0)).
\end{split}\]
\end{proof}

Finally we construct a version of the Maslov triple index in our context.
The Maslov triple index as defined in (cf. \cite[Sec. 8]{CapLeeMil:OMI})
cannot be generalized to the present infinite--dimensional setting. The reason
is simply that interesting triples of Lagrangian subspaces $L_1,L_2,L_3$
such that $(L_1,L_2),(L_2,L_3),(L_3,L_1)$ are all Fredholm pairs are
hard to find.

However, motivated by \cite[Prop. 8.4]{CapLeeMil:OMI}
we can construct a variant $\tau_\mu$ of the Maslov triple index which is
related to the double index $\tau_w$: consider continuous paths
$f,g,h:[0,1]\to\Gr(H)$ such that $(f,g), (g,h), (f,h)$ map into
$\Gr_{\Fred}^{(2)}(H)$ and such that $f-g$ or $g-h$ or $f-h$ maps into the
set of compact operators. If, say, $f(t)-g(t)$ is compact for all $t$ then,
of course, it suffices to assume that $(f(t),h(t))$ is Fredholm for all $t$.
The Fredholmness of $(f(t),g(t)), (g(t),h(t))$ then follows.
Now in view of Theorem \plref{MaslovWinding} and Proposition
\plref{doubleindex} we find
\begin{equation}\label{eq6.21}\begin{split}
      &\Mas(f,g)+\Mas(g,h)-\Mas(f,h)\\
&=-\wind(\Phi(f)\Phi(g)^*)-\wind(\Phi(g)\Phi(h)^*)
                                 +\wind(\Phi(f)\Phi(h)^*)\\
                &=-\tau_w(\Phi(f(1))\Phi(g(1))^*,\Phi(g(1))\Phi(h(1))^*)\\
                &\qquad+ \tau_w(\Phi(f(0))\Phi(g(0))^*,\Phi(g(0))\Phi(h(1))^*).
        \end{split}
\end{equation}

This motivates the following definition.

\begin{dfn}\label{tripleinddef} Let $P,Q,R\in\Gr(H)$ such that
$(P,Q),(Q,R),(P,R)$ are Fredholm and at least one of the differences
$P-Q,Q-R,P-R$ is compact. Then we set
\begin{equation}
     \tau_\mu(P,Q,R):= -\tau_w(\Phi(P)\Phi(Q)^*,\Phi(Q)\Phi(R)^*).
\end{equation}
and call $\tau_\mu$ the {\it triple index} of   $(P,Q,R)$.
\end{dfn}

The triple index $\tau_\mu$ inherits properties from $\tau_w$ in a
straightforward way.  For example, one has the following.

\begin{lemma}\label{tripleindex}
Let $P,Q,R\in\Gr(H)$ such that $P-Q,Q-R$ are trace class. Then
\begin{equation}\begin{split}
   \tau_\mu(P,Q,R)=\frac{1}{2\pi
   i}\bigl(&\tr\log(\Phi(P)\Phi(Q)^*)+\tr\log(\Phi(Q)\Phi(R)^*)\\
&-\tr\log(\Phi(P)\Phi(R)^*)\bigr).
        \end{split}
\end{equation}\qed
\end{lemma}

We will have occasion below to use the homotopy invariance of the triple index.
\begin{lemma}\label{triplehomotopy}
Let $P,Q,R:[0,1]\to \Gr(H)$ be paths in $\Gr(H)$ so that
$(P,Q), (Q,R), (P,R)$ map into $\Gr_{\Fred}^{(2)}(H)$ and at least one of the
differences
$P-Q, Q-R, P-R$ maps into the set of compact operators.
Suppose further that $\dim(\ker P(t)\cap \im Q(t))$, $\dim(\ker Q(t)\cap \im
R(t))$, and $\dim(\ker P(t)\cap \im R(t))$ are independent of $t$. Then
$$ \tau_\mu(P(0),Q(0),R(0))= \tau_\mu(P(1),Q(1),R(1)).$$
\end{lemma}
\begin{proof} By \eqref{eq6.21} we have
\begin{equation}\begin{split}
    \tau_\mu&(P(0),Q(0),R(0))-\tau_\mu(P(1),Q(1),R(1))\\
      &=\Mas(P,Q)+\Mas(Q,R)-\Mas(P,R).
        \end{split}
\end{equation}
Now the claim follows immediately from the homotopy invariance of the Maslov
index.
\end{proof}

We defined the triple index in such a way that formulas become short. A
drawback of this is that $\tau_\mu$ is not antisymmetric in the variables, as
the following proposition shows.
\begin{prop}\label{vanishing} Let $P,Q,R\in\Gr(H)$ such that $(P,Q), (Q,R),
(P,R)$ are Fredholm and at least one of the differences is compact. Then
\begin{equation}\begin{split} \tau_\mu(P,R,Q)=& -\tau_\mu(P,Q,R)+\dim(\ker Q\cap\im R),\\
    \tau_\mu(Q,P,R)=& -\tau_\mu(P,Q,R)+\dim(\ker P\cap\im Q),\\
    \tau_\mu(R,Q,P)=& -\tau_\mu(P,Q,R)+\dim(\ker P\cap\im Q) \\
                   &\quad \quad+\dim(\ker Q\cap\im  R)-\dim(\ker P\cap\im R).
        \end{split}\label{ML-G6.29}
\end{equation}
Moreover,
\begin{equation}
      \tau_\mu(P,P,Q)=\tau_\mu(Q,P,P)=0  \text{ and }  \tau_\mu(P,Q,P)=\dim(\ker
P\cap\im Q).\label{ML-G6.30}
\end{equation}
\end{prop}
\begin{proof} We prove \eqref{ML-G6.30} first. From Proposition
\plref{doubleindex}
and the definition of $\tau_\mu$ we infer
\begin{equation}\begin{split}
     \tau_\mu(P,P,Q)&=-\tau_w(I,\Phi(P)\Phi(Q)^*)=0,\\
     \tau_\mu(Q,P,P)&=-\tau_w(\Phi(Q)\Phi(P)^*,I)=0,\\
     \tau_\mu(P,Q,P)&=-\tau_w(\Phi(P)\Phi(Q)^*,(\Phi(P)\Phi(Q)^*)^{-1})=\dim(\ker
     P\cap\im Q).
        \end{split}
\end{equation}
To prove \eqref{ML-G6.29} we assume, without loss of generality,  that $Q-R$ is
compact. Let
$f(t):=(1-t)Q+t R, 0\le t\le 1$. Then we obtain from Corollary
\plref{Maslovorientation} and \eqref{ML-G6.30}
\begin{equation}\begin{split}
    \tau_\mu(P,R,Q)&=\tau_\mu(P,R,Q)-\tau_\mu(P,Q,Q)=\Mas(P,f)+\Mas(f,Q)\\
            &=\Mas(P,f)-\Mas(Q,f)+\dim(\ker Q\cap\im R)\\
            &= -\tau_\mu(P,Q,R)+\dim(\ker Q\cap\im R),\\
    \tau_\mu(Q,P,R)&=\Mas(P,f)-\Mas(Q,f)+\tau_\mu(Q,P,Q)\\
                   &=-\tau_\mu(P,Q,R)+\dim(\ker P\cap\im Q),\\
    \tau_\mu(R,Q,P)&=\Mas(f,Q)-\Mas(f,P)\\
     &=-\Mas(Q,f)+\Mas(P,f)+\dim(\ker P\cap\im Q) \\
                   &\quad \quad+\dim(\ker Q\cap\im  R)-\dim(\ker P\cap\im R).
        \end{split}
\end{equation}
\end{proof}

\subsection{Symplectic reduction}\label{Symplectic reduction}

We conclude this section with  a discussion of symplectic reduction in our
infinite--dimensional context. We will use symplectic reduction in Section
\plref{sec8}.

\newcommand{\Ann}{\operatorname{Ann}}
Let $(H,\scalar{.}{.},\gamma)$ be a Hermitian symplectic Hilbert space
with symplectic
form $\go(x,y)=\scalar{x}{\gamma y}.$ For a subspace $U\subset H$ the {\it
annihilator of $U$} is defined to be
\[\Ann(U):=\bigsetdef{x\in
H}{\forall y\in U\; \go(x,y)=0}=(\gamma U)^\perp.\] A subspace $U\subset H$ is
called {\it isotropic} if $U\subset \Ann(U)$.

Assume for the moment that
$H$ is finite--dimensional and that $\Ann(U)\subset U$. Then $\go$ induces a
symplectic structure on the quotient $U/\Ann(U)$ in a natural way. Moreover,
if $L\subset H$ is Lagrangian then $R_U(L):=L\cap U/L\cap\Ann(U)$ is Lagrangian
in $U/\Ann(U)$. $R_U(L)$ is called the {\it symplectic reduction} of $L$.

\begin{prop} \label{symplreduction} Let $(H,\scalar{.}{.},\gamma)$ be
a Hermitian symplectic Hilbert space, $U\subset H$ a closed subspace with
$\Ann(U)\subset U$.

Suppose that $L\subset H$ is a Lagrangian subspace such that $L+\Ann(U)$ is
a closed subspace of $H$. Then $(U\cap\gamma U,\scalar{.}{.},\gamma)$ is a Hermitian
symplectic Hilbert space and the orthogonal projection
\[
   P_{L,U}:=\operatorname{proj}_{U\cap\gamma U}:L\cap U\longrightarrow
U\cap\gamma
   U\]
has closed range isomorphic to
$L\cap U/L\cap\Ann(U)$. Moreover,
$R_U(L)=\im P_{L,U}$ is
   Lagrangian
in $U\cap \gamma U$.

$R_U(L)$ is called the {\em symplectic reduction  of $L$ with
respect to } $U$.
\end{prop}
\begin{remark}\hfill
\begin{enumerate}
\item  For $(U\cap\gamma U,\scalar{.}{.},\gamma)$ to be Hermitian
symplectic it is crucial that there is at least one Lagrangian  subspace
$L\subset H$ with $L+\Ann(U)$ closed. To illustrate the problem start with
an infinite--dimensional Hermitian symplectic Hilbert space
$(H,\scalar{.}{.},\gamma)$. Let $\tilde H:=H\oplus H_1$, where $H_1$ is another
Hilbert space, and put $\tilde\gamma:=\gamma\oplus i$. Furthermore, pick a
symplectic subspace $L\subset H$
and put $U:=L\oplus H_1\subset \tilde H$. Then $\Ann(U)=L\oplus 0$ and $U\cap\gamma
U=0\oplus H_1$. Since $\tilde \gamma$ acts by multiplication by $i$ on
$U\cap\gamma U$ we conclude that $(U\cap\gamma U,\scalar{.}{.},\tilde \gamma)$ is
not Hermitian symplectic. From the proposition we infer that for each
Lagrangian subspace $K\subset \tilde H$ the space $K+\Ann(U)$ is not
closed.

\item Proposition \plref{symplreduction} in particular applies if $(L,\Ann(U))$ is
a Fredholm pair of subspaces.

\item The assignment $L\mapsto R_U(L)$ is not continuous,
but is continuous along paths
$L_t$ so that $\dim (L_t\cap \Ann(U))$ is constant. These facts are
well--known and we omit the examples.
\end{enumerate}
\end{remark}
\begin{proof} Certainly $U\cap\gamma U$ is a Hilbert space and $\gamma$
leaves $U\cap \gamma U$ invariant. If we can prove that $\im P_{L,U}$ is
Lagrangian in $U\cap\gamma U$ then from Lemma \plref{hermitiansymplectic} we
infer
$\dim(\ker(\gamma+i)\cap U\cap\gamma U)=\dim(\ker(\gamma-i)\cap U\cap \gamma
U)$.

What remains, therefore, is to prove the second part of   Proposition
\plref{symplreduction} without using the fact that
$\dim(\ker(\gamma+i)\cap U\cap\gamma U)=\dim(\ker(\gamma-i)\cap U\cap \gamma
U)$.

We note first that we have an orthogonal direct sum decomposition
\begin{equation}
\Ann(U)\oplus (U\cap\gamma U)=U.
\label{ML-G6.31}
\end{equation}

Also, $\im P_{L,U}$ is an isotropic subspace of $U\cap\gamma U$. In fact, if
$x\in L\cap U$
then, since $L$ is Lagrangian, $\scalar{x}{\gamma x}=0$. Writing $x=\xi+\eta,
\xi\in U\cap\gamma U, \eta\in \Ann(U)$ then
$0=\scalar{x}{\gamma x}=\scalar{\xi}{\gamma
x}=\scalar{\xi}{\xi}=\scalar{P_{L,U}x}{\gamma P_{L,U}x}$.

Next consider $\xi\in U\cap\gamma U$ such that $\gamma(\xi)\perp \im P_{L,U}$.
Thus for all $x\in L\cap U$ we have
$\scalar{\gamma(\xi)}{x}=\scalar{\gamma(\xi)}{P_{L,U}x}=0$.
Hence $\gamma(\xi)\in (L\cap
U)^{\perp}=\ovl{L^\perp+U^\perp}=\ovl{\gamma(L)+U^\perp}$
and consequently, since $L+\Ann(U)$ is closed, $\xi\in L+\Ann(U)$.
We may write $\xi=l+\eta, l\in L, \eta\in\Ann(U)$. From $\xi\in U,
\eta\in\Ann(U)\subset U$ we infer $l\in L\cap U$ and hence $\xi=P_{L,U}(l)\in\im
P(L)$.

Summing up we have proved $\gamma((\im P_{L,U})^\perp)\subset \im P_{L,U}$. Since
$\im P_{L,U}$ is isotropic this implies $\im P_{L,U}=\gamma((\im
P_{L,U})^\perp)$. Thus
$\im P_{L,U}$ is a Lagrangian (in particular closed) subspace of $U\cap\gamma U$.

From \eqref{ML-G6.31} it is now clear that $\im P_{L,U}$ is isomorphic to $L\cap
U/L\cap \Ann(U)$.\end{proof}

\section{Splittings of manifolds and the $\eta$--invariant II}\label{sec7}

For the proof of Theorem \plref{symsplit},  Lemma \plref{swconc} was
crucial. The proof of Lemma \plref{swconc} depends on the Scott--Wojciechowski
theorem \plref{ML-S4.1}. In this section we want to give   proofs
of Lemma \plref{swconc} and Theorem \plref{symsplit} which are independent
of the Scott--Wojciechowski theorem and which apply  to all $P\in\Gr(A)$.
We only use (a mild generalization of) Theorem \plref{brule}.
Moreover, we derive generalizations of two results due to L. Nicolaescu
\cite{Nic:MIS}. This in turn   leads to a nicer version of the splitting
formula for the $\eta$--invariant which involves our version of the Maslov
triple index.

We first introduce a setting which slightly generalizes the one
described in Section \plref{sec2}. Let $X$ be a compact Riemannian manifold
with boundary $\partial X=Y\coprod Z$, i.e. the boundary is a disjoint
union of two (not necessarily connected) manifolds. We assume that
in collars $U=U_Y$ and $U_Z$ we have $D=\gamma_Y(\frac{d}{dx}+A_Y)$
(resp. $D=\gamma_Z(\frac{d}{dx}+A_Z)$) and that the $\pm i$--eigenspaces
of $\gamma_Y$ ($\gamma_Z$) acting on $\ker A_Y$ ($\ker A_Z$) have the
same dimension. The latter does not follow as in section \plref{sec2};
rather it is an assumption. We fix once and for all a $P_Z\in\Gr(A_Z)$.
Then we can define the \Calderon\ projector (relative to $P_Z$)
in $\Gr(A_Y)$. Write again $A$ instead of $A_Y$.
It will be convenient to address $Y, Z$ as boundary components
although $Y, Z$ are not assumed to be connected.

The results of Sections \plref{sec2} to \plref{sec5} generalize verbatim
to this more general setting. Also Theorem \plref{brule} applies to this
setting since all proofs work locally in a collar of the separating
hypersurface.
The advantage of this setting is that it allows in particular to glue
cylinders of the form $[0,\eps]\times N$ with different boundary
conditions on the ends to a manifold.

\begin{lemma}[{cf. \cite[Lemma 2.5]{LesWoj:IGA}}]\label{ML-S6.1} Let
$M=[0,\eps]\times N$
and $D=\gamma(\frac{d}{dx}+A)$ as before. Moreover, let $P,Q\in\Gr(A)$
and denote by $D_{P,Q}$ be the operator obtained by imposing the
boundary condition $P$ at $\{0\}\times N$ and $I-Q$ at $\{\eps\}\times N$.
Then $\gl\in\spec D_{P,Q}$ if and only if $-\gl\in\spec D_{Q,P}$.
In particular,
\[ \eta(D_{P,Q})=-\eta(D_{Q,P}),\quad \dim\ker D_{P,Q}=\dim\ker D_{Q,P}.\]
\end{lemma}
\begin{proof} The proof is exactly the same as the proof of \cite[Lemma
2.5]{LesWoj:IGA}.
Namely, the isometry
\[ T:L^2([0,\eps],L^2(E_{|N}))\longrightarrow L^2([0,\eps],L^2(E_{|N})), \quad
Tf(x):=\gamma f(\eps-x)\]
maps the domain of $D_{P,Q}$ onto the domain of $D_{Q,P}$ and it anticommutes
with $D$. Hence $T^*D_{P,Q}T=-D_{Q,P}$ and we are done.
\end{proof}

Now let $M$ be a Riemannian manifold with boundary containing a separating
hypersurface $N\subset (M\setminus \partial M)$. Let $D$ be a Dirac operator
as in Section \plref{sec5}; i.e. in a collar neighborhood
$[-\eps,\eps]\times N$ of
$N$,
$D$ has the form $D=\gamma(\frac{d}{dx}+A)$ as in \eqref{ML-G2.1}.
Moreover, we assume
that the $\pm i$--eigenspaces of $\gamma$ acting on $\ker A$ have the same
dimension. Define $\mcut, M^\pm$ as in Section \plref{sec5}. We assume that
on the boundary components of $(\partial M^\pm)\setminus N$
self--adjoint boundary
projections have been fixed once and for all.

\begin{lemma}\label{ML-S6.2}
For any $P\in\Gr(A)$ we have
$\etab(D,M)-\etab(D_P,M^+)-\etab(D_{I-P},M^-)\in \Z$.
\end{lemma}
\begin{proof} Denote by $\mcut_\eps$ the manifold with boundary obtained
by removing $[-\eps,\eps]\times N$ from $M$. As in Lemma \plref{ML-S6.1}
for $P,Q\in\Gr(A)$ we denote by $\eta(D_{P,Q},[-\eps,\eps]\times N)$
the $\eta$--invariant of the operator on $[-\eps,\eps]\times N$ obtained
from $D$ by imposing the boundary condition $P$ at $\{-\eps\}\times N$
and the boundary condition $I-Q$ at $\{\eps\}\times N$.
The $\mod \Z$ gluing formula for the $\eta$--invariant
\eqref{p+eqn} then implies
\begin{equation}
\etab(D,M)\equiv \etab(D_{P^+\oplus
(I-P^+)},\mcut_\eps)+\etab(D_{P^+,P^+},[-\eps,\eps]\times N)\mod \Z
\label{ML-G6.1}
\end{equation}
for $P^+=P^+(L)$ the Atiyah--Patodi--Singer projection with respect to a
Lagrangian subspace $L\subset \ker A$.  One easily checks that $\ker
D_{P^+,P^+}=\{0\},$ hence Lemma
\plref{ML-S6.1} implies
\begin{equation}
     \etab(D_{P^+,P^+},[-\eps,\eps]\times N)=0.
\end{equation}
Also by Lemma \plref{ML-S6.1}
\begin{equation}\begin{split}
    \etab(D_{P^+,P}&,[-\eps,0]\times N)+\etab(D_{P,P^+},[0,\eps]\times N)\\
    &=\frac 12\dim\ker(D_{P^+,P},[-\eps,0]\times N)+\frac 12
    \dim\ker(D_{P,P^+},[0,\eps]\times N)\\
    & =  \dim\ker (D_{P,P^+},[0,\eps]\times N)\in\Z.
        \end{split}
\end{equation}
Plugging this into \eqref{ML-G6.1} and applying again the $\mod \Z$
splitting formula
for the $\eta$--invariant we get
\begin{equation}\begin{split}
   \etab(D,M)&\equiv  \etab(D_{P^+\oplus (I-P^+)},\mcut_\eps) +
   \etab(D_{P^+,P},[-\eps,0]\times N)+\etab(D_{P,P^+},[0,\eps]\times N)\\
    &\equiv \etab(D_P,M^+)+\etab(D_{I-P},M^-)\mod\Z.
        \end{split}
\end{equation}
\end{proof}

\begin{lemma}\label{swconc1} Lemma \plref{swconc} holds for all
$P_0,P_1\in\Gr(A)$.
\end{lemma}
\begin{proof} We freely use the notations of Lemma \plref{swconc} and its
   proof.
By    Lemma \plref{ML-S6.2} we have for all $t$
\begin{equation}\begin{split}
  &\etab(D_{Q_t},\mcut)-\etab(D_{Q_0},\mcut)\\
  &=(\etab(D_{Q_t},\mcut)-\etab(D,M))-(\etab(D_{Q_0},\mcut)-\etab(D,M))\in\Z.
        \end{split}
\end{equation}
Hence
\begin{equation}
   \frac{d}{dt}\etab(D_{Q_t},\mcut)=0
\end{equation}
and the assertion follows from Lemma \plref{sf-eta} and Lemma \plref{splitbundles}.
\end{proof}

Now we can prove the following considerable generalization of the
splitting formula for the $\eta$--invariant. In Theorem
\plref{symsplit} we assumed
that $P\in \Gr_\infty(A)$. In the following theorem we only
require $P\in
\Gr(A)$.


\begin{theorem}\label{symsplit1}
The statement of Theorem \plref{symsplit} remains valid if $P\in\Gr(A)$ and
$P_t$ is a smooth path in $\Gr(A)$ from $P$ to the \Calderon\ projector.
\end{theorem}
\begin{proof} The proof is exactly the same as the one of Theorem
\plref{symsplit}. One only has to invoke Lemma \plref{swconc1} instead
of Lemma \plref{swconc}.
\end{proof}

\vskip.4in

We next present generalizations of two results
due to L. Nicolaescu \cite{Nic:MIS}.

\begin{theorem}\label{ML-S6.5} Let $X$ be a manifold with boundary  and $D(t),
a\le t\le b,$ a smooth family of Dirac operators. We assume that in a collar of
the boundary
$D$ takes the form $\gamma(\frac{d}{dx}+A(t))$ as before. Let
$P(t)\in\Gr(A(t))$
be a smooth family. Denote by $P_X(t)$ the \Calderon\ projectors of $D(t)$, and
$L_X(t)=\im P_X(t)$ the Cauchy data spaces. Then
\[ \SF(D_{P(t)}(t))_{t\in [a,b]}=  \Mas(P(t),P_X(t))_{t\in[a,b]}= \Mas(\ker
P(t), L_X(t))_{t\in [a,b]}.\]
\end{theorem}
Note that $\gamma$ is assumed to be constant. This is essential.
Note that in \cite[Theorem 4.3]{DanKir:GSF} the collar of $\partial X$ was
parametrized as $(-\eps,0]\times\partial X$. Their formula is obtained by
invoking
Corollary \plref{Maslovorientation}.
\begin{proof}  We first consider the case $P(t)\in\Gr_{\infty}(A(t))$.
Since $D_{P_X(t)}(t)$ is invertible, its spectral flow vanishes. We apply
Lemma \plref{sf-eta}, Theorem \plref{invertible}, \eqref{ML-G6.10},  and Theorem
\plref{MaslovWinding} to calculate
\[\begin{split}
     \SF(D_{P(t)}(t))_{t\in [a,b]}&= \SF(D_{P(t)}(t))_{t\in [a,b]}-
     \SF(D_{P_X(t)}(t))_{t\in [a,b]}\\
     &=
\etab(D_{P(b)}(b))-\etab(D_{P_X(b)}(b))-\etab(D_{P(a)}(a))+\etab(D_{P_X(a)}(a))\\
     &\quad
 -\int_a^b\frac{d}{dt}\bigl(\etab(D_{P(t)}(t))-\etab(D_{P_X(t)}(t))\bigr)dt\\
     &=\tfrac{1}{2\pi i}\tr\log(\Phi(P(b))\Phi(P_X(b)^*))-\tfrac{1}{2\pi
i}\tr\log(
     \Phi(P(a))\Phi(P_X(a)^*))\\
     &\quad -\int_a^b\tfrac{1}{2\pi i} \frac{d}{dt}\tr\log
     (\Phi(P(t))\Phi(P_X(t))^*)dt\\
     &= -\wind(\Phi(P(t))\Phi(P_X(t))^*)_{t\in[a,b]}=\Mas(P(t),P_X(t))_{t\in[a,b]}.
   \end{split}
\]

  Now suppose that   $P(t)$ is arbitrary. Choose smooth paths $P_0(t)$ in
$\Gr(A(0))$
    and $P_1(t)\in\Gr(A(1))$ such that $P_0(0)\in\Gr_\infty(A(0)), P_0(1)=P(0),
    P_1(0)=P(1), P_1(1)\in\Gr_\infty(A(1))$
and such that
\begin{equation}
\dim (\ker P_0(t)\cap\im P_X(0)) \text{ and }  \dim (\ker P_1(t)\cap\im
P_X(1))\label{ML-G7.7}
\end{equation}
    are independent of $t$. The existence of $P_0, P_1$ follows from Lemma
    \plref{ML2-S6.1} by considering $\Phi(P(j))\Phi(P_X(j))^*, j=0,1$.
In view of \eqref{ML-G7.7} and Lemma \plref{kerdp} the  dimension of
the kernels
of
$D_{P_0(t)}(0)$ and $ D_{P_1(t)}(1)$ are constant and hence the spectral flow
of $D_{P_0(t)}(0) $ and $D_{P_1(t)}(1)$ vanishes. We may therefore
compose the paths
$D_{P_0(t)}(0), D_{P(t)}(t), D_{P_1}(t)$ without changing the spectral flow.
Also $\Mas(P_0(t),P_X(0))=\Mas(P_1(t),P_X(1))=0$ in view of
    \eqref{ML-G7.7}.
In sum, without loss of generality  we may assume that the family $P(t)$ satisfies
    $P(0)\in\Gr_\infty(A(0)), P(1)\in\Gr_{\infty}(A(1))$.
Now consider the path $\Phi(P(t))\Phi(P_X(t))^*$ in $\cu_{\Fred}$. In view of
Lemma
    \plref{ML2-S6.1} this path is homotopic rel endpoints
to a path $f(t)\in\cu_\infty$. Putting $\tilde
    P(t):=\Phi^{-1}(f(t)\Phi(P_X(t)))\in\Gr_{\infty}(A(t))$
we see that $(P(t),P_X(t))$ is homotopic rel endpoints
to the path $(\tilde P(t),P_X(t))$. Since homotopies with fixed endpoints
neither change the spectral flow nor the Maslov index we find
\[\begin{split}
\SF(D_{P(t)}(t))_{t\in[a,b]}&=\SF(D_{\tilde P(t)}(t))_{t\in[a,b]}=\Mas(\tilde
    P(t),P_X(t))_{t\in[a,b]}\\
&=\Mas(P(t),P_X(t))_{t\in[a,b]}.\end{split}\]
\end{proof}

We also give a generalization of Nicolaescu's theorem for closed manifolds. The
result in the following form was first proven in \cite{Dan:ETN}.

\begin{theorem}\label{ML-S6.6} Let $M$ be a split manifold as in Section
\plref{sec5} and let
$D(t), a\le t\le b,$ be a smooth path of Dirac type operators such that in a
collar of the separating hypersurface we have $D(t)=\gamma(\frac{d}{dx}+A(t))$.
Then
\[   \SF(D(t))_{t\in[a,b]}=\Mas_{\gamma}(P_{M^-}(t),I-P_{M^+}(t))_{t\in[a,b]}=
\Mas(L_{M^-}(t),L_{M^+}(t))_{t\in[a,b]} .\]
\end{theorem}
\begin{proof}   Corollary \plref{ML-S5.6} states that
\begin{equation}\label{ML-G812}
\SF(D(t))_{t\in[a,b]}=\SF(D_{I-P_{M^+}(t)}(t),M^-)_{t\in[a,b]}.\end{equation}

Applying Theorem \plref{ML-S6.5} to the right hand side of
\eqref{ML-G812} and using Corollary \plref{Maslovorientation} yields
\[\begin{split}
\SF(D_{I-P_{M^+}(t)}(t),M^-)_{t\in[a,b]}&=
\Mas_{-\gamma}(I-P_{M^+}(t),P_{M^-}(t))_{t\in[a,b]}\\ &=
\Mas_{\gamma}(P_{M^-}(t),I-P_{M^+}(t))_{t\in[a,b]},
\end{split}\]
finishing the proof.
\end{proof}

Notice that the proof of Theorem \plref{ML-S6.6} does not rely on Theorem
\plref{symsplit}, and in particular does not use the result  of
\cite{BruLes:EIN}.

\vskip.4in

Finally, we state the following nicer version of the gluing formula for the
$\eta$--invariant.
We emphasize that the term
$\tau_\mu(I-P_{M_-},P,P_{M^+})$, which was defined in Subsection \plref{secsixtwo}, is
an integer invariant which is defined completely in terms of the Hermitian
symplectic structure on
$L^2(E_{|N})$.

\begin{theorem}\label{symsplit2} In the situation of Theorem \plref{symsplit},
let
$P\in\Gr(A)$.
Then
\[\etab(D,M)=\etab(D_{P},M^+)+\etab(D_{I-P},M^-)-\tau_\mu(I-P_{M_-},P,P_{M^+}).\]
\end{theorem}


\begin{proof}  We note again that $I-P_{M^-}-P_{M^+}$ is trace class
(cf. the proof of Theorem \plref{invertiblecase}). Thus $I-P_{M^-}-P_{M^+}$
is compact and hence the triple index $\tau_\mu(I-P_{M^-},P,P_{M^+})$
is well--defined for any $P\in\Gr(A)$.

Let $P_t, 0\le t\le 1,$ be a smooth path in $\Gr(A)$
from $P$ to the \Calderon\
projector
$P_{M^+}$. Notice that $\Mas_{\gamma}(I-P_{M^-},P_{M^+})=0$ since $I-P_{M^-}$
and $P_{M^+}$ are constant paths.  From Theorem
\plref{symsplit1}, Theorem
\plref{ML-S6.5},
\eqref{eq6.21}, Corollary
\plref{Maslovorientation}, and Proposition \plref{vanishing} we infer
\[\begin{split}
      &\etab(D,M)-\etab(D_{P},M^+)-\etab(D_{I-P},M^-)\\&=
  \SF(D_{P_t},M^+)_{t\in[0,1]}+\SF(D_{I-P_t},M^-)_{t\in[0,1]}\\
 &=\Mas_{\gamma}(P_t,P_{M^+})_{t\in[0,1]}+
\Mas_{-\gamma}(I-P_t,P_{M^-})_{t\in[0,1]}\\
 &=\Mas_\gamma(P_t,P_{M^+})_{t\in[0,1]}+\Mas_{\gamma}(I-P_{M^-},P_t)_{t\in[0,1]}
-\Mas_{\gamma}(I-P_{M^-},P_{M^+})_{t\in[0,1]}\\
       &= \tau_\mu(I-P_{M^-},P_{M^+},P_{M^+})-\tau_\mu(I-P_{M_-},P,P_{M^+})\\
       &=-\tau_\mu(I-P_{M_-},P,P_{M^+}).
   \end{split}
\]\end{proof}

\newcommand\ep{\varepsilon}
\newcommand\mapright[1]{\smash{ \mathop{\longrightarrow}\limits^{#1}}}




\section{Adiabatic stretching and applications to the Atiyah-Patodi-Singer
$\rho$-invariant}
\label{sec8}

For the purpose of computation, one
weakness of the splitting formulas of Theorems
\plref{symsplit}, \plref{invertiblecase}, and \plref{symsplit2}
is that it is difficult in practice to identify the
   Calder\'on projector. In many applications it is more convenient to work
with the Atiyah--Patodi--Singer projection
$P^+(L)$, or at least some finite rank perturbation of $P^+(L)$, as a
boundary condition.  According to Theorem \plref{symsplit}, this requires
knowing the spectral flow of $D_{P_t}$  and $D_{I-P_t}$ along a path
$P_t$ starting at the Calder{\'o}n projector and ending at $P^+(L)$.

A natural choice of such a path is the path obtained by stretching  the
collar neighborhood of the separating surface. According to a theorem of
Nicolaescu \cite{Nic:MIS}, the Calder{\'o}n projector limits to a projection
    of the form $ P_{>\nu}+ \proj_L$, where
$P_{>\nu}$ is the projection to the span of the eigenvectors of $A$ with
eigenvalues greater than $\nu$ and $L$ is a Lagrangian subspace of the
(finite--dimensional) span of eigenvectors of $A$ with eigenvalues in the
range
$[-\nu,\nu]$. The number $\nu$ is the {\em non--resonance level}
\cite{Nic:MIS} of
$D$ acting on $M^+$ and in particular is zero if  and only if there are no
$L^2$ solutions to
$D\phi=0$ on the manifold obtained from $M^+$ by adding an infinite
collar. If $\nu=0$, then the limit of the Calder{\'o}n projector is an
Atiyah--Patodi--Singer projection $P^+(V)$ for a particular Lagrangian
$V\subset \ker A$.

This approach works particularly well to study the odd
signature operator and the Atiyah--Patodi--Singer
$\rho_\alpha$ invariant \cite{AtiPatSin:SARII}, since the effect of the
Riemannian metric is minimized in this important case. We present the details.
The approach can be made to work    for arbitrary Dirac operators as well,
however additional correction terms appear
  corresponding to the 1--parameter  family of operators acting on $M$ and
$M^\pm$ as the collar of the separating hypersurface is stretched to
infinity. We will make some comments about the  case of general Dirac operators
at the end of this section.
\vskip.4in

\subsection{The odd signature operator.}
Let $X$ be a compact manifold of dimension $2n+1$, with (possibly empty)
boundary $\del X^{2n}$. Assume a collar of $\partial X$ is isometric to
$[0,\ep)\times \del X$.  Let
$\alpha:\pi_1(X)\to U(n)$ be  a representation.  To $\alpha$ one can
assign a flat vector bundle, that is, a $\C^n$ bundle $E\to X$ together
with a flat connection
$B$ on
$E$ so that the holonomy representation of $B$ is equal to $\alpha$.  If
$\partial X$ is non--empty, we may assume, by gauge transforming $B$ if
necessary, that
$B$ is in temporal gauge on the collar, in other words   there is a flat
$U(n)$ connection
$b$ on $E_{|\del X}$ so that the restriction of $B$ to the collar
$ [0,\ep)\times \del X$ is of the form
$$B_{[0,\ep)\times \del X}=q^*(b),$$
where $q:[0,\ep)\times \del X\to \del X$ is the projection to the second factor.

Let $d_B:\Omega^p(X;E)\to \Omega^{p+1}(X;E)$ and
$d_b:\Omega^p(\del X;E_{|\del X})\to \Omega^{p+1}(\del X;E_{|\del X})$ denote the associated
coupled DeRham operators. Note that $d_B^2$ and $d_b^2$ are zero since
$B$ and $b$ are flat.  The cohomology of the  complex  $(\Omega^*(X;E),d_B)$
(resp. $(\Omega^*(\del X;E_{|\del X}),d_b)$) is identified via the DeRham theorem with
the singular cohomology $H^*(X;\C^n_\alpha)$ (resp. $H^*(\del X;\C^n_\alpha)$),
where
$\C^n_\alpha$ denotes the local coefficient system determined by the
representation $\alpha$.

The {\it odd signature operator on $X$ coupled to the flat connection $B$}
is the operator
$$D_B:\moplus_{p}\Omega^{2p}(X;E)\to \moplus_{p}\Omega^{2p}(X;E)$$
defined by
$$D_B(\beta)=i^{n+1}(-1)^{p-1}(*d_B-d_B*)(\beta)\ \ \text{for} \ \
\beta\in \Omega^{2p}(X;E),$$
where $*:\Omega^{k}(X;E)\to \Omega^{2n+1-k}(X;E)$ denotes the Hodge $*$
operator (see \cite{AtiPatSin:SARII}).

The operator $D_B$ is a symmetric Dirac operator.  Its square is the
twisted Laplacian acting on even bundle--valued forms:
$$D_B^2=d_B^*d_B+d_Bd_B^*.$$
   In particular $D_B$  is self--adjoint if $X$ has
empty boundary and in that case its kernel can be identified with the
twisted cohomology group
$\oplus_p H^{2p}(X;\C_\alpha^n)$ by the Hodge and DeRham theorems. This
implies that the dimension of the kernel of $D_B$ is independent of the
choice of Riemannian metric   if $X$ is closed.

Define a restriction map
$$r:\moplus_{p}\Omega^{2p}(X;E)\to \moplus_{k}\Omega^{k}(\del X; E_{|\del X})$$
by the formula
$$r(\beta)= i^*(\beta) + i^*(*\beta)$$
where $i:\del X\hookrightarrow X$ denotes the inclusion of the boundary.

   To avoid confusion we denote the Hodge $*$ operator on the boundary by
$\hat{*}$, thus
$$\hat{*}:\Omega^k(\del X; E_{|\del X})\to\Omega^{2n-k}(\del X; E_{|\del X}).$$
We use $\hat{*}$ to define
$$\gamma:\moplus_{k}\Omega^{k}(\del X; E_{|\del X})\to \moplus_{k}\Omega^{k}(\del X;
E_{|\del X})$$ by
$$ \gamma(\beta)=\begin{cases}
i^{n+1}(-1)^{p-1}\hat{*}\ \beta  & \text{if } \beta\in  \Omega^{2p}(\del X;
E_{|\del X}), \\
i^{n+1}(-1)^{n-q}\hat{*}\ \beta & \text{if }  \beta\in  \Omega^{2q+1}(\del X;
E_{|\del X}). \end{cases}$$

Finally, we define the operator $$A_b:\moplus_k\Omega^{k}(\del X;E_{|\del X})\to
\moplus_k\Omega^{k}(\del X;E_{|\del X})$$
by $$A_b(\beta)=
\begin{cases} -(d_b\hat{*}+\hat{*}d_b)\beta & \text{if } \beta\in
\moplus_k\Omega^{2k}(\del X;E_{|\del X}),\\
\quad (d_b\hat{*}+\hat{*}d_b)\beta & \text{if } \beta\in
\moplus_k\Omega^{2k+1}(\del X;E_{|\del X}).\end{cases}$$

The following facts are routine to verify.

\begin{enumerate}
\item $A_b$ is a self--adjoint Dirac operator on $\del X$.
\item $r$ induces an identification
$\Phi:\moplus_p\Omega^{2p}([0,\eps)\times\partial X;E)\to
\cinf{[0,\eps),\moplus_k\Omega^k(\partial X;E)}$ which is isometric with
respect to the $L^2$--structures. Moreover,
\begin{equation}\label{cylinder}
\Phi D_B\Phi^*=\gamma(\tfrac{\del}{\del x} + A_b),\end{equation}
where $x$ denotes the collar coordinate.
\item $\gamma A_b=-A_b\gamma$.
\item $\gamma^2=-I$.
\item $A_b$ reverses the parity of forms.
\item $A_b d_b=-d_b A_b$ and $A_b d_b^*=-d_b^* A_b$, where
$d_b^*=-\hat{*}d_b\hat{*}$ is the $L^2$--adjoint of $d_b$.
\item  $A_b^2$ preserves the subspace $\Omega^k(\del X;E_{|\del X})$ for each $k$,
and equals the twisted Laplacian on $k$--forms, $A_b^2=\Delta_b=
d_bd_b^*+d_b^*d_b$.
\item The kernel of $A_b$ equals $\ker A_b^2=\ker \Delta_b$, which is
identified using the Hodge theorem with the DeRham cohomology of the complex
$(\Omega^k(\del X;E_{|\del X}), d_b)$. The DeRham isomorphism identifies the DeRham
cohomology with twisted
cohomology
$  H^*(\del X;\C^n_\alpha)$, where $\alpha:\pi_1\del X\to U(n)\subset GL(\C^n)$ is the
holonomy representation of the flat connection $b$.
\end{enumerate}

The first 5 facts do not depend on $B$ being a flat connection, and hold
for any $U(n)$ connection in temporal gauge near the boundary. The last
three  depend on $b$ being flat.

For convenience we simplify the notation as follows. Let
$\Omega^{\text{\tiny even}}_X$ denote $\oplus_{p}\Omega^{2p}(X;E)$ and let
$\Omega^*_{\del X}$ denote $
\oplus_{k}\Omega^{k}(\del X;E_{|\del X})$. The $L^2$ completion  of
$\Omega^*_{\del X}$
will be denoted by $L^2(\Omega^*_{\del X})$.   We will often drop the
subscripts ``$B$'' and ``$b$''  and, for example, write $D$ for $D_B$,
$A$ for
$A_b$, and  $d$ for $d_B$ or $d_b$.

The self--adjoint operator $A$ induces a spectral decomposition of
$L^2(\Omega^*_{\del X})$. We denote the $\mu$--eigenspace of $A$ by $E_\mu$.
Given $\nu\ge 0$ we will also use the notation
$$F^+_\nu=\rmspan_{L^2}\{\psi_\mu\ | \ A\psi_\mu=\mu\psi_\mu, \ \mu>\nu\}
=\moplus_{\mu>\nu}E_\mu,$$
$$F^-_\nu=\rmspan_{L^2}\{\psi_\mu\ | \ A\psi_\mu=\mu\psi_\mu, \ \mu<-\nu\}
=\moplus_{\mu<-\nu}E_\mu,$$
$$E^+_\nu=\moplus_{0<\mu\leq\nu}E_\mu,$$
and
$$E^-_\nu=\moplus_{-\nu\leq\mu<0}E_\mu.$$

Thus $E^-_\nu$ is the finite--dimensional span of the eigenvectors of $A$
with eigenvalues $\mu$ in the range $-\nu\leq\mu<0$, $E^+_\nu$
corresponds to the range $0<\mu\leq\nu$ (if $\nu=0$, then $E^\pm_\nu=0$).
Similarly
$F^-_\nu$ is the infinite--dimensional space spanned by eigenvectors with
eigenvalues
$\mu$ satisfying $\mu<-\nu$, and $F^+_\nu$ corresponds to $\mu>\nu$.
In particular $F_0^+$ denotes the positive eigenspan and $F_0^-$ the
negative eigenspan of $A$.
This
gives an orthogonal  decomposition
\begin{equation}\label{eigenspace}L^2(\Omega^*_{\del X})=
F^-_\nu\oplus E_\nu^-\oplus \ker A\oplus E_\nu^+\oplus
F^+_\nu.\end{equation}

Another orthogonal decomposition of $L^2(\Omega^*_{\del X})$ is the Hodge
decomposition:
\begin{equation}\label{hodge}L^2(\Omega^*_{\del X})=\im d \oplus
\ker A\oplus \im d^* .\end{equation}

We introduce a notational convention: the decomposition \plref{eigenspace}
is compatible with the operators $d, d^*$ in the sense that
we have decompositions of domains:
\begin{equation}\begin{split}
     \cd(d)&=(F^-_\nu\cap\cd(d))\oplus E_\nu^-\oplus \ker A\oplus E_\nu^+\oplus
(F^+_\nu\cap\cd(d)),\\
     \cd(d^*)&=(F^-_\nu\cap\cd(d^*))\oplus E_\nu^-\oplus \ker A\oplus E_\nu^+\oplus
(F^+_\nu\cap\cd(d^*).
        \end{split}
\end{equation}
Note that $E_\nu^-\oplus \ker A\oplus E_\nu^+$ consists of smooth sections
hence $(E_\nu^-\oplus \ker A\oplus E_\nu^+)\cap\cd(d)\cap\cd(d^*)=E_\nu^-\oplus
\ker A\oplus E_\nu^+$. By slight abuse of notation we will write in the sequel
$d^{(*)}(F_\nu^{\pm})$ for
the image of $d^{(*)}$ on $F_\nu^{\pm}\cap\cd(d^{(*)})$.

The relations between the decompositions \eqref{eigenspace}
and \eqref{hodge} are summarized in the following useful lemma.

\begin{lemma} \label{lemma1}\hfill
\begin{enumerate}
\item $d(F_\nu^\pm)\subset F_\nu^\mp$ and $d^*(F_\nu^\pm)\subset
F_\nu^\mp$.
\item  $F_\nu^+=d(F_\nu^-)\oplus d^*(F_\nu^-)=(\ker d:F_\nu^+\to F_\nu^-)
\oplus (\ker d^*:F_\nu^+\to
F_\nu^-).$
\item $F_\nu^-=d(F_\nu^+)\oplus d^*(F_\nu^+)=(\ker d:F_\nu^-\to F_\nu^+)
\oplus (\ker d^*:F_\nu^-\to
F_\nu^+).$
\item  $d(E_{-\mu})\subset E_\mu$ and $d^*(E_{-\mu})\subset
E_\mu$, and for $\mu\ne 0$,
$E_\mu=d(E_{-\mu})\oplus d^*(E_{-\mu})$.
\item $\gamma(\ker d)= \ker d^*$ and $\gamma(\ker d^*)= \ker
d$.
\end{enumerate}
\end{lemma}

\begin{proof}  If $A\beta=\mu\beta$, then $A
d\beta=-dA\beta=-\mu d\beta$, and similarly $A
d^*\beta= -\mu d^*\beta$.  This proves the first assertion and the
first part of  the fourth assertion.

If $\beta\in F_\nu^+$, then $\beta$ is orthogonal to $\ker A$,  since
the decomposition \eqref{eigenspace} is orthogonal. Since the decomposition
\eqref{hodge} is also orthogonal, $\beta$ has the orthogonal
decomposition
$\beta=d\tau+d^*\sigma$.  Write $\tau=\tau_-+\tau_+\in F_\nu^-\oplus
F_\nu^+$, and similarly $\sigma=\sigma_-+\sigma_+$. Then
$$\beta =d\tau_- +d\tau_+ +d^*\sigma_-+d^*\sigma_+.$$
Since $\beta\in F_\nu^+$, the first assertion implies that
$d\tau_+=0=d^*\sigma_+$, so that $\beta=d\tau_-+d^*\sigma_-$. The
second assertion follows from this and the consequence
$d(F_\nu^-\oplus F_\nu^+)=\ker\  d:F_\nu^-\oplus F_\nu^+\to  F_\nu^-\oplus
F_\nu^+$ of the DeRham theorem.  The third assertion is proved similarly,
as is the second part of the fourth assertion.

The last assertion follows from the identity $d^*=-\hat{*}d\hat{*}$ and
the fact that $\gamma$ equals $\hat{*}$ up to a non--zero constant.
\end{proof}

Of particular concern will be the symplectic structure on $\ker A$.  The
isomorphism
$\gamma$ preserves
$\ker A$, satisfies $\gamma^2=-I$, and acts with signature zero, since
$(\del X,A)$ bounds
$(X,D)$. Therefore $\ker A$ is a finite--dimensional Hermitian symplectic
subspace  of $L^2(\Omega^*_{\del X})$.

Notice that the restrictions of $\langle\ , \ \rangle, \ \gamma$, and
$\omega$ to
$\ker A$ induce these structures on the cohomology  $H^*(\del X;\C^n_\alpha)$
via the Hodge and DeRham isomorphisms. The inner product $\langle\ , \
\rangle$ and complex structure $\gamma$ depend on the choice of
  Riemannian metric on $\del X$, but the symplectic structure $\omega$ does
not: if $\beta_1 \in \ker A$ is a $p$-form  and $\beta_2\in\ker A$ is a $2n-p$
form, then
\begin{equation}\label{cupproduct}\omega(\beta_1,\beta_2)=\langle \beta_1,
\gamma(\beta_2)\rangle=i^r
\int_{\del X}
\beta_1\wedge  \beta_2
\end{equation} where the constant $i^r$ depends only on $p$ and $n$
(and we have
suppressed the notation for the inner product in the flat $\C^n$ bundle
$E_{|\del X}$). Since wedge products and cup products correspond via the
DeRham isomorphism, $\omega$ coincides with the cup product up to a
power of $i$, and in particular is a homotopy invariant. To put this
differently, The cup product, together with the standard $U(n)$--invariant
Hermitian inner product on $\C^n$, induces a skew--hermitian form
$$\omega:H^*(\del X;\C^n)\times H^*(\del X;\C^n)\to \C, \ \
\omega(\beta_1,\beta_2)=i^r\big(\beta_1\cup \beta_2\big)\cap [\del X]$$
which is a homotopy invariant of the pair $(\del X,\alpha_{|\del X})$. Fixing a
Riemannian metric on $\del X$ induces a positive definite Hermitian inner product
and an isomorphism $\gamma$ on $\ker A$. The Hodge and DeRham theorems define an
isomorphism $\ker A\to H^*(\del X;\C^n_\alpha)$ which takes the form
$\langle x,\gamma(y)\rangle$ to the form $\omega(x,y)$.

The following lemma collects some useful information about symplectic subspaces
and symplectic reduction. For more details about symplectic reduction in this
setting the reader should consult Section \plref{Symplectic reduction}
and
\cite{Nic:MIS}.

\begin{lemma} \label{symred} \hfill
\begin{enumerate}
\item Let $S\subset L^2(\Omega^*_{\del X})$ be a closed subspace satisfying
$\gamma(S)\perp S$. Then $S\oplus \gamma(S)$ is a Hermitian symplectic
subspace of $ L^2(\Omega^*_{\del X})$, and $S$ is a Lagrangian subspace of
$S\oplus
\gamma(S)$.

\item If $\nu\ge 0$, then $F^-_\nu\oplus F^+_\nu$, $E^-_\nu\oplus E^+_\nu$,
$E^-_\nu\oplus \ker A\oplus E^+_\nu$, and $d(E^{\pm}_\nu)\oplus d^*(E^\mp_\nu)$
are Hermitian symplectic subspaces of $ L^2(\Omega^*_{\del X})$.

\item Given a Lagrangian subspace $L\subset   L^2(\Omega^*_{\del X})$ so that
$(L, F_0^-)$ form a Fredholm pair of subspaces, then
\begin{equation} R_\nu(L):={{L\cap(F^-_\nu\oplus E^-_\nu\oplus \ker A\oplus
E^+_\nu)}\over{ L\cap F^-_\nu}}\subset E^-_\nu\oplus \ker A\oplus E^+_\nu
\end{equation}
is a Lagrangian  subspace, called the {\em symplectic reduction  of $L$ with
respect to }$F_\nu^-$.
\end{enumerate}
\end{lemma}
\begin{proof}

1.  Notice that $\gamma$ preserves $S\oplus \gamma(S)$. Let $K_{\pm i}$ denote
the $\pm i$ eigenspaces of $\gamma$ acting on $S\oplus \gamma(S)$. It is easy
to check that since $\gamma(S)$ is
orthogonal to
$S$,   the projections   $ S\oplus
\gamma(S)\to  K_{\pm i}$ restrict  to   isomorphisms on $S$.   Thus the
$\pm i$ eigenspaces of
$\gamma$ on $S\oplus \gamma(S)$ have the same dimension (or are both infinite).
This shows that $S\oplus \gamma(S)$ is a symplectic subspace of $
L^2(\Omega^*_{\del X})$.  Clearly $S$ is a Lagrangian subspace of $S\oplus
\gamma(S)$.

2.  For  $F^-_\nu\oplus F^+_\nu$, take $S=F^-_\nu$ and apply the
first assertion.
For
$E^-_\nu\oplus E^+_\nu$, take $S=E^-_\nu$.  For $d(E^{\pm}_\nu)\oplus
d^*(E^\mp_\nu)$, take $S=d(E^{\pm}_\nu)$; then $\gamma(S)=\hat{*}S
=\hat{*}d(E^{\pm}_\nu)=\hat{*}d(\hat{*}E^{\mp}_\nu)=d^*(E^{\mp}_\nu) $.
That $\ker A$ is
symplectic was discussed above; hence the direct sum $E^-_\nu\oplus
\ker A\oplus
E^+_\nu$ is symplectic.

3. We apply Proposition \plref{symplreduction} with $U=F_\nu^-\oplus
   E_\nu^-\oplus
\ker A\oplus E_\nu^+$. We have $\Ann(U)=F_\nu^-$ and $U\cap\gamma U=
   E_\nu^-\oplus\ker A\oplus E_\nu^+$. Since $(L,F_0^-)$ form a Fredholm pair
and $F_0^-/F_\nu^-$ is finite--dimensional, also $(L,\Ann(U))=(L,F_\nu^-)$
is Fredholm. Consequently $L+\Ann(U)$ is closed and we reach the
desired conclusion using Proposition \plref{symplreduction}.
\comment{
3.  Write $F^\pm$ for $F^\pm_\nu$ and $K$ for $E^-_\nu\oplus \ker A\oplus
E^+_\nu$. Hence we must show that $R_\nu(L)\subset K$ is a Lagrangian
subspace.

Since
$(L, F_0^-)$ form a Fredholm pair, so is
$(L, F^-\oplus K)$. Hence $L\cap (F^-\oplus K)$ is finite--dimensional and
$\gamma(L)+F^+$ has finite codimension.  A straightforward argument shows that
$$\gamma(L)+ F^+=(L\cap (F^-\oplus K))^\perp.$$

Suppose that $k\in K$ such that $\gamma(k)$ is perpendicular to
$R_\nu(L)$. Then
for each $\ell\in L\cap(F^+\oplus K)$,
$\langle\ell,\gamma(k)\rangle=\langle\proj_K(\ell),\gamma(k)\rangle=0$.  Thus
$\gamma(k)\in (L\cap (F^-\oplus K))^\perp=\gamma(L)+ F^+$.  Write
$\gamma(k)=\gamma(m)+f^+$. Then $k=-m+\gamma(f^+)$ so that
$-m= -\gamma(f^+)+k $ lies in $L\cap(F^-\oplus K)$.  Thus $k\in
R_\nu(L)$.  This
shows that $\gamma(R_\nu(L)^\perp)\subset R_\nu(L)$ and for dimensional reasons
this implies that $\gamma (R_\nu(L))=R_\nu(L)^\perp$ (the orthogonal
complements taken in $K$), so that $R_\nu(L)$ is a Lagrangian.

The assertions about the continuity of symplectic reduction are
well--known, and
we omit the examples.}
\end{proof}

In preparation for what follows we define the following enlargements of
$X$.  Given $r\ge 0$ define
$$X_r=([-r,0]\times \del X)\cup X$$
and
$$X_\infty= ((-\infty,0]\times \del X)\cup X.$$
Thus $X_r$ has a collar of length $r$ attached to $X$ and $X_\infty$ is
obtained from $X$ by attaching an infinitely long collar.  Equation
\eqref{cylinder} can be used on the collar to define a natural extension
of
$D$ to
$X_r$ and $X_\infty$.

  The key to
identifying the adiabatic limit of the Calder{\'o}n projector is the
following result.

\begin{proposition}\label{kerneld} Suppose that the boundary of $X$ is
non--empty, and  suppose that $\beta\in \Omega^{\text{\rm \tiny even}}_X$
satisfies
   $D \beta =0$ and $r(\beta)\in F_0^-\oplus\ker A =\rmspan\{\psi_\mu\ |
\
\mu\le 0\}$. Then
$d\beta=0$,
$d(*\beta)=0$, and $d (r(\beta)) =0$.
\end{proposition}
\begin{proof}  Naturality of the exterior derivative implies that
$d(i^*(z))=i^*(dz)$ for any $z\in \Omega^k_X$. It suffices, therefore,
to show
    that $d\beta=0$ and $d(*\beta)=0$, since
    $r(\beta)=
i^*(\beta) + i^*(*\beta)$ and hence
$$d(r(\beta))=d(i^*(\beta) + i^*(*\beta))=i^*(d\beta+d*\beta)=0.$$

Following \cite{AtiPatSin:SARI}, since $D\beta=0$ and $r(\beta)\in
F_0^-\oplus \ker A$,
$\beta$ has a Fourier expansion on the collar $[0,\ep)\times \del X$ of the
form
\begin{equation}\label{Fourier}\beta_{|[0,\ep)\times \del X}=  \sum_{\mu<0}
c_\mu e^{-x\mu}\psi_\mu + k,
\end{equation}
where $k\in \ker A$, $x\in [0,\ep)$, and $ \psi_\mu\in E_\mu$.
Equation \plref{Fourier} can be used to extend $\beta$ to a bounded form
on
$X_\infty$  so that the extension still
satisfies $D\beta=0$.

Notice that $dk=0$   since $k\in \ker A$ and $k$ is
independent of the collar parameter. Thus   $d\beta$
decays exponentially on the infinite collar $ (-\infty,0]\times \del X$. Write
   $\beta=\sum \beta_{2p}$. Then
\begin{eqnarray*}
\langle d*\beta_{2p},*d\beta_{2(p-1)}\rangle_{L^2(\Omega^*_{X_r})} &=&
\pm\int_{X_r} d*\beta_{2p}\wedge
\beta_{2(p-1)} \\
&=&\pm  \int_{X_r} d(*\beta_{2p}\wedge d\beta_{2(p-1)}) \\
&=&\pm \int_{\del X\times \{-r\}} i^*(*\beta_{2p})\wedge
i^*(d\beta_{2(p-1)}).
\end{eqnarray*}
The last step follows from Stokes's theorem. As $r$ increases to
infinity, the last integral converges to zero since $ *\beta_{2p} $
is bounded on $X_\infty$ and $d\beta_{2(p-1)}$ exponentially decays.
It
follows that   $d*\beta_{2p}$ and $*d\beta_{2(p-1)}$ are orthogonal in
$L^2(\Omega^{2(n-p-1)}_{X_\infty})$.
Now
$$0=D\beta=i^{n+1}\sum_p(-1)^p(d*\beta_{2p} + *d\beta_{2(p-1)}),$$
with this sum expressed as a sum of homogeneous components. Thus
$d*\beta_{2p}$ and $*d\beta_{2(p-1)}$ both vanish for each $p$, and
therefore     $d*\beta$ and $*d\beta$  both vanish.
\end{proof}

As an  application of Proposition  \plref{kerneld} we can identify the {\it
limiting values of extended
$L^2$ solutions of $D\beta=0$} in the sense of \cite{AtiPatSin:SARI}.
Recall that this is the subspace of $\ker A$ defined  by
\begin{equation}\label{vee1}
V_\alpha =\Bigsetdef{k}{\begin{array}{l} \text{there exists a
}\beta\in\Omega^{\text{\tiny even}}_X
\text{ with } D\beta=0\\
\text{and } r(\beta)=f_-+k\in F^-_0\oplus \ker A\end{array}}.\end{equation}
The terminology is justified by the Fourier expansion \eqref{Fourier}.
In light of the unique continuation property for $D$ (which says that
for each $\ell \in L_X$ there exists a unique $\beta$ with
$D\beta=0$ and $r(\beta)=\ell$), it is easy to see that
$V_\al$ has the alternative description as a symplectic reduction:
\begin{equation}\label{vee}V_\alpha=R_0(L_X)={{L_X\cap (F_0^-\oplus
\ker A)}\over {L_X\cap  F_0^- }}
\subset \ker A.\end{equation}

Equation
  \eqref{vee} says that
$V_\alpha$ is   the symplectic reduction of the Cauchy data space
$L_X$ with respect to   subspace $F_0^-$. Using Lemma \plref{symred} it
follows that    $V_\alpha$ is   a Lagrangian subspace of $\ker A$.

\vskip.3in

The kernel of $A$ is identified via the Hodge and DeRham theorems with
the cohomology $H^*(\del X,\C^n_\alpha)$. The next result identifies
$V_\alpha$.

\begin{corollary}\label{veehomology}  The space $V_\al$ of limiting values of
extended
$L^2$ solutions to $D\beta=0$ on $X_\infty$ is identified via the Hodge and
DeRham theorems with the image of the cohomology of $X$ in the cohomology of
$\del X$ (with local coefficients in the corresponding flat $\C^n$ bundle):
$$V_\alpha=\im i^*:H^*(X;\C^n_\alpha)\to
H^*(\del X;\C^n_\alpha).$$
\end{corollary}

\begin{proof} Proposition \plref{kerneld} shows  that if
$\beta\in\Omega^{\text{\tiny even}}_X
$ satisfies $D\beta=0$ and $r(\beta) \in F_0^-\oplus \ker A $, then $
\beta $ and $ *\beta $ are closed forms.  Thus they represent classes
in $H^*(X;\C^n_\alpha)$. Since
$r(\beta)=i^*(\beta)+i^*(*\beta)$, it follows that
   $r(\beta)$ is a closed form on $\del X$ representing a class in $ \im
i^*:H^*(X;\C^n_\alpha)\to H^*(\del X;\C^n_\alpha).$  The identification of cohomology
with harmonic forms   takes
$[r(\beta)]=[f_-+k]$ to $k$ and so
$$V_\alpha\subset \im i^*:H^*(X;\C^n_\alpha)\to H^*(\del X;\C^n_\alpha).$$
The space $V_\alpha$ is a Lagrangian subspace, as is
$ \im i^*:H^*(X;\C^n_\alpha)\to H^*(\del X;\C^n_\alpha) $ by a standard
argument using Poincar\'e duality. Since any two Lagrangian subspaces of a
finite--dimensional symplectic vector space have the same dimension,
$V_\al=\im i^*$.
\end{proof}

It follows from  Lemma
\plref{lemma1} that  $E^\pm_\nu =d(E^\mp_\nu)\oplus d^*(E^\mp_\nu)$, and so
the decomposition \plref{eigenspace} can be refined to
\begin{equation}\label{bigsum}
L^2(\Omega^*_{\del X})=
F^-_\nu\oplus d(E^+_\nu)\oplus d^*(E^+_\nu)\oplus \ker A\oplus
d(E^-_\nu)\oplus d^*(E^-_\nu)\oplus F^+_\nu.
\end{equation}
The terms in this decomposition are arranged according to increasing
eigenvalues.   We will find it  convenient   to rewrite this in a
different order, as a symplectic direct sum of     symplectic subspaces:
\begin{equation}\label{bigsum2}
L^2(\Omega^*_{\del X})=
(F^-_\nu\oplus F^+_\nu)\oplus(d(E^+_\nu)\oplus d^*(E^-_\nu))\oplus
(d^*(E^+_\nu)\oplus d(E^-_\nu))\oplus\ker A.
\end{equation}

We will refer to the decomposition \eqref{bigsum2} frequently. Notice that
$F^-_\nu\oplus F^+_\nu$ is infinite--dimensional and the other three
symplectic summands in this decomposition have finite dimension.
\vskip.3in
There exists a $\nu\ge 0$ so that
the Cauchy data space $L_X$ of $D$ is transverse to $F^-_\nu$. This is
because $L_X\cap F_0^-$ is finite--dimensional, and as $\nu$
increases,
   $L_X\cap F^-_\nu$ decreases to zero. Nicolaescu calls the smallest
such $\nu$ the {\it non--resonance level} for $D$.

We can now state and prove a theorem identifying the limit of the
Calder{\'o}n projectors of $D$ acting on $X_r$ as $r$ goes to infinity.
Denote by $L^r_{X}$   the Cauchy data space  (i.e. the image of
the Calder\'on projector)  of $D$ acting on
$X_r=  ([-r,0]\times \del X)\cup X$.

\begin{theorem}\label{thmonadlim} Let $X$ be an odd--dimensional manifold with
boundary  and $D$ the odd signature operator coupled to a flat connection $B$
acting on
$X$ as above. Let
$\nu\ge 0$ be any number greater than or equal  to the non--resonance level for
$D$.

Then there exists a subspace
$$W_\al\subset d(E_\nu^+)\subset F^-_0$$ isomorphic to the image of
$$H^{\text{\rm \tiny even}}(X,\del X;\C^n_\alpha)\to H^{\text{\rm \tiny
even}}(X;\C^n_\alpha)$$ so that if $W_\al^\perp$ denotes the orthogonal
complement of
$W_\al$ in
$d(E^+_\nu)$, then with respect to the decomposition \eqref{bigsum2} of
$L^2(\Omega^*_{\del X})$ into symplectic subspaces, the adiabatic limit of the
Cauchy data spaces decomposes as a direct sum of Lagrangian subspaces:
\begin{equation}\label{adiablim}\lim_{r\to
\infty}L^r_{X} = F^+_\nu\oplus (W_\al\oplus \gamma(W_\al^\perp))\oplus
  d(E^-_\nu) \oplus V_{\al}.\end{equation} where
$V_\alpha\subset\ker A=H^*(\del X;\C^n_\alpha)$ denotes the image of
$H^*(X;\C^n_\alpha)\to H^*(\del X;\C^n_\alpha)$.

\end{theorem}
\begin{proof}
Lemma \plref{symred} shows that the finite--dimensional  vector space
$E_\nu^-\oplus
\ker A\oplus E^+_\nu$  is a   symplectic
subspace of
$L^2(\Omega^*_{\del X})$.

 Let $R_\nu(L_X)\subset
E_\nu^-\oplus
\ker A\oplus E^+_\nu$ be the symplectic reduction of $L_X$ with respect to
the isotropic subspace $F^-_\nu$ as in Lemma \plref{symred}:
\begin{eqnarray*}R_\nu(L_X)&=&{{L_X\cap (F^-_\nu\oplus E_\nu^-\oplus
\ker A\oplus E^+_\nu)}
\over {L_X\cap F^-_\nu}}\\
&=&
\text{proj}_{E_\nu^-\oplus \ker A\oplus
E^+_\nu}\left( L_X\cap (F^-_\nu\oplus E_\nu^-\oplus
\ker A\oplus E^+_\nu)\right).\end{eqnarray*}
Then $R_\nu(L_X)$ is a Lagrangian subspace of $E_\nu^-\oplus \ker A\oplus
E^+_\nu$.

Nicolaescu's theorem \cite[Theorem 4.9]{Nic:MIS} says
\begin{equation}\label{nicoadiablim}\lim_{r\to\infty}L^r_{X}=
\left(\lim_{r\to
\infty} e^{rA}R_\nu(L_X)\right) \oplus F^+_\nu\end{equation}
(The sign in the exponent $e^{rA}$ differs from \cite{Nic:MIS} because in
that paper the collar of $X_r$ is parameterized as $\del X\times [0,r]$.) Thus
we need only to identify the limit of
$e^{rA}R_\nu(L_X)$. To help with the rest of the argument the reader
should observe that the dynamics of $e^{rA}$ favor  the vectors with a
non--zero component in
    eigenspaces corresponding to positive eigenvalues.

   Let $\mu_1<\mu_2<\cdots<\mu_q$ denote the complete list of eigenvalues
of $A$ in the range $[-\nu,\nu]$. Thus
$$E_\nu^-\oplus \ker A\oplus
E^+_\nu=E_{\mu_1}\oplus E_{\mu_2}\oplus\cdots\oplus E_{\mu_q}.$$
Given $\ell\in R_\nu(L_X)$, we use this decomposition to write
$$\ell=(\ell_1,\ell_2,\cdots, \ell_q).$$

Let $\mu(\ell)$ denote the  largest $\mu_i$ so that $\ell_i$ is non--zero
(and hence $\ell_{\mu(\ell)+1}=\cdots=\ell_{q}=0$). Then
$$\lim_{r\to\infty}e^{rA} \left(\frac{1}{e^{r\mu(\ell)}} \ell\right)
=(0,0,\cdots, 0,\ell_{\mu(\ell)},0,\cdots,0).$$
This shows that
$$\lim_{r\to\infty}e^{rA}R_\nu(L_X)=
L_{\mu_1}\oplus L_{\mu_2}\oplus\cdots\oplus L_{\mu_q} \subset
E_{\mu_1}\oplus E_{\mu_2}\oplus\cdots\oplus E_{\mu_q},$$
where
\begin{eqnarray}\label{defofLi}
L_{\mu_i}&=&
\text{proj}_{E_{\mu_i}}\left(R_\nu(L_X)\cap (E_{\mu_1}\oplus
\cdots\oplus E_{\mu_i})\right) \nonumber\\
&=&\text{proj}_{E_{\mu_i}}\left(L_X\cap (F_\nu^-\oplus E_{\mu_1}\oplus
\cdots\oplus E_{\mu_i})\right).
\end{eqnarray}
Write
$$L^-=\moplus_{\mu_j<0}L_{\mu_j}\subset E^-_\nu,$$
$$L^0=\moplus_{\mu_j=0}L_{\mu_j}\subset \ker A,$$
and
$$L^+=\moplus_{\mu_j>0}L_{\mu_j}\subset E^-_\nu,$$
so that
$$\lim_{r\to\infty} e^{rA}R_\nu(L_X)=L^-\oplus L^0\oplus L^+\subset E_\nu^-\oplus
\ker A\oplus E^+_\nu.$$
Set \begin{equation}W_\al:=L^-.\end{equation}

\begin{lemma}\label{technical}\hfill
\begin{enumerate}
\item $L_0=V_\alpha$.
\item The spaces
\begin{enumerate}
\item $W_\alpha=L^-$,
\item the image of $H^{\text{\rm \tiny even}}(X,\partial X;\C^n_\alpha)\to
H^{\text{\rm \tiny even}}(X;\C^n_\alpha)$,
\item $L_X \cap F_0^-$, and
\item the $L^2$ solutions of $D x=0$ on $X_\infty$
\end{enumerate}
are all isomorphic.
\item $L^-\subset d(E^+_\nu)$.
\end{enumerate}
\end{lemma}

Assuming these three facts, the rest of the proof of Theorem
\plref{thmonadlim} is completed as follows.

Note that $W_\al\subset d(E^+_\nu)\subset d(E^+_\nu)\oplus
d^*(E^+_\nu)=E^-_\nu$. We define  $W_\al^\perp$ to be the orthogonal
complement of $W_\al$ in $d(E^+_\nu)$. Since
$$W_\al\oplus L^+=L^-\oplus L^+\subset E_\nu^- \oplus
E^+_\nu=\left(d(E^+_\nu)\oplus
d^*(E^+_\nu)\right)\oplus\left(d(E^-_\nu)\oplus d^*(E^-_\nu)\right)$$
is a Lagrangian subspace (obtained by modding  out $L^0$ and $\ker A$),
it follows from Lemma
\plref{lemma1} that
$$L^+=d(E_\nu^-)\oplus \gamma(W_\al^\perp) \subset  d(E^-_\nu)\oplus
d^*(E^-_\nu),$$
completing the proof of Theorem \plref{thmonadlim}.\end{proof}
\vskip.3in
{\sc Proof of Lemma } \plref{technical}.
The first assertion follows immediately by comparing Equations
\plref{vee} and \plref{defofLi}.

For the second    assertion, Equation \eqref{defofLi} shows
that if $m\in L^-$, there exists a $\mu_i<0$ and an
$$\ell=(f, \ell_{\mu_1},\ell_{\mu_{2}}
    \cdots,\ell_{\mu_{i}})\in
(F^-_\nu\oplus E_{\mu_1}\oplus
E_{\mu_{2}}\oplus\cdots\oplus E_{\mu_i})\cap L_X$$
with $m=\ell_{\mu_i}$.  This sets up an identification of $L^-$ with
$F_0^-\cap L_X$. The unique continuation property identifies this (via
the restriction map $r$)  with the kernel of $D$ with $P_{\ge 0}$ boundary
conditions, which, by Proposition \plref{kerneld} and Equation \plref{Fourier}
(with $k=0$), is the same as the space of $L^2$ harmonic forms in
$\Omega^{\text{\tiny even}}_{X_\infty}$. The space of $L^2$ harmonic $p$-forms
is shown to be isomorphic to the image of
$H^{p}(X,\partial X;\C^n_\alpha)\to H^{p}(X;\C^n_\alpha)$ in
\cite[Proposition 4.9]{AtiPatSin:SARI}.

The third assertion also follows, since if
$\ell=r(\beta)$, then   Proposition
\plref{kerneld} says $d\ell=d(r(\beta))=0$. But
$$0=d\ell=df +d\ell_{\mu_1}+d\ell_{\mu_{2}}+\cdots +d\ell_{\mu_{i}} $$
and since $d(E_\mu)\subset E_{-\mu}$,
$$0=d\ell_{\mu_i}=dm.$$
Hence (since $\mu_i< 0$)
$$m\in \ker (d:E_{\mu_i}\to E_{-\mu_i})=
d(E_{-\mu_i})\subset d(E_\nu^+),$$
completing the proof of Lemma \plref{technical}
\qed

\begin{remark}
Notice that $W_\al^\perp$ denotes the orthogonal complement to $W_\al$ in
the finite--dimensional space $d(E^+_\nu)$, not in $L^2(\Omega^*_{\del X})$.
\end{remark}

We adopt the following notation in the rest of this section to deal with
boundary conditions. Given a manifold with boundary $X$, the odd signature
operator $D=D_B$ coupled to  a flat connection $B$ on $X$ as above,   and a
Lagrangian subspace
$L\in
\mathcal{L}_{\text{\rm Fred}}$, then let
$\eta(D, X; L)$  denote the
$\eta$--invariant of the Dirac operator
$D$ with boundary conditions given by the orthogonal projection to $L$. Thus,
$$\eta(D, X; L):= \eta(D_{\proj_{L}}, X)$$
in the previous notation.  The same notation applies to
the reduced $\eta$--invariant $\etab$.

In a similar manner, given appropriate Lagrangian subspaces $L,M,N$
of a Hermitian symplectic Hilbert space we will use
$\tau_\mu(L,M,N)$ to denote the triple index of the corresponding projections
$\tau_\mu(\proj_L,\proj_M,\proj_N)$ (cf. Section \plref{secsixtwo}).

\vskip.3in
  Suppose that $B$ and $B'$ are flat connections on $X$ in temporal gauge
near $\partial X$ such that the holonomy representations
$\alpha,\alpha':\pi_1X\to U(n)$ of
$B, B'$ are conjugate. Then there exists a gauge transformation $g$ so that
on a collar $[0,\epsilon)\times \partial X$, $g=\pi^*(h)$ for a gauge
transformation $h$ on $\partial X$ satisfying $B'= g\cdot B$. Hence the
restrictions $b, b'$ of $B,B'$ to the boundary satisfy $b'=h\cdot b$. We have
$$D_{B'} =D_{gB}=gD_Bg^{-1},$$ and $$A_{b'}=A_{hb}=hA_bh^{-1}. $$

In particular, $h$ takes the positive (resp. negative) eigenspan of $A_b$ to
the positive (resp. negative) eigenspan of $A_{b'}$, and gives an isomorphism
$\ker A_b\to\ker A_{b'}$ which coincides via the Hodge and DeRham theorems
with the isomorphism $H^*(\partial X; \C^n_\alpha)\to H^*(\partial X;
\C^n_{\alpha'})$  induced by conjugating the holonomies $\alpha,\alpha'$.
Thus if $K\subset \ker A_b$ is a Lagrangian subspace, the
$\lambda$--eigenspace  of $D_B$ with $F^+_0(b)\oplus K$ boundary conditions
is sent by $g$ to the $\lambda$-eigenspace  of $D_{B'}$ with $F^+_0(b')\oplus
h(K)$ boundary conditions.

Since any representation $\alpha:\pi_1X\to U(n)$
is the holonomy representation of a flat connection $B$, we conclude that
given a representation $\alpha$ and a Lagrangian subspace $K\subset
H^*(\partial X;\C^n_\alpha)$ (recall that the symplectic structure $\omega$
on $H^*(\partial X;\C^n_\alpha)$ is defined by the cup product),  the quantity
$\eta(D,X; F^+_0\oplus K)$ is unambiguously defined, i.e. it is
independent of the choice of flat connection $B$ in temporal gauge with
holonomy conjugate to $\alpha$,  and the Lagrangian  in $\ker A_b$
corresponding to $K$ via the Hodge and DeRham theorems is well defined.  Of
course $\eta(D,X; F^+_0\oplus K)$ may  depend  on the choice of Riemannian
metric on $X$.

\vskip.3in

We can now turn to the splitting problem for the $\eta$--invariant of the odd
signature operator. As in earlier sections suppose that $M=M^+\cup M^-$ is a
closed manifold decomposed into 2 submanifolds along a separating hypersurface
$N$. Assume that $N$ has a
    neighborhood isometric to $N\times[-1,1]$. Suppose that $B$ is a
flat connection on $M$ in temporal gauge on $N\times[-1,1]$.

  As we have seen, because the outward normal for $M^+$ is the inward
normal for $M^-$, the operators $\gamma$ and $A$ for $M^-$ are related
to those for $M^+$ by a change in signs.  This has the following
consequences. First, whereas  the conclusion of Theorem \plref{thmonadlim}
identifies  the limit of the Cauchy data spaces $L^r_{M^+}$ of $D$ acting on
$M^+_r$,  $
\lim\limits_{r\to
\infty}L^r_{ M^+ }$ as
\begin{equation}\label{plusside}F^+_\nu\oplus (W_{+,\al}\oplus
\gamma(W_{+,\al}^\perp))\oplus
  d(E^-_\nu) \oplus V_{+,\al}\end{equation}
(in the decomposition \eqref{bigsum2})
for $W_{+,\al}\subset d(E_\nu^+)\subset F^-_0$ a space isomorphic to
the image
$$\im \bigl(H^{\text{\tiny even}}(M^+,\partial M^+;\C^n_\alpha)\to
H^{\text{\tiny even}}(M^+;\C^n_\alpha)\bigr)$$
and $V_{+,\al}\subset \ker A$ a space isomorphic to
$$\im \bigl(H^{\text{\tiny even}}(M^+;\C^n_\alpha)\to H^{\text{\tiny
even}}(N;\C^n_\alpha)\bigr).$$
For $M^-$ the conclusion is that
$ \lim\limits_{r\to \infty}L^r_{ M^- }$ is
\begin{equation}\label{minusside}F^-_\nu\oplus  d(E^+_\nu) \oplus
(\gamma(W_{-,\al}^\perp)\oplus W_{-,\al})
\oplus V_{-,\alpha},\end{equation}
where
$W_{-,\al}\subset d(E^-_\nu)\subset F_0^+$ is
isomorphic to
the image
$$\im \bigl(H^{\text{\tiny even}}(M^-,\partial M^-;\C^n_\alpha)\to
H^{\text{\tiny even}}(M^-;\C^n_\alpha)\bigr)$$
and $V_{-,\al}\subset \ker A$ is a space isomorphic to
$$\im \bigl(H^{\text{\tiny even}}(M^-;\C^n_\alpha)\to H^{\text{\tiny
even}}(N;\C^n_\alpha)\bigr).$$ (We assume that $\nu$ has been chosen greater than or
equal to the non-resonance level for $D$ acting on both $M^+$ and $M^-$.)

Theorem \plref{symsplit2} calculates the $\eta$--invariant of $D$ acting on $M$
in terms of the
$\eta$--invariants of $D$ on $M^+$ and $M^-$.   Take $P$ to be the
Atiyah--Patodi--Singer boundary projection $P=P^+(V)$ for some Lagrangian
subspace $V\subset \ker A$.  Then Theorem \plref{symsplit2} says
\begin{equation}\label{eq8.1}\begin{split}
\etab(D,M)&=\etab(D,M^+; V\oplus F^+_0)
+\etab(D,M^-;F^-_0\oplus \gamma(V) )\\
&\quad-\tau_\mu(I-P_{M^-},P^+(V),P_{M^+}).
                 \end{split}
\end{equation}
(Recall that $P_{M^\pm}$ denotes the Calder\'on \ projectors onto the Cauchy
data spaces $L_{M^\pm}$.)

\begin{theorem}\label{triplebaby} Let $D$ denote the odd signature operator
coupled to a flat connection.  For any Lagrangian subspace
$V\subset  \ker A $,
$$\etab(D,M)=\etab(D,M^+; V\oplus F^+_0)
+\etab(D,M^-;F^-_0\oplus \gamma(V)
)-\tau_\mu(\gamma(V_{-,\alpha}),V,V_{+,\alpha}),$$
where
$\tau_\mu(\gamma(V_{-,\alpha}),V,V_{+,\alpha})$ refers to the triple index in
the finite--dimensional space $\ker A \cong H^*(N;\C^n_\alpha) $.

  The triple index  $\tau_\mu(\gamma(V_{-,\alpha}),V,V_{+,\alpha})$  vanishes if
$V=V_{+,\alpha}$ or $V=\gamma(V_{-,\alpha})$ and so
\[ \begin{split}\etab(D,M)&=
\etab(D,M^+; V_{+,\alpha}\oplus F^+_0)
+\etab(D,M^-;F^-_0\oplus \gamma(V_{+,\alpha}) )\\
&=\etab(D,M^+; \gamma(V_{-,\alpha})\oplus F^+_0)
+\etab(D,M^-;F^-_0\oplus V_{-,\alpha} ).\end{split}\]

In particular, if $H^*(N,\C^n_\alpha)=0$, then
$$\etab(D,M)=\etab(D,M^+;F^+_0)
+\etab(D,M^-;F^-_0).$$
\end{theorem}
\vskip.4in

The
main advantage that Theorem \plref{triplebaby} has over Theorem \plref{symsplit2} is
that the \Calderon \ projectors have been replaced by the Atiyah--Patodi--Singer
projections.

We postpone the proof of Theorem \plref{triplebaby} until after two lemmas are
in place.   The basic idea is to apply Lemma \plref{triplehomotopy} to the paths
obtained by stretching the Cauchy data spaces to their adiabatic limits.

\begin{lemma}\label{algtoplemma} $W_{+,\al}\oplus W_{-,\al} \oplus
(V_{+,\al}\cap V_{-,\alpha})$ is isomorphic to
   $H^{\text{\rm \tiny even}}(M;\C^n_\al)$.
\end{lemma}
\begin{proof}
For each integer $k$  let (twisted coefficients in $\C^n_\al$ are to be
understood for all cohomology groups)
$$W_\pm^k=\im \bigl(H^k(M^\pm,\del M^\pm)\to
H^k(M^\pm)\bigr)=\ker \bigl(i^*_\pm:H^k(M^\pm)\to H^k(N)\bigr)$$  and let
$$V_\pm^k= \im \bigl(i^*_\pm:H^k(M^\pm)\to
H^k(N)\bigr) .$$

Consider the map $$\Psi^k:H^{k}(M^+)\oplus H^{k}(M^-)\to H^{k}(N),
\ \ \Psi^k(m_+,m_-)=  i_+^*(m_+)-i_-^*(m_-)$$ in
   the Mayer--Vietoris sequence for $M=M^+\cup_N M^-$. Then there is a
short exact sequence
   \begin{equation}\label{MV1}
0\to W^{k}_+\oplus W^{k}_-\to \ker \Psi^{k}\mapright{\beta}
V^{k}_+\cap V^{k}_-\to 0,
\end{equation}
   where $\beta(m_+, m_-)=i_1^*(m_+)=i_2^*(m_-)$.
Moreover, the Mayer--Vietoris sequence gives a short exact sequence
   \begin{equation}\label{MV2}
0\to\text{coker }\Psi^{k-1}\to
H^{k}(M)\to \ker
\Psi^{k}
\to 0.
\end{equation}
Thus
\begin{equation}\label{no1}
\dim H^k(M) =\dim \text{coker }\Psi^{k-1}
+\dim (V_+^k\cap V_-^k)+\dim W^k_+ +\dim W^k_-.\end{equation}
Also
\begin{equation}\label{no2}\begin{split}
\dim &\text{ coker }\Psi^{k-1}=\dim
H^{k-1}(N)/(V_+^{k-1}+V_-^{k-1})\\ &=\dim
H^{k-1}(N)-\dim V_+^{k-1}-\dim V_-^{k-1}+ \dim (V_+^{k-1}\cap V_{-}^{k-1}).
               \end{split}
\end{equation}

Combining \eqref{no1} and \eqref{no2} and summing up over $k$ even
one obtains
\begin{equation}\label{no3}\begin{split}
\sum \dim H^{2k}(M)&=\sum\dim(V_+^k\cap V_-^k) + \sum \dim W^{2k}_+
+\sum \dim W^{2k}_-\\ &\quad+ \sum \dim H^{2k-1}(N)- \sum\dim
V_+^{2k-1}-\sum\dim V_-^{2k-1}.\end{split}\end{equation}

The symplectic space $H^*(N)$ decomposes as a symplectic sum $H^{\text{\tiny
even}}(N) \oplus  H^{\text{\tiny odd}}(N)$ (one way to see this is to
notice that
$\hat{*}$
  and hence $\gamma$ preserves the parity of a harmonic form  since $N$
is
$2n$--dimensional). The Lagrangian subspace $V^*_+=\sum V^k_+$ decomposes
accordingly into a sum of Lagrangian subspaces $\oplus V^{2k}_+\oplus
V^{2k-1}_+$. Hence
$\dim (\oplus V^{2k-1}_+)= \tfrac{1}{2}\dim H^{\text{\tiny
odd}}(N) $.  Similarly $\dim (\oplus V^{2k-1}_-)= \tfrac{1}{2}\dim
H^{\text{\tiny
odd}}(N)
$. Thus the last three terms in \eqref{no3} cancel. Since $W_{\pm, \al}=\oplus
W^{2k}_\pm$ and $V_{\pm,\al}=\oplus V_\pm^k$,
$$\dim H^{\text{\tiny even}}(M)=\dim (V_{+,\al}\cap V_{-,\al}) + \dim W_{+,\al}
+\dim W_{-,\al},$$  completing the proof of Lemma \plref{algtoplemma}.\end{proof}

\begin{lemma} \label{Lemma-8.9} Let $V\subset \ker A$ be a Lagrangian subspace.
Let
$L_{M^\pm}^r$  denote the Cauchy data space for $D$ acting on
$ M^\pm_r = ([-r,0]\times N)\cup M^\pm $ when $r<\infty$
and let $L_{M^\pm}^\infty$ be the adiabatic limit
$\lim\limits_{r\to \infty}L_{M^\pm}^r$ which was identified in Theorem
\plref{thmonadlim}.
\begin{enumerate}
\item The dimension of the intersection $ L_{M^-}^r \cap L_{M^+}^r$ is
independent of $r\in[0,\infty]$.
\item The dimension of the intersection $L_{M^-}^r\cap (F_0^+\oplus V)$ is
independent of $r\in [0,\infty]$.
\item The dimension of the intersection $  (F_0^-\oplus
\gamma(V))\cap L_{M^+}^r$
is independent of $r\in [0,\infty]$.
\item $\tau_\mu(\gamma(L_{M^-}^\infty), F_0^+\oplus V, L_{M^+}^\infty)=
\tau_\mu(\gamma(V_{-,\alpha}), V, V_{+,\alpha})$.

\end{enumerate}
\end{lemma}

  \begin{proof}

1.  For $r<\infty$, the intersection $ L_{M^-}^r \cap L_{M^+}^r$ is
isomorphic to
the kernel of $D_B$ acting on the closed manifold $M_r= M^+_r\cup  M^-_r$
obtained by stretching the collar of $N$.  But this kernel is a homotopy
invariant, isomorphic to $H^{\rm{\scriptstyle even}}(M;\C^n_\alpha)$, and in
particular its dimension is independent of $r$.

To compute $L_{M^-}^\infty\cap L_{M^+}^\infty$, we use Theorem
\plref{thmonadlim},
or, more conveniently, its consequences \eqref{plusside} and \eqref{minusside}.
These show that $L_{M^-}^\infty\cap L_{M^+}^\infty=W_{+,\al}\oplus W_{-,\al}
\oplus (V_{+,\al}\cap V_{-,\alpha})$, which by
Lemma \plref{algtoplemma} is also isomorphic to $H^{\rm{\scriptstyle
even}}(M;\C^n_\alpha)$. Notice that the full conclusion of Theorem
\plref{thmonadlim} is used here.

2. and 3. are proven by the same argument. We prove 3.  From the definition of
$V_{+,\alpha}$ \eqref{vee} there is an exact sequence
$$0\to L_{M^+}\cap F^-_0\to L_{M^+}\cap (F^-_0+\ker A)\to V_{+,\alpha}\to
0.$$ It follows easily that for any subspace $V\subset \ker A$ there
is an exact
sequence
\begin{equation}\label{exact1}0\to L_{M^+}\cap F^-_0\to L_{M^+}\cap
(F^-_0+\gamma(V))\to
  V_{+,\alpha}\cap \gamma(V) \to 0.\end{equation}
  Lemma \plref{technical} identifies $L_{M^+}\cap F^-_0$ with the image
  $$H^{\rm{\scriptstyle even}}(M^+,\del M^+;\C^n_\alpha)\to H^{\rm{\scriptstyle
even}}(M^+;\C^n_\alpha),$$  and with $W_{+,\alpha}$.   Thus the dimension of
$L_{M^+}\cap F^-_0$  is  independent of the length of the  collar of
$M^+$.   Corollary
\plref{veehomology} identifies
$ V_{+,\alpha}$ with the image of $H^*(M^+;\C^n_\alpha)\to
H^*(\del M^+;\C^n_\alpha)$, hence its intersection with $\gamma(V)$ is
independent of  the length of the collar as well.  Thus the middle term in the
exact sequence
\eqref{exact1} is isomorphic to
$W_{+,\alpha}\oplus (V_{+,\alpha}\cap\gamma(V))$ and in particular
its dimension is independent of
the length of the collar; this shows that
$  (F_0^-\oplus \gamma(V))\cap L_{M^+}^r$
is independent of $r$ for $r<\infty$.

 Now consider the case when $r=\infty$. In the decomposition \eqref{bigsum2}
of $L^2(\Omega^*_N)$,
$$F^-_0\oplus \gamma(V)= F_\nu^-\oplus d(E^+_\nu)\oplus d^*(E^+_\nu)\oplus
\gamma(V)$$
and so
 using \eqref{plusside}  shows that
$  (F_0^-\oplus \gamma(V))\cap L_{M^+}^\infty$ equals
\[  \big(F_\nu^-\oplus d(E^+_\nu)\oplus d^*(E^+_\nu)\oplus
\gamma(V)\big)\cap
\big(F^+_\nu\oplus (W_{+,\al}\oplus
\gamma(W_{+,\al}^\perp))\oplus
  d(E^-_\nu) \oplus V_{+,\al}\big)\]
\[=W_{+,\al}\oplus (\gamma(V)\cap V_{+,\al}).
\]

Therefore, $(F_0^-\oplus \gamma(V))\cap L_{M^+}^\infty$ equals to
$W_{+,\alpha}\oplus (V_{+,\alpha}\cap\gamma(V))$.

4.  It follows immediately from the definition that the triple index is
additive in the following sense: let $H, H'$ be Hermitian symplectic Hilbert
spaces and $P, Q, R$ (resp. $P', Q', R'$) be projections in $\Gr(H)$
(resp. $\Gr(H')$) such that the triple indices $\tau_\mu(P,Q,R),
\tau_\mu(P',Q',R')$ are well--defined. Then the triple index of $(P\oplus
P',Q\oplus Q',R\oplus R')$ is well--defined in the Hermitian symplectic
Hilbert space $H\oplus H'$ and we have
$$\tau_\mu(P\oplus P',Q\oplus
Q', R\oplus
R')=\tau_\mu(P ,Q , R )+\tau_\mu(P', Q', R').$$

In the decomposition \eqref{bigsum2} we have
\[ \begin{split} \gamma(L_{M^-}^\infty)&=
\gamma(F^-_\nu\oplus  d(E^+_\nu) \oplus
\big(\gamma(W_{-,\al}^\perp)\oplus W_{-,\al}\big)
\oplus V_{-,\alpha})\\
&=  F^+_\nu \oplus d^*(E^-_\nu)\oplus \big(W^\perp_{-,\alpha}\oplus
\gamma(W_{-,\alpha})\big)\oplus \gamma(V_{-,\alpha}),\\
  F_0^+\oplus V&=
F_\nu^+\oplus d^*(E^-_\nu)\oplus d(E^-_\nu)\oplus V,  \text{ and }\\
  L_{M^+}^\infty&= F^+_\nu\oplus \big(W_{+,\al}\oplus
\gamma(W_{+,\al}^\perp)\big)\oplus
  d(E^-_\nu) \oplus V_{+,\al}.\end{split}\]
Using the additivity of $\tau_\mu$  we see that
$
\tau_\mu(\gamma(L_{M^-}^\infty), F_0^+\oplus V, L_{M^+}^\infty)$ equals
\[\begin{split}& \tau_\mu(F^+_\nu,F^+_\nu,F^+_\nu)_{F^-_\nu\oplus F^+_\nu}+
\tau_\mu(d^*(E_\nu^-), d^*(E_\nu^-),W_{+,\al}\oplus
\gamma(W_{+,\al}^\perp))_{d(E^+_\nu)\oplus
d^*(E^-_\nu)}\\
&\quad \quad +\tau_\mu(\gamma(W_{-,\al}^\perp)\oplus W_{-,\al},
d(E^-_\nu),d(E^-_\nu))_{d^*(E^+_\nu)\oplus
d(E^-_\nu)} +\tau_\mu(\gamma(V_{-,\alpha}), V,
V_{+,\alpha})_{\ker A}\end{split}\] which equals $\tau_\mu(\gamma(V_{-,\alpha}),
V, V_{+,\alpha})_{\ker A}$ by Proposition
\plref{vanishing}.

\end{proof}

\begin{proof}[Proof of Theorem \plref{triplebaby}]  Combining Lemma \plref{Lemma-8.9}
and Lemma \plref{triplehomotopy} to the paths of Cauchy data spaces obtained by
stretching the collars of $M^{\pm}$ to their adiabatic limit, we see that
(switching from projection to Lagrangian notation)
 $$\tau_\mu(I-P_{M^-}, P^+(V), P_{M_+})=
\tau_\mu(\gamma(L_{M^-}^\infty), F_0^+\oplus V, L_{M^+}^\infty)=
\tau_\mu(\gamma(V_{-,\alpha}), V, V_{+,\alpha}).$$
It then follows from Theorem \plref{symsplit2} (and in particular
\eqref{eq8.1}) that
$$\etab(D,M)=\etab(D_B,M^+; V\oplus F^+_0)
+\etab(D,M^-;F^-_0\oplus \gamma(V) )-\tau_\mu(\gamma(V_{-,\alpha}), V,
V_{+,\alpha}).$$

Proposition \plref{vanishing} shows that $\tau_\mu(\gamma(V_{-,\alpha}), V,
V_{+,\alpha})=0$ if $V=\gamma(V_{-,\alpha}) $ or $
V_{+,\alpha}$.
\end{proof}

\begin{remark} The Riemannian metric on the separating hypersurface $N$ enters
into the formula of Theorem \plref{triplebaby}  via the map $\gamma:\ker A\to \ker
A$, since $\gamma$ equals the Hodge--$*$ operator up to a power of $i$. It
follows that the correction term $\tau_\mu(\gamma(V_{-,\alpha}), V,
V_{+,\alpha})$ is not a homotopy invariant. In fact, suppose that varying the
Riemannian metric moves the space $\gamma(V_{-,\alpha})$ slightly (we use the
Hodge and DeRham theorems to identify this as a subspace of the fixed space
$H^*(N;\C^n_\alpha)$).  Then by choosing $V\subset H^*(N;\C^n_\alpha)$ carefully
so that $\gamma(V_{-,\alpha})$ passes through $V$ as the metric is varied one
can change $\tau_\mu(\gamma(V_{-,\alpha}), V,
V_{+,\alpha})$. \end{remark}

   To complete our  analysis of the odd signature operator, we derive a formula
which calculates $\etab(D,M)$ in terms of intrinsic invariants of the pieces
and a ``correction term'' which only depends on finite--dimensional
  data, namely the subspaces $V_{\pm,\al}\subset \ker A$.

First we define the analogue of the map $\Phi:\Gr(A)\to
\cu(\mathcal{E}_{i},\mathcal{E}_{-i})$ of Equation
\eqref{phimap} in the finite--dimensional space $\ker A$.
  We use the Lagrangian notation,  so that
to any Lagrangian subspace
$K\subset
\ker A$ we assign the unitary map
\begin{equation}\label{smallphimap}\phi(K):(\mathcal{E}_i\cap \ker
A)\to(\mathcal{E}_{-i}\cap \ker A)\end{equation} by the formula
$$K=\bigsetdef{x+\phi(K)(x)}{x\in (\mathcal{E}_i\cap \ker A)}.$$

\newcommand{\even}{\operatorname{even}}
\begin{theorem}\label{invariantform} For the odd signature operator $D$ coupled
to a flat connection $B$ with holonomy $\alpha:\pi_1M\to U(n)$ acting on a split
manifold
$M=M^+\cup_N M^-$,
\[\begin{split}\etab(D,M)&=\etab(D,M^+;V_{+,\alpha}\oplus F_0^+)+
\etab(D,M^-;F_0^-\oplus V_{-,\alpha})\\
&\quad  +\dim ( V_{+,\alpha}\cap V_{-,\alpha})
-\tfrac{1}{2\pi i}\tr\log(-\phi(V_{+,\alpha})\phi(V_{-,\alpha})^*).
\end{split}\]
\end{theorem}
\begin{remark}   Corollary \plref{veehomology} implies that  $\dim
(V_{+,\alpha}\cap V_{-,\alpha})$ depends only on the homotopy type of the triple
$(M,M^+,M^-)$ and the representation $\al:\pi_1M\to U(n)$.\end{remark}

\begin{proof}
 Theorem \plref{invertible}  implies that
\begin{equation}\label{8.10.1}\begin{split}
 \etab(D,M^+;P^+(\gamma(V_{-,\alpha})))&=\etab(D,M^+;P_{M^+})\\ &\quad +
\tfrac{1}{2\pi i}\tr\log
\big(\Phi(P^+(\gamma(V_{-,\alpha})))\Phi(P_{M^+})^*\big) \end{split}
\end{equation}
and that
\begin{equation}\label{8.10.2}\begin{split}
 \etab(D,M^+;P^+(V_{+,\alpha}))&=
\etab(D,M^+;P_{M^+})\\ &\quad +
\tfrac{1}{2\pi i}\tr\log \big(\Phi(P^+(V_{+,\alpha}))\Phi(P_{M^+})^*\big).
                  \end{split}
\end{equation}
Subtracting \eqref{8.10.2} from \eqref{8.10.1}   yields
\begin{equation}\label{8.10.5}
\begin{split}
 &\etab(D,M^+;P^+(\gamma(V_{-,\alpha})))=
\etab(D,M^+;P^+( V_{+,\alpha}) )  \\
 & +
\tfrac{1}{2\pi i}\tr\log \big(\Phi(P^+(\gamma(V_{-,\alpha})))\Phi(P_{M^+})^*\big)
-\tfrac{1}{2\pi i}\tr\log \big(\Phi(P^+(V_{+,\alpha}))\Phi(P_{M^+})^*\big).
\end{split}\end{equation}

The difference
$$\tfrac{1}{2\pi i}\tr\log
\big(\Phi(P^+(\gamma(V_{-,\alpha})))\Phi(P_{M^+})^*\big)
- \tfrac{1}{2\pi i}\tr\log \big(\Phi(P^+(V_{+,\alpha}))\Phi(P_{M^+})^*\big)$$
is equal to
$$\tau_\mu( P^+(V_{+,\alpha}), P^+(\gamma(V_{-,\alpha})), P_{M^+})
- \tfrac{1}{2\pi i}\tr\log
\big(\Phi(P^+(V_{+,\alpha}))\Phi(P^+(\gamma(V_{-,\alpha})))^*\big)$$
by Lemma \plref{tripleindex}. Moreover, it follows easily from the definitions
that
\[ \tr\log
\big(\Phi(P^+(V_{+,\alpha}))\Phi(P^+(\gamma(V_{-,\alpha})))^*\big)
=\tr\log
\big(\phi(V_{+,\alpha})\phi(\gamma(V_{-,\alpha}))^*\big),\]
and  since $\phi(\gamma(V))=-\phi(V)$,  that
\[\tr\log
\big(\phi(V_{+,\alpha})\phi(\gamma(V_{-,\alpha}))^*\big)=
\tr\log
\big(-\phi(V_{+,\alpha})\phi(V_{-,\alpha})^*\big).\]

Hence \eqref{8.10.5} reduces to

\begin{equation}\label{8.10.3}
\begin{split}
 \etab(D,M^+;P^+(\gamma(V_{-,\alpha})))&=
\etab(D,M^+;P^+( V_{+,\alpha}) )  \\
 & \hskip-1.3in + \tau_\mu( P^+(V_{+,\alpha}), P^+(\gamma(V_{-,\alpha})),
P_{M^+}) -\tfrac{1}{2\pi i}
\tr\log \big(-\phi(V_{+,\alpha})\phi(V_{-,\alpha})^*\big).
\end{split}\end{equation}

We will show that
\begin{equation}\label{8.10.4}
\tau_\mu( P^+(V_{+,\alpha}), P^+(\gamma(V_{-,\alpha})), P_{M^+})
=\dim ( V_{+,\alpha}\cap V_{-,\alpha}).
\end{equation}

Assuming \eqref{8.10.4}, the proof of Theorem \plref{invariantform} is completed
by combining
\eqref{8.10.3} and
 Theorem \plref{triplebaby}, taking $V=\gamma(V_{-,\alpha})$.

It remains, therefore, to prove \eqref{8.10.4}.
The proof is similar to the proof of Theorem \plref{triplebaby}.
Lemma \plref{Lemma-8.9} implies that as the collar of $M^+$ is stretched to its
adiabatic limit, the dimension of the  intersection of $L^r_{M^+}$ with
$F^+_0\oplus V_{+,\alpha}$ is independent of $r\in[0,\infty]$, as is the
dimension of the intersection  of $L^r_{M^+}$ with $F^+_0\oplus
\gamma(V_{-,\alpha})$.

Lemma
\plref{triplehomotopy} then implies that
$\tau_\mu( P^+(V_{+,\alpha}), P^+(\gamma(V_{-,\alpha})), P_{M^+})$
is equal to $\tau_\mu( P^+(V_{+,\alpha}), P^+(\gamma(V_{-,\alpha})),
P^\infty_{M^+})$. Using additivity of the triple index with respect to  the
decomposition
\eqref{bigsum2},    the calculation of $L^\infty_{M^+}$ \eqref{plusside}, and
Proposition \plref{vanishing}, we conclude that
$$\tau_\mu( P^+(V_{+,\alpha}), P^+(\gamma(V_{-,\alpha})),
P^\infty_{M^+})=\tau_\mu(V_{+,\alpha},\gamma(V_{-,\alpha}),
V_{+,\alpha}).$$
Proposition \plref{vanishing} then implies that
\[ \begin{split}\tau_\mu(V_{+,\alpha},\gamma(V_{-,\alpha}),
V_{+,\alpha})&=\dim( V_{+,\alpha}\cap V_{-,\alpha}).\end{split}\]
\end{proof}

It is convenient to introduce the following notation.
\begin{definition}\label{defofm}
 Let $(H,\langle\ , \ \rangle, \gamma)$ be a finite-dimensional
Hermitian symplectic space (cf. Def. \plref{defofsymp}).    Define a
function of pairs of Lagrangian subspaces
$$ m_H:\mathcal{L}(H)\times \mathcal{L}(H) \to \R$$ by
\[\begin{split}m_H(V,W)&=-\tfrac{1}{\pi
i}\tr\log(-\phi(V)\phi(W)^*)+  \dim(V\cap W)\\
&=-\tfrac{1}{\pi i}
\sum_{\begin{array}{c}\SST\lambda\in\spec(-\phi(V)\phi(W)^*)\\
\SST\lambda\neq -1\end{array}} \log \lambda.
  \end{split}
\]
Here $\phi(V)$ is the unitary map from the $+i$ eigenspace $E_i$ of
 $\gamma$ to the $-i$ eigenspace $E_{-i}$ of $\gamma$  whose graph is $V$.
(If $H=0$ then  set $m_H(V,W)=0$.) Recall that $V\cap W$ is isomorphic
to the $-1$--eigenspace of $\phi(\gamma V)\phi(W)^*=-\phi(V)\phi(W)^*$ (cf. Lemma \plref{pairs}).

The function $m$ has been investigated before, the notation is taken from 
\cite{Bun:GPE}.

 Given an even-dimensional Riemannian
manifold
$(N,g)$  and a representation $\alpha:\pi_1N\to U(n)$, define
$$m(V_+,V_-, \alpha,g)=m_{H^*(N;\C^n_\al)}(V_+,V_-),$$
where we have used the Hodge Theorem (and hence the metric $g$ on $N$)
to identify
$H^*(N;\C^n_\alpha)$ with $\ker A$ (so that $\gamma$ and hence
$\phi$ make sense).
\end{definition}

Thus Theorem \plref{invariantform} says that
\begin{equation}\label{dobeedo}\begin{split}\etab(D,M)&=\etab(D, M^+;
V_{+,\al}\oplus F^+_0)+
\etab(D, M^-;F^-_0\oplus V_{-,\al})\\ &\hskip.5in +
\tfrac{1}{2}\dim(V_{+,\al}\cap V_{-,\al})+\tfrac{1}{2} \ m(V_{+,\al},
V_{-,\al},\alpha,g).\end{split}\end{equation}
Using $\eta$--invariants instead of $\etab$--invariants, \eqref{dobeedo}
 can
 be put in the more compact form 
\begin{equation}\label{dobeedo2}
\eta(D,M)=\eta(D, M^+;
V_{+,\al}\oplus F^+_0)+
\eta(D, M^+;F^-_0\oplus V_{-,\al}) +
 m(V_{+,\al},V_{-,\al},\alpha,g)\end{equation}
using  Equation \eqref{exact1} and 
Lemma  \plref{algtoplemma}.

The function $m_H(V,W)$ has some useful properties which we list
in the following proposition.

\begin{proposition}\label{mproperties}\hfill
\begin{enumerate}
\item $m_H(W,V)=-m_H(V,W)$.
\item $m_{H_1\oplus H_2}(V_1\oplus V_2,W_1\oplus
W_2)=m_{H_1}(V_1,W_1) + m_{H_2}(V_2,W_2)$.
\item If $h_t:H\to H, t\in [0,1]$ is a continuous family of
symplectic automorphisms, then $m_H(h_t(V), h_t(W))$ is
continuous  in $t$.\qed
\end{enumerate}
\end{proposition}
\begin{proof}
The first assertion follows immediately from the definition of $m_H$.
The second assertion is clear.  For the 
third, notice that $\dim(h_t(V)\cap h_t(W))$ is independent of $t$, 
and that the $-1$ eigenspace of $ -\phi(h_t(V))\phi(h_t(W))^*$ is isomorphic
to 
$h_t(V)\cap h_t(W)$. In particular the $-1$ eigenspace of
 $ -\phi(h_t(V))\phi(h_t(W))^*$  is constant dimensional, and so 
$t\mapsto \log( -\phi(h_t(V))\phi(h_t(W))^*)$ is continuous. These facts  
imply that $m_H(h_t(V),h_t(W))$ is continuous in $t$.
\end{proof}

\vskip.3in

\subsection{The Atiyah--Patodi--Singer $\rho_\alpha$--invariant for manifolds
with boundary}  We  apply  the previous results to obtain information
about the Atiyah--Patodi--Singer $\rho_\al$--invariant \cite{AtiPatSin:SARII}.
Consider two flat connections: $B$ with holonomy $\alpha$, odd signature
operator $D_B$  and  tangential operator
$A_b$, and the trivial connection $\Theta$ on the bundle $\C^n\times
X\to X$ with (trivial) holonomy
$\tau$, odd signature operator $D_\Theta$,  and tangential operator $A_\theta$.

In expressions like $\eta(D,X;V\oplus F_0^+)$ the notation
$F_0^+$ is to be understood
as the positive eigenspan of the tangential operator $A$ of $D$ and $V$ as
a Lagrangian in $\ker A$. In particular, in a ``mixed'' expression
like $\etab(D_B,X;V_\al\oplus F_0^+ )-
\etab(D_\Theta,X;V_\tau\oplus F_0^+ )$
the reader should understand that the first $F_0^+$ refers to the positive
eigenspan of $A_b$ and the second the positive eigenspan of $A_\theta$.
These are in general unrelated since $A_b$ acts on the bundle
$E|{\partial X}$ and $A_\theta$ acts on the trivial bundle.

\begin{lemma} \label{independence} Let $X$ be a Riemannian manifold with
boundary, whose collar is isometric to $ [0,\ep)\times \partial X$. Let
$B$ be a flat connection on a compact manifold
$X$ in temporal gauge near the boundary with  holonomy $\alpha$, and let
$\Theta$ denote the trivial connection, with trivial holonomy
$\tau:\pi_1(X)\to U(n)$.

    Then the difference
\begin{equation}\label{rhoinvt}\etab(D_B,X;V_\al\oplus F_0^+ )-
\etab(D_\Theta,X;V_\tau\oplus F_0^+ )\end{equation}
depends only on the diffeomorphism type of $X$, the conjugacy class of the
holonomy representation of
$B$ and the restriction of the Riemannian metric to $\partial X$.
\end{lemma}

\begin{proof}
We explained above why the $\eta$--invariant depends   on the flat connection
$B$ only through the  conjugacy class of its holonomy representation.

  By taking the double of $X$ we
obtain a closed Riemannian manifold
$M=X\cup -X=M^+\cup M^-$ over which the connections $B$  and $\Theta$ extend
flatly.

Letting
$D_B$ denote the extension to $M$, we know from \cite{AtiPatSin:SARII}
that  the difference
\begin{equation}\label{bebop}\etab(D_B,M)- \etab(D_\Theta,M)\end{equation}
is independent of the metric on $M$ and depends only on the conjugacy
class of the holonomy representation of $B$ (see the paragraph following this
proof).

 Theorem
\plref{triplebaby}  shows  that
\begin{equation}\label{eq8-12-1}\begin{split}
\etab(D_B,M)- \etab(D_\Theta,M)&=\etab(D_B,M^+; V_\al\oplus F_0^+)+
\etab(D_B,M^-;F_0^-\oplus \gamma(V_\al)) \\
&\quad-\etab(D_\Theta,M^+;
V_\tau\oplus F_0^+)-\etab(D_\Theta,M^-;F^-_0\oplus\gamma(V_\tau)).\end{split}
\end{equation}

Notice that by Corollary \plref{veehomology} the subspaces $V_\alpha$ and $V_\tau$
are independent of the Riemannian metric on $M^+$. Hence
solving for $ \etab(D_B,M^+;V_\al\oplus F_0^+)-\etab(D_\Theta,M^+;
V_\tau\oplus F_0^+)$ in
\eqref{eq8-12-1} yields an expression which is unchanged when the Riemannian
metric is altered on the interior of
$X=M^+$.
\end{proof}

We can now  extend the definition of the Atiyah--Patodi--Singer
$\rho_\al$--invariant to manifolds with boundary.  Recall that
the $\rho_\al$--invariant is defined in \cite{AtiPatSin:SARII} for a
closed manifold $M$  and a representation $\al:\pi_1(M)\to
U(n)$ by
$$\rho (M,\al)=\eta(D_B)-\eta(D_\Theta)$$
where $B$ is a flat connection on $M$ with holonomy $\al$ and
$\Theta$ denotes a trivial $U(n)$ connection. It is
a diffeomorphism invariant of the pair $(M,[\al])$ where
$[\al]$ denotes the conjugacy class of $\al$.  In terms of
reduced $\eta$-invariants $\rho(M,\al)$ can be written:
$$\rho (M,\al)= 2(\etab(D_B)-\etab(D_\Theta))-\dim\ker D_B +
\dim \ker D_\Theta.$$
Since $\ker D_B$ is isomorphic to $H^{\text{\tiny even}}(M;\C^n_\al)$
this is the
same as
   \begin{equation}
\label{rhomclosed}
\rho(M,\al)= 2(\etab(D_B,M)-\etab(D_\Theta,M))-\dim H^{\text{\tiny
even}}(M;\C^n_\al) +
\dim   H^{\text{\tiny even}}(M;\C^n).\end{equation}

\begin{definition}\label{rhoXnotclosed} Given a triple
$(X,\al,g)$, where
\begin{enumerate}
\item $X$ is a compact odd--dimensional manifold with boundary,
\item $\al:\pi_1(X)\to U(n)$ is a representation, and
\item $g$ is a Riemannian metric on $\del X$,
\end{enumerate}
choose   a Riemannian
metric on $X$  isometric to $[0,\epsilon)\times \del X $ on a collar of $\del
X$ and a flat connection $B$ with holonomy $\al$  in temporal gauge near
the boundary. Then define
\[  \rho(X,\al,g) := \eta(D_B,X,F_0^+\oplus  V_\al )-
\eta(D_\Theta,X,F_0^+\oplus
  V_\tau). \]
\end{definition}

Reversing the
orientation of $X$ changes the sign of $\rho(X,\al,g)$, since the
$\eta$--invariant changes sign when the orientation is reversed.

In terms of reduced $\eta$ invariants $\rho(X,\al,g)$ can be expressed as:
\[ \rho(X,\al,g) =2\left(\etab(D_B,X,F_0^+\oplus  V_\al )-
\etab(D_\Theta,X,F_0^+\oplus
  V_\tau)\right)
-\dim W_\al +\dim W_\tau,\]
where
$$W_\al\cong\im\bigl( H^{\text{\tiny even}}(X,\partial
X;\C^n_\al)\to H^{\text{\tiny even}}(X;\C^n_\al)\bigr) $$
and
$$W_\tau\cong \im \bigl(H^{\text{\tiny even}}(X,\partial
X;\C^n )\to H^{\text{\tiny even}}(X;\C^n )\bigr).$$
This is because the kernel of $D_B$ acting on $X$ with boundary conditions given
by the Atiyah--Patodi--Singer projection $P^+(V_\al)$  is isomorphic to
$W_\al$ by Lemma \plref{technical} and Equation \eqref{exact1} (with
$V=V_{+,\al}$), and similarly for $D_\Theta$.

Lemma \plref{independence} shows that   $\rho(X,\al,g)$ is independent
of the choice of Riemannian metric on the interior of $X$ (as long as the
metric is a product in  some collar of the boundary) and the choice of flat
connection $B$ with holonomy $\alpha$.

\vskip.3in
We now can state the main result of this section.

\begin{theorem}\label{there}
Suppose the closed manifold $M$ contains a hypersurface $N$ separating
$M$ into $M^+$ and $M^-$. Fix a Riemannian metric $g$ on $N$.  Suppose
that
$\al:\pi_1(M)\to U(n)$ is a representation, and let $\tau:\pi_1(M)\to
U(n)$ denote
the trivial  representation.

Then
\begin{equation*}
\rho (M,\al)
    =\rho(M^+,\al,g)+\rho(M^-,\al,g)
   +
m(V_{+,\al},
V_{-,\al},\alpha,g)-m(V_{+,\tau},
V_{-,\tau},\tau,g).
\end{equation*}
\end{theorem}
\begin{proof}

This follows by
applying Equation \eqref{dobeedo2}  to $B$ and $\Theta$ and subtracting.
\end{proof}

It can be shown that the invariants $\rho(M^\pm,\al,g)$ and $m(V_+,V_-,\al,g)$
depend in general on the choice of Riemannian metric $g$ on  the hypersurface
$N$. We leave an an intriguing open problem to determine exactly how they depend
on the metric $g$, and in particular, how these invariants change if $g$ is
replaced by the pulled--back metric $f^*(g)$  for a diffeomorphism $f:N\to N$.

\subsection{Relationship to Wall's non-additivity theorem}

Theorems \plref{there}  and \plref{invariantform} should be viewed as
odd--dimensional counterparts of Wall's non--additivity theorem for the
signature
\cite{Wal:NAS}. Indeed these theorems give formulas which express the
``non-additivity of the signature defect''.  The relationship between
splitting theorems for the
$\eta$--invariant and Wall non--additivity  is explored in Bunke's article
\cite{Bun:GPE} and also in
\cite{HMM:SFM}.

To clarify the relationship between Wall's theorem and Theorem \plref{there},
 consider the
following situation. Suppose we are given  two $4k$--dimensional manifolds
$Z^+$ and
$Z^-$  with
$\partial
Z^\pm=M^\pm\cup_N M^0$. Suppose that $\del M^0=N$ and that $Z=Z^+\cup_{M^0} Z^-$.
Finally suppose that $\al:\pi_1Z\to U(n)$ is a representation and let
$\tau:\pi_1Z\to U(n)$ denote the trivial representation. The
Atiyah--Patodi--Singer signature theorem
\cite[Theorem 2.4]{AtiPatSin:SARII} says that
$$\text{Sign}_\tau(Z)-\text{Sign}_\al(Z)=\rho(M,\al).$$  Similarly
$\text{Sign}_\tau(Z^+)-\text{Sign}_\al(Z^+)=\rho(M^+\cup M^0,\al)$
and
$\text{Sign}_\tau(Z^-)-\text{Sign}_\al(Z^-)=\rho(-M^0\cup M^-,\al)$.
On the other hand Wall's theorem says that
$$\text{Sign}_\al(Z)=\text{Sign}_\al(Z^+)+\text{Sign}_\al(Z^-)
-\sigma(V_{+,\al},V_{-,\al},V_{0,\al},\al),$$
where $\sigma$ is a correction term   which depends on
the relative positions of the subspaces $V_{+,\al},V_{-,\al}$ and $V_{0,\al}$
in
$H^*(N;\C^n_\al)$. Similarly
$
\text{Sign}_\tau(Z)=\text{Sign}_\tau(Z^+)+\text{Sign}_\tau(Z^-)-
\sigma(V_{+,\tau},V_{-,\tau},V_{0,\tau},\tau).
$

Hence
\begin{equation}\label{signature}\begin{split}
\sigma(V_{+,\al},V_{-,\al},V_{0,\al},\al)-
\sigma(V_{+,\tau},V_{-,\tau},V_{0,\tau},\tau)\hskip-1.8in&\\
&=\text{Sign}_\tau(Z)-\text{Sign}_\tau(Z^+)-\text{Sign}_\tau(Z^-)\\
&\quad
-\big(\text{Sign}_\al(Z)-\text{Sign}_\al(Z^+)-\text{Sign}_\al(Z^-)\big)\\
&=\rho(M^+\cup M^-,\alpha)-\rho(M^+\cup M^0,\al)-\rho(-M^0\cup
M^-,\al)\end{split}\end{equation}
Applying Theorem   \plref{there} we see that \eqref{signature}
is equal to  $\tilde{\sigma}_\al- \tilde{\sigma}_\tau$, where
\[\tilde{\sigma}_\al :=m(V_{+,\al},V_{-,\al},\al,g)-m(V_{+,\al},V_{0,\al},\al,g)
-m(V_{0,\al},V_{-,\al},\al,g)\]
and
\[\tilde{\sigma}_\tau :=m(V_{+,\tau},V_{-,\tau},\tau,g)-m(V_{+,\tau},V_{0,\tau},\tau,g)
-m(V_{0,\tau},V_{-,\tau},\tau,g).\]

This motivates introducing the following notation. Given a
Hermitian symplectic space $H$, define the function of triples of
Lagrangian subspaces
\[\tsig_H:\mathcal{L}(H)\times
\mathcal{L}(H)\times \mathcal{L}(H)\to \Z\]
 by
\begin{equation}\label{defoftsig}\tilde{\sigma}_H(V,W,U):=m_H(V,W)+m_H(W,U)
+m_H(U,V).\end{equation}
By definition
$\tsig_\alpha=\tsig_{H^*(N;\C^n_\al)}(V_{+,\al},V_{-,\al},V_{0,\al})$
and similarly for $\tsig_\tau$.
 That $\tsig_H$ is an integer can be seen
by exponentiating
 and using the multiplicativity of the determinant:
 \begin{equation*}
 \begin{split}\exp(2\pi i
\ \tsig_H(V,W,U))&\\
&\hskip-1in =\big(\exp(\tr\log(-\phi(V)\phi(W)^*)+
 \tr\log(-\phi(W)\phi(U)^*)+\tr\log(-\phi(U)\phi(V)^*)\big)^2\\
&=\det\big(
(-1)^3\phi(V)\phi(W)^*\phi(W)\phi(U)^*\phi(U)\phi(V)^*)^2\big)\\
&=1.
\end{split}
\end{equation*}

\begin{proposition}\label{tsigproperties}
The function $\tsig_H$ satisfies the following properties.
\begin{enumerate}
\item
Given a permutation $\beta$,
 $\tsig_H(V_{\beta(1)},V_{\beta(2)},V_{\beta(3)})=\text{\rm sign}(\beta)\tsig_H(V_1,V_2,V_3)$.
\item $\tsig_{H_1\oplus H_2}(V_1\oplus V_2, W_1\oplus W_2,
U_1\oplus U_2)=\tsig_{H_1}(V_1,W_1,U_1)+\tsig_{H_2}(V_2,W_2,U_2)$.

\item If $h:H\to H$ is a symplectic automorphism, then
$\tsig_H(h(V),h(W),h(U))=\tsig_H(V,W,U)$.

\item Take $H=\C^2$ with $\gamma=\begin{pmatrix}
0&-1\\1&0\end{pmatrix}$. Then
$\tsig_H(\C(1,0),\C(1,1),\C(0,1))=1.$
\end{enumerate}
\end{proposition}
\begin{proof} The first and second assertions follow from the
first and second assertions of Proposition \plref{mproperties}.

 For
the third assertion,  we first claim that the group $Sp(H)$ of symplectic
automorphisms  of $H$ is path connected.  To see this, fix a
Lagrangian subspace $V$ of $H$. The map $Sp(H)\to \mathcal{L}(H)$
taking $g$ to $g(V)$ is a fibration with fiber the subgroup
$Sp_V(H)$ consisting of those symplectic automorphisms which leave
$V$ invariant. Next, $Sp_V(H)$ fibers over $GL(V)$ by mapping
$g\in Sp_V(H)$ to the restriction $g|_V$.  The fiber $F$ of this map
 consists of those symplectic transformations $g$
so that $g$ restricts to the identity on $V$.  Writing $H=V\oplus
\gamma(V)$ it is easy to see that $F$ consists of all matrices of the
form \[\begin{pmatrix}
I&A\\0&I
\end{pmatrix}\]
with $A$ an arbitrary real matrix. Thus $F$ is contractible, and
since $GL(V)$ is path connected,
$Sp_V(H)$ is also path connected. Finally, since $Sp(H)$ fibers
over the path connected space $\mathcal{L}(H)\cong U(n)$ with path
connected fiber $Sp_V(H)$, it is itself path connected.

Choose a path $h_t$ from the identity matrix to $h$.   The third
assertion of Proposition \plref{mproperties} shows that $m_H(h_t(V),
h_t(W))$ varies continuously in $t$. Thus the same is true of the
integer-valued function $t\mapsto \tsig_H(h_t(V),h_t(W),h_t(U))$.
Therefore this function is constant, completing the proof of the
third assertion.

To prove the last fact, Notice that $\C^2=E_i\oplus E_{-i}$, where
$E_i$ is the span of $(1,-i)$ and $E_{-i}$ is the span of $(1,i)$.
It is easy to calculate that with respect to these bases,
$$\phi(\C(1,0))=1, \ \phi(\C(1,1))=-i, \ \text{ and }\
\phi(\C(0,1))=-1.$$
Thus
\[\begin{split}
\tsig_H(\C(1,0),\C(1,1),\C(0,1))&=-\tfrac{1}{\pi i}(\log(-1\cdot
i) + \log(-(-i)\cdot(-1))+\log(-(-1)(1)))\\
&=1
\end{split}
\]
\end{proof}

It follows from Proposition \plref{tsigproperties} and
 \cite[Theorem 8.1]{CapLeeMil:OMI} (suitably generalized to the 
complex Hermitian case)
that $\tsig_H$
 is equal to the Maslov triple index $\tau_H$ defined in loc. cit.
Therefore,  $\tsig_H$
 depends only on the underlying symplectic form
 $\omega(x,y)=\langle x,\gamma y\rangle$, and not on the Hermitian
 metric. In particular $\tsig_\al$ and $\tsig_\tau$ are
 independent of the Riemannian metric on $N$.

$\tsig_H$ and the Maslov triple index $\tau_\mu$ defined in Section
\plref{secsixtwo} are (of course) intimately related. $\tau_\mu$
is, up to normalization, what Bunke \cite{Bun:GPE} called the twisted
Maslov triple index. A direct calculation shows the following:
\begin{align}
  \tsig_H(V,W,U)&=-\tau_\mu(V,W,U)-\tau_\mu(\gamma V,W,U)-\tau_\mu(V,\gamma W,U)\\
          &\quad-\tau_\mu(V,W,\gamma U)+\dim(V\cap W)+\dim(W\cap U)+\dim(V\cap U),\nonumber \\
  \tau_\mu(V,W,U)&=\frac 14\bigl(\tsig_H(V,W,U)-\tsig_H(\gamma V, W, U)-\tsig_H(V,\gamma W, U)\\
      &\quad-\tsig_H(V,W,\gamma U)+2\dim(\gamma V\cap W)+2\dim(W\cap \gamma U)+2\dim(V\cap \gamma U)\bigr).\nonumber
\end{align}

Using Proposition 8.2 of loc. cit. we conclude that $\tsig_\al$
equals  Wall's correction term $\sigma(V_{+,\al}, V_{-,\al},
V_{0,\al})$ and similarly $\tsig_\tau$ equals $\sigma(V_{+,\tau},
V_{-,\tau}, V_{0,\tau})$.

Using the argument outlined   in \cite{HMM:SFM} one can analyze
Wall's theorem using  Theorem \plref{invariantform} as follows.

 In the
context described above, the Atiyah--Patodi--Singer signature theorem states that
$$\text{Sign}_\al(Z)=n\int_Z L - \eta(D_B,M^+\cup M^-),$$
where $L$ denotes the $L$--polynomial of the Riemannian curvature tensor on $Z$.
Similarly one obtains formulas for $\text{Sign}_\al(Z^\pm)$
$$\text{Sign}_\al(Z^+)=n\int_{Z^+} L -  \eta(D_B,M^+\cup M^0) $$
and
$$\text{Sign}_\al(Z^-)=n\int_{Z^+} L -  \eta(D_B,-M^0\cup M^-).$$
Applying Theorem \plref{invariantform} and using Lemma \plref{tripleindex}  as
before, one obtains
\[\text{Sign}_\al(Z) -
\text{Sign}_\al(Z^+)-\text{Sign}_\al(Z^-)
 =\tilde{\sigma}_\al
+n\big(\int_Z L-  \int_{Z^+} L- \int_{Z^-} L\big).\]

At this point one can  invoke Wall's theorem and the
identification  of $\tsig_\al$  with Wall's correction term given
above to conclude that the integrals cancel.  (This is not
immediate since the Riemannian metrics on $Z^+$ and $Z^-$ need to be 
in cylindrical form near the boundary to apply the Atiyah--Patodi--Singer
 theorem, but these  do not
glue to give a smooth metric on $Z$ in cylindrical form near the
boundary.)

On the other hand,  Wall's theorem can be proven by showing  that
the integrals cancel.
 This is discussed in   \cite{HMM:SFM} and so one obtains  an analytic proof of
 Wall's theorem.

 More importantly, Equation \eqref{defoftsig} establishes
  a  direct relationship between the correction terms $\tsig_\al$
 and $\tsig_\tau$
 for the non-additivity of the signature to the correction term
 $m(V_{+,\al}, V_{-,\al},\al,g)-m(V_{+,\tau}, V_{-,\tau},\tau,g)$
 for the non-additivity of the $\rho$ invariant.

\subsection{Adiabatic stretching and general Dirac operators.}
Some of the preceding exposition  for the odd signature  operator extends to  the
more general context of arbitrary Dirac operators, and we discuss   aspects of
this now. The new feature of this approach is that the role of adiabatic
stretching in the splitting formula for the $\eta$--invariant is clarified.

Suppose we are given an arbitrary  Dirac operator $D$ on a split manifold
$M=M^+\cup_N
M^-$. Assume as usual that $D=\gamma(\frac{d}{dx} +A)$ on a
collar of $N$.  Let $M_r$ denote the  manifold obtained by replacing the collar
$[-1,1]\times N$ of $N$ by the stretched collar $[-r,r]\times N$.  Thus $M_0=M$.
Given $0\leq r<\infty$, let
$L_{M^\pm}^r$ denote the Cauchy data space of the   operator $D$ acting
on $M^\pm_r=M^\pm\cup([-r,0]\times N)$, and let $L_{M^\pm}^\infty$ denote the
adiabatic limit  $\lim\limits_{r\to\infty}L^r_{M^\pm}$.  Lemma 3.2 of
\cite{DanKir:GSF} states that the path $[0,\infty]\to \Gr(A)$ given by $r\mapsto
L^r_{M^\pm}$ is continuous.

 We let $F_\nu^\pm$ denote the span of
$\lambda$--eigenvectors of $A$ for $\pm \lambda>\nu$, and $E^\pm_\nu$ the
span of
$\lambda$--eigenvectors of $A$ for $0<\pm \lambda\leq\nu$, so that the
$L^2$--sections over $N$ decompose as
$$F^-_\nu\oplus E^-_\nu\oplus \ker A\oplus E^+_\nu\oplus F^+_\nu $$
or, as a symplectic direct sum
\begin{equation}\label{8-A1}
(F^-_\nu\oplus F^+_\nu)\oplus(E^-_\nu\oplus E^+_\nu)\oplus \ker A.
\end{equation}

Theorem \plref{thmonadlim} has a counterpart  for general Dirac operators,
but the conclusion is slightly weaker. The following theorem has a
similar but  simpler proof than Theorem \plref{thmonadlim}. It is
  implied by Theorem 6.5 of \cite{DanKir:GSF}.

\begin{theorem}\label{thmonadlimgen} Let
$V_+\subset
\ker A$ denote the limiting values of extended $L^2$ solutions on $M^+$, so
$V_+=\proj_{\ker A}(L_{M^+}\cap (F^-\oplus \ker A)).$ Let $\nu\ge 0$ be a
number in the non-resonance range of $D$, i.e. $L_{M^+}\cap  F^-_\nu=0$.
Then there exists a subspace $W_+\subset E^-_\nu$ isomorphic to the space of
$L^2$ solutions to $D\beta=0$ on $M^+_\infty$ so that letting
$W_+^\perp\subset E^-_\nu$ denote the orthogonal complement of $W_+$ in
$E^-_\nu$,
$$L_{M^+}^\infty=F^+_\nu\oplus (W_+\oplus
\gamma(W_+^\perp))
\oplus   V_+$$ in the decomposition \eqref{8-A1}.  Moreover,
$\gamma(L^\infty_{M^+})\cap L_{M^+}=0$.
\qed
\end{theorem}

 Then we have the following theorem.

\begin{theorem}\label{generaldirac} With notation as above,  for any $r_0\ge 0$,
\[\begin{split} \etab(D,M_{r_0})-\etab(D,M^+;L_{M^+}^\infty )-
\etab(D, M^-; \gamma(L_{M^+}^\infty))\quad\quad \\
\quad\quad =
\Mas(L_{M^-}^r,L_{M^+}^\infty)_{r\in [r_0,\infty]}
-\Mas( L^r_{M^-},L^r_{M^+})_{r\in [r_0,\infty]}
 \end{split}\]

\end{theorem}
\begin{remark}\label{rem8} In light of Theorems \plref{ML-S6.5} and
\plref{ML-S6.6} the  term  $\Mas(L_{M^-}^r, L_{M^+}^r)_{r\in [r_0,\infty]}$ in
Theorem
\plref{generaldirac} can be thought of as the spectral flow of the family of
operators on $M$ obtained by stretching the collar from $r_0$ to infinity.
Similarly   the term $\Mas( L_{M^-}^r, L_{M^+}^\infty)_{r\in [r_0,\infty]}$ can
be thought of as the spectral flow of the family on $M^-$ obtained by using the
projection to $L^\infty_{M^+}$ as boundary conditions and stretching the collar
of $M^-$ from $r_0$ to infinity.
\end{remark}

\begin{proof} We prove this for $r_0=0$, the general case is obtained by
reparameterizing. Let $P_t$ denote the path of projections to the Cauchy data
space
$L_{M^+}^r$, where $t=1/(r+1)$.
Applying Theorem \plref{symsplit} and Proposition \plref{splitbundles} to    the
path
$P_t$ we see that
$\etab(D,M_{r_0})-\etab(D,M^+;L_{M^+}^\infty )-
\etab(D, M^-; \gamma(L_{M^+}^\infty))$ equals $\SF(D_{P_t},M^+)  +
\SF(D_{I-P_t},M^-)_{[0,1]}$, which by Theorem \plref{ML-S6.5} equals
\begin{equation}\label{8.3.10}
\Mas( {P_t}, P_{M^+} )  + \Mas(P_{M^-}, {I-P_t}) .
\end{equation}
Switching to Lagrangian notation and parameterizing by $r$ instead of $t$ we can
rewrite
\eqref{8.3.10} as
\begin{equation}\label{8.3.1}
-\Mas(\gamma(L^r_{M^+}), L_{M^+} )_{r\in[0,\infty]} -
\Mas( L_{M^-},L^r_{M^+})_{r\in[0,\infty]}.
\end{equation}

We use the homotopy invariance of the Maslov index to simplify these terms.
Consider first $\Mas(\gamma(L^r_{M^+}), L_{M^+} )_{r\in[0,\infty]}$. We will
show this term vanishes.

 The path
$r\mapsto \gamma(L^r_{M^+})$ is homotopic to the composite of $r\mapsto
\gamma(L^r_{M^+})$ and the constant map at $  \gamma(L^\infty_{M^+})$, and the
constant path at $L_{M^+}$ is homotopic to the  composite of
$r\mapsto  L^r_{M^+} $ and its inverse.  Since $\gamma(L^r_{M^+})\cap
L^r_{M^+}=0$ for all $r$, $\Mas(\gamma(L^r_{M^+}), L^r_{M^+} )=0 $ and
so by path additivity of the Maslov index,
$\Mas(\gamma(L^r_{M^+}), L_{M^+} )= -  \Mas(\gamma(L^\infty_{M^+}), L^r_{M^+} ).$
 Theorem \plref{thmonadlimgen} says that $\gamma(L^\infty_{M^+})\cap L^r_{M^+}=0$
for
$r=0$, but by reparameterizing we see that the intersection is zero for all
$r<\infty$; obviously $\gamma(L^\infty_{M^+})\cap L^\infty_{M^+}=0$. Hence
$\Mas(\gamma(L^\infty_{M^+}), L^r_{M^+} )=0$ and so
\begin{equation}\label{8.3.2}\Mas(\gamma(L^r_{M^+}), L_{M^+}
)_{r\in[0,\infty]}=0.\end{equation}

Consider now the term $\Mas( L_{M^-},L^r_{M^+})_{r\in[0,\infty]}.$  The constant
path at $L_{M^-}$ is homotopic to the composite of $r\mapsto L_{M^-}^r$ and its
inverse, and the path $r\mapsto L^r_{M^+}$ is homotopic to its composite with
the constant path at $L^\infty_{M^+}$. Therefore,
\begin{equation}\label{8.3.3}\Mas( L_{M^-},L^r_{M^+})= \Mas(L_{M^-}^r,L_{M^+}^r)
-
\Mas(L_{M^-}^r,L_{M^+}^\infty).\end{equation}
Substituting \eqref{8.3.2} and \eqref{8.3.3} into \eqref{8.3.1} finishes the
proof.
 \end{proof}

We finish this article by outlining a few ways to use Theorem
\plref{generaldirac} to obtain   other useful splitting formulas for the
$\eta$--invariant.  We will not give an exhaustive list, but
we note  that many other useful formulas can be derived from these using the
results of Sections
\plref{sec5}, \plref{sec6}, and
\plref{sec7}.  One can, of course, obtain other formulas by reversing the roles
of
$M^+$ and $M^-$ in Theorem   \plref{generaldirac} and in these examples.

\begin{example} Suppose that $L^\infty_{M^-}\cap L^\infty_{M^+}=0$. Then there
exists an $r_0\ge 0$ so that
$L^r_{M^-}\cap L^r_{M^+}=0$ and $L^r_{M^-}\cap L^\infty_{M^+}=0$ for all $r\ge
r_0$. Applying Theorem \plref{generaldirac} we see that if $r\ge r_0$ then
\[ \begin{split}\etab(D,M_{r}) &=\etab(D,M^+;L^\infty_{M^+} )+
\etab(D, M^-;   \gamma( L^\infty_{M^+}))\\
&\hskip-.5in  =\etab(D,M^+;F^+_\nu\oplus W_+\oplus \gamma(W_+^\perp)\oplus V_+
)+
\etab(D, M^-;  F^-_\nu\oplus W^\perp_+\oplus \gamma(W_+)\oplus \gamma(V_+))
. \end{split}\]
The hypothesis $L^\infty_{M^-}\cap L^\infty_{M^+}=0$ is a technically simpler
replacement of the hypothesis ``no exponentially small eigenvalues'' which
appears in related results in the literature.
\end{example}

\begin{example} Suppose that $D\beta=0$ has no $L^2$ solutions   on
$M^+_\infty $; i.e.  that $W^+=0$ in Theorem
\plref{thmonadlimgen}. Then $L^\infty_{M^+}=F^+\oplus V_+$ and so
\[\begin{split} \etab(D,M)-\etab(D,M^+;F^+\oplus V_+ )-
\etab(D, M^-; F^-\oplus \gamma(V_+))\quad\quad \\
\quad\quad =
\Mas(L_{M^-}^r ,F^+\oplus V_+)_{r\in [ 0,\infty]}
-\Mas( L^r_{M^-},L^r_{M^+})_{r\in [0,\infty]}.
 \end{split}\]
In other words, with respect to the Atiyah--Patodi--Singer boundary conditions
given by the projection to $F^+\oplus V_+$ on $M^+$ and the projection to
$F^-\oplus \gamma(V_+)$ on $M^-$, the failure of the additivity of the
$\etab$--invariants is measured by $\Mas(L_{M^-}^r,F^+\oplus V_+)_{r\in
[0,\infty]} -\Mas( L^r_{M^-},L^r_{M^+})_{r\in [0,\infty]}$.  As remarked
above this is the difference of the spectral flow of $D$ on $M^-$ with $F^+\oplus
V_+$ conditions as the collar of $M^-$ is stretched to infinity, and the
spectral flow of $D$ on $M$ as the collar is stretched to infinity.
\end{example}
\begin{example} We can combine the previous two examples as follows. Suppose that
there are no $L^2$ solutions on $M_\infty^+$ and $M_\infty^-$ (i.e. $W_+=0=W_-$)
and that the limiting values of extended $L^2$ solutions are transverse  (i.e.
$V_+\cap V_-=0$ in $\ker A$;  this happens for example if $\ker A=0$).  Then
$L^\infty_{M^\pm}= F^{\pm}\oplus V_{\pm}$ and so both of the previous examples
apply.

Hence there exists an
$r_0\ge 0$ so that
$L^r_{M^-}\cap L^r_{M^+}=0$ and $L^r_{M^-}\cap L^\infty_{M^+}=0$ for all $r\ge
r_0$.  Therefore,
\[\begin{split} \etab(D,M)-\etab(D,M^+;F^+\oplus V_+ )-
\etab(D, M^-; F^-\oplus \gamma(V_+))\quad\quad\quad\quad \quad\quad\\
\quad\quad\quad\quad\quad\quad =
\SF(D,  M^-_r ;F^+\oplus V_+)_{r\in [0,r_0]}
-\SF( D, M_{r})_{r\in [0,r_0]}.
 \end{split}\]
This says that the failure of additivity of the $\etab$--invariants with
Atiyah--Patodi--Singer boundary conditions is measured by the difference of the spectral flow of $D$ on $M^-$ with $F^+\oplus
V_+$ conditions as the length of the collar of $M^-$ is stretched to $r_0$,
and the spectral flow of $D$ on $M$ as the collar is stretched to $r_0$.

In particular, if $r\ge r_0$,
\[  \etab(D,M_r)=\etab(D,M^+;F^+\oplus V_+ )+
\etab(D, M^-; F^-\oplus \gamma(V_+)).\]
This last formula appears in Bunke's article \cite{Bun:GPE}.  The reader should
compare this formula with the formula of Theorem \plref{triplebaby} (with
$V=V_{+,\al}$) which,  by contrast, holds in complete generality for the odd
signature operator.
\end{example}

These examples, together with Theorem \plref{thmonadlimgen}, underscore the
point that   difficulties in establishing simple splitting formulas
for the
$\eta$--invariant using Atiyah--Patodi--Singer boundary conditions  arise
from  the existence of
$L^2$ solutions on the two parts of the decomposition of $M$. To put this in a
positive perspective, the failure of the additivity of the
$\eta$--invariant with Atiyah--Patodi--Singer boundary conditions is measured by
the spectral flow terms discussed in Remark \plref{rem8} and  symplectic
invariants of the Lagrangian subspaces
$W_\pm\oplus
\gamma(W^\perp_\pm)$ in the finite--dimensional symplectic space $E^-_\nu\oplus
E^+_\nu$ consisting of the span of the
$\mu$--eigenvectors of $A$ with $-\nu\leq\mu\leq\nu$, $\mu\ne 0$. In our analysis
of the odd signature operator the formulas simplify because we can control these
terms; the spectral flow terms vanish   for topological reasons and the
symplectic invariants  of  the Lagrangian subspaces
$W_\pm\oplus
\gamma(W^\perp_\pm)$ vanish  because of the additional control on
$W_\pm$ that Theorem \plref{thmonadlim} provides over   Theorem
\plref{thmonadlimgen}.

\bibliography{mlabbr,books99,papers99,local.bib,paperPaul0512.bbl}

\begin{thebibliography}{10}

\bibitem{AtiPatSin:SARI}
\textsc{M.~F. Atiyah}, \textsc{V.~K. Patodi}, \protect\BIBand{} \textsc{I.~M.
  Singer}: \emph{Spectral asymmetry and Riemannian geometry I}.
\newblock Math. Proc. Camb. Phil. Soc. \textbf{77} (1975), 43--69

\bibitem{AtiPatSin:SARII}
\textsc{M.~F. Atiyah}, \textsc{V.~K. Patodi}, \protect\BIBand{} \textsc{I.~M.
  Singer}: \emph{Spectral asymmetry and Riemannian geometry II}.
\newblock Math. Proc. Camb. Phil. Soc. \textbf{78} (1975), 405--432

\bibitem{AvrSeiSim:IPP}
\textsc{J.~Avron}, \textsc{R.~Seiler}, \protect\BIBand{} \textsc{B.~Simon}:
  \emph{The index of a pair of projections}.
\newblock J. Funct. Anal. \textbf{120} (1994), 220--237

\bibitem{BooFur:MIF}
\textsc{B.~Booss-Bavnbek} \protect\BIBand{} \textsc{K.~Furutani}: \emph{The
  Maslov index: a functional analytical definition and the spectral flow
  formula}.
\newblock Tokyo J. Math. \textbf{21} (1998), 1--34

\bibitem{BooWoj:EBP}
\textsc{B.~Boo{\ss}-Bavnbek} \protect\BIBand{} \textsc{K.~P. Wojciechowski}:
  \emph{Elliptic Boundary Problems for Dirac Operators}.
\newblock Birkh\"auser, Basel (1993)

\bibitem{BruLes:BVPI}
\textsc{J.~Br{\"u}ning} \protect\BIBand{} \textsc{M.~Lesch}: \emph{On boundary
  value problems for Dirac type operators: I. Regularity and
  self--adjointness}.
\newblock Preprint, 55p. (1999).
\newblock Math.FA/9905181

\bibitem{BruLes:EIN}
\textsc{J.~Br{\"u}ning} \protect\BIBand{} \textsc{M.~Lesch}: \emph{On the
  eta--invariant of certain non--local boundary value problems}.
\newblock Duke Math. J. \textbf{96} (1999), 425--468.
\newblock Dg-ga/9609001

\bibitem{BruLes:STB}
\textsc{J.~Br{\"u}ning} \protect\BIBand{} \textsc{M.~Lesch}: \emph{Spectral
  theory of boundary value problems for Dirac type operators}.
\newblock In: B.~Boo{\ss}-Bavnbek \protect\BIBand{} K.~P. Wojciechowski (eds.),
  \emph{Geometric Aspects of PDE - Spectral Invariants (Roskilde Conf. Sept
  1998)}, vol. 242 of \emph{Contemp. Math.} AMS, Providence, RI (1999), pp.
  203--215.
\newblock Math.DG/9902100

\bibitem{Bun:GPE}
\textsc{U.~Bunke}: \emph{On the gluing problem for the $\eta$--invariant}.
\newblock J. Diff. Geom. \textbf{41} (1995), 397--448

\bibitem{CapLeeMil:OMI}
\textsc{S.~Cappell}, \textsc{R.~Lee}, \protect\BIBand{}
\textsc{E.~Miller}:
  \emph{On the Maslov index}.
\newblock Comm. Pure Appl. Math. \textbf{47} (1994), 121--186

\bibitem{CapLeeMil:SA2}
\textsc{S.~Cappell}, \textsc{R.~Lee}, \protect\BIBand{}
\textsc{E.~Miller}:
  \emph{Self-adjoint elliptic operators and manifold decompositions II:
  Spectral flow and Maslov index}.
\newblock Comm. Pure Appl. Math. \textbf{49} (1996), 869--909

\bibitem{DaiFre:EID}
\textsc{X.~Dai} \protect\BIBand{} \textsc{D.~S. Freed}:
  \emph{$\eta$--invariants and determinant lines}.
\newblock J. Math. Phys. \textbf{35} (1994), 5155--5194

\bibitem{Dan:ETN}  \textsc{M.~Daniel}:
  \emph{An extension of a theorem of Nicolaescu on
spectral flow and the Maslov index}.
\newblock Proc. Amer. Math. Soc. \textbf{128}
(2000), no. 2, 611--619 


\bibitem{DanKir:GSF}
\textsc{M.~Daniel} \protect\BIBand{} \textsc{P.~Kirk}: \emph{A general
  splitting formula for the spectral flow (with an appendix by K. P.
  Wojciechowski)}.
\newblock Mich. Math. J. \textbf{46} (1999), 589--617

\bibitem{DouWoj:ALE}
\textsc{R.~G. Douglas} \protect\BIBand{} \textsc{K.~P. Wojciechowski}:
  \emph{Adiabatic limits of the $\eta$-invariants the odd--dimensional
  Atiyah--Patodi--Singer problem}.
\newblock Commun. Math. Phys. \textbf{142} (1991), 139--168

\bibitem{Gru:PZE}
\textsc{G.~Grubb}: \emph{Poles of zeta and eta functions for perturbations of
  the Atiyah-Patodi-Singer problem}.
\newblock University of Copenhagen, Preprint No. 14 (1999)

\bibitem{Gru:TEP}
\textsc{G.~Grubb}: \emph{Trace expansions for pseudodifferential boundary
  problems for Dirac--type operators and more general systems}.
\newblock Arkiv f. Matematik \textbf{37} (1999), 45--86

\bibitem{GruSee:WPP}
\textsc{G.~Grubb} \protect\BIBand{} \textsc{R.~Seeley}: \emph{Weakly parametric
  pseudodifferential operators and Atiyah--Patodi--Singer boundary problems}.
\newblock Invent. Math. \textbf{121} (1995), 481--529

\bibitem{GruSee:ZEF}
\textsc{G.~Grubb} \protect\BIBand{} \textsc{R.~T. Seeley}: \emph{Zeta and eta
  functions for Atiyah--Patodi--Singer operators}.
\newblock J. Geom. Anal. \textbf{6} (1996), 31--77

\bibitem{HMM:SFM}
\textsc{A.~Hassell,  R.~Mazzeo}, \protect\BIBand{} \textsc{R.~Melrose}:
\emph{A
signature formula for manifolds with corners of codimension two}.
\newblock Topology
\textbf{36} (1997), no. 5, 1055--1075.

\bibitem{Kat:PTL}
\textsc{T.~Kato}: \emph{Perturbation theory for linear operators}, vol. 132 of
  \emph{Grundlehren der Mathematischen Wissenschaften}.
\newblock 2nd edn. Springer--Verlag, Berlin--Heidelberg--New York (1976)

\bibitem{Kui:HTU}
\textsc{N.~H. Kuiper}: \emph{The homotopy type of the unitary group of Hilbert
  space}.
\newblock Topology \textbf{3} (1965), 19--30

\bibitem{Les:IP}
\textsc{M.~Lesch}: \emph{Determinants of Dirac type operators on manifolds with
  boundary}.
\newblock In preparation

\bibitem{LesWoj:IGA}
\textsc{M.~Lesch} \protect\BIBand{} \textsc{K.~P. Wojciechowski}: \emph{On the
  $\eta$--invariant of generalized Atiyah--Patodi--Singer boundary value
  problems}.
\newblock Ill. J. Math. \textbf{40} (1996), 30--46

\bibitem{Mul:IDO}
\textsc{W.~M\"uller}: \emph{On the $L^2$--index of Dirac operators on manifolds
  with corners of codimension two. I.}
\newblock J. Diff. Geom. \textbf{44} (1996), 97--177

\bibitem{Nic:MIS}
\textsc{L.~Nicolaescu}: \emph{The Maslov index, the spectral flow, and
  decompositions of manifolds}.
\newblock Duke Math. J. \textbf{80} (1995), 485--533

\bibitem{Pal:SAS}\textsc{R.~S. Palais}: \emph{Seminar on the Atiyah-Singer
Index Theorem.} \newblock Annals of Mathematics Studies, No. 57 Princeton
University Press, Princeton, N.J. 1965 x+366 pp.

\bibitem{Phi:SAF}
\textsc{J.~Phillips}: \emph{Self-adjoint Fredholm operators and spectral flow}.
\newblock Canad. Math. Bull. \textbf{39} (1996), 460--467

\bibitem{Sco:DDB}
\textsc{S.~G. Scott}: \emph{Determinants of Dirac boundary value problems over
  odd-dimensional manifolds}.
\newblock Commun. Math. Phys. \textbf{173} (1995), 43--76

\bibitem{ScoWoj:DQD}
\textsc{S.~G. Scott} \protect\BIBand{} \textsc{K.~P. Wojciechowski}: \emph{The
  $\zeta$--determinant and Quillen determinant for a Dirac operator on a
  manifold with boundary}.
\newblock Preprint, IUPUI (1999)

\bibitem{See:TPO}
\textsc{R.~T. Seeley}: \emph{Topics in pseudo--differential operators}.
\newblock C.I.M.E., Conference on pseudo--differential operators 1968, Edizioni
  Cremonese, Roma, 1969, pp. 169--305

\bibitem{Wal:NAS}
\textsc{C.~T.~C. Wall}: \emph{Non-additivity of the signature}.
\newblock Invent. Math. \textbf{7} (1969), 269--274

\bibitem{Woj:AEII}
\textsc{K.~P. Wojciechowski}: \emph{The additivity of the $\eta$--invariant.
  The case of an invertible tangential operator}.
\newblock Houston J. Math. \textbf{20} (1994), 603--621

\bibitem{Woj:AEIII}
\textsc{K.~P. Wojciechowski}: \emph{The additivity of the $\eta$--invariant.
  The case of a singular tangential operator}.
\newblock Commun. Math. Phys. \textbf{109} (1995), 315--327

\bibitem{Woj:DAI}
\textsc{K.~P. Wojciechowski}: \emph{The $\zeta$-determinant and the additivity
  of the $\eta$-invariant on the smooth, self-adjoint Grassmannian}.
\newblock Commun. Math. Phys. \textbf{201} (1999), 423--444

\end{thebibliography}
\bibliographystyle{lesch}
\end{document}